\documentclass[11pt,a4paper]{article}
\usepackage[english]{babel}
\usepackage{latexsym}
\usepackage{amssymb,amsbsy,amsmath,amsfonts,amssymb,amscd}
\usepackage{amsmath, amsthm}
\usepackage[utf8]{inputenc}
\usepackage{epsfig, graphicx}
\usepackage{color}
\usepackage{fullpage}
\usepackage{graphics}
\usepackage{wrapfig}
\usepackage{enumitem}
\usepackage{mathrsfs}
\usepackage{environ}
\usepackage{lipsum}
\usepackage{enumitem} 
\usepackage{sectsty}
\usepackage{tcolorbox}
\usepackage{dsfont}
\usepackage[colorlinks=true,linkcolor=red,citecolor=blue]{hyperref}
\usepackage{mathtools,accents}

\def\X{\mathrm{X}}
\def\V{\mathrm{V}}

\def\R{\mathbb R}

\def\N{\mathbb N}

\def\T{\mathbb T}

\def\det{\mathrm{det}}

\def\Ld{\mathrm{L}}
\def\H{\mathrm{H}}
\def\M{\mathrm{M}}
\def\B{\mathrm{B}}
\def\E{\mathrm{E}}
\def\W{\mathrm{W}}
\newcommand{\enstq}[2]{\left\{#1~\middle|~#2\right\}}

\newtcolorbox{dev}{arc=0pt}
\newcounter{compteur}
\counterwithin{compteur}{section}

\definecolor{thmcolor}{rgb}{0.8,0.14,0.2}
\definecolor{defcolor}{rgb}{0.0,0.50,0.0}
\definecolor{excolor}{rgb}{0.50,0.0,0.990}
\definecolor{applicolor}{rgb}{0.50,0.0,0.990}

\newtheoremstyle{thm}
  {\topsep}
  {\topsep}
  {\itshape}
  {0pt}
  {\bfseries}
  {.}
  { }
  {\textcolor{black!100}{\thmname{#1}\thmnumber{ #2}}\thmnote{ (#3)}}

\newtheoremstyle{def}
  {3pt}
  {3pt}
  {}
  {0pt}
  {\bfseries}
  {.}
  { }
  {\textcolor{black!100}{\thmname{#1}\thmnumber{ #2}}\thmnote{ (#3)}}
  
\newtheoremstyle{ex}
  {1pt}
  {1pt}
  {}
  {0pt}
  {\bfseries}
  {.}
  { }
  {\textcolor{black!100}{\thmname{#1}\thmnumber{ #2}\thmnote{ (#3)}}}
  
\theoremstyle{thm}
\newtheorem{theoreme}[compteur]{Theorem}
\newtheorem{corollaire}[compteur]{Corollary}
\newtheorem{proposition}[compteur]{Proposition}
\theoremstyle{appli}

\newtheorem{lemme}[compteur]{Lemma}
\theoremstyle{def}
\newtheorem{definition}[compteur]{Definition}
\newtheorem{notation}[compteur]{Notation}
\theoremstyle{ex}

\newtheorem{remarque}[compteur]{Remark}

\numberwithin{equation}{section}

\subsectionfont{\normalfont\Large\bfseries}
\sectionfont{\normalfont\LARGE\bfseries}
\newcommand{\nocontentsline}[3]{}
\newcommand{\tocless}[2]{\bgroup\let\addcontentsline=\nocontentsline#1{#2}\egroup}

\title{Concentration versus absorption for the Vlasov-Navier-Stokes system on bounded domains}

\vspace{2mm}

\author{Lucas Ertzbischoff\footnote{Centre de Math\'ematiques Laurent Schwartz (UMR 7640), Ecole Polytechnique, Institut Polytechnique de Paris, 91128 Palaiseau Cedex, France (\href{mailto:lucas.ertzbischoff@polytechnique.edu}{lucas.ertzbischoff@polytechnique.edu})}, \, Daniel Han-Kwan\footnote{Centre de Math\'ematiques Laurent Schwartz (UMR 7640), Ecole Polytechnique, Institut Polytechnique de Paris, 91128 Palaiseau Cedex, France (\href{mailto:daniel.han-kwan@polytechnique.edu}{daniel.han-kwan@polytechnique.edu})} \, and Ayman Moussa\footnote{Sorbonne Université, Université Paris-Diderot, CNRS, INRIA, Laboratoire Jacques-Louis Lions (LJLL), F-75005 Paris, France  (\href{mailto:ayman.moussa@sorbonne-universite.fr}{ayman.moussa@sorbonne-universite.fr})}}

\begin{document}
\date{}
\maketitle

\begin{abstract}
We study the large time behavior of small data solutions to the Vlasov-Navier-Stokes system set on $\Omega \times \R^3$, for a smooth bounded domain $\Omega$ of $\R^3$, with homogeneous Dirichlet boundary condition for the fluid and absorption boundary condition for the kinetic phase. We prove that the fluid velocity homogenizes to $0$ while the distribution function concentrates towards a Dirac mass in velocity centered at $0$, with an exponential rate. The proof, which follows the methods introduced in \cite{HKMM}, requires a careful analysis of the boundary effects. We also exhibit examples of classes of initial data leading to a variety of asymptotic behaviors for the kinetic density, from total absorption to no absorption at all.
\end{abstract}

%

\section{Introduction}\label{Section:intro}
Let $\Omega$ be a smooth connected and bounded open set of $\R^3$. We consider the Vlasov-Navier-Stokes system in $\Omega \times \R^3$
\begin{align}
\partial_t u + (u \cdot\nabla)u-\Delta u + \nabla p &=j_f-\rho_f u,  \ \ (t,x) \in \R^{*}_{+} \times \Omega, \label{eq:NS} \\ 
{\rm div} \, u &=0, \ \ (t,x) \in \R^{*}_{+} \times \Omega, \label{eq:NS2}\\
\partial_t f +v\cdot \nabla_x f + {\rm div}_v((u-v)f)&=0, \ \ (t,x,v) \in \R^{*}_{+} \times \Omega \times \R^3 \label{eq:Vlasov},
\end{align}
where we define
\begin{align*}
\rho_f(t,x)&:=\int_{\R^3} f(t,x,v) \,\mathrm{d}v,\\ 
 j_f(t,x)&:=\int_{\R^3} vf(t,x,v)\,\mathrm{d}v.
\end{align*}

This system of nonlinear PDEs aims at describing the transport of small particles (the \textit{dispersed phase}) immerged within a Newtonian viscous and incompressible fluid (the \textit{continous phase}). Here, we describe the particles thanks to a kinetic distribution function $f(t,\cdot,\cdot)$ on the phase space $\Omega \times \R^3$ while the fluid is described thanks to its velocity $u(t,\cdot)$ and pressure $p(t,\cdot)$. The function $f$ satisfies absorption boundary condition and we impose no-slip boundary condition for the vector-field $u$. We refer to Subsection \ref{subsec:notations} for more details.

In the large gallery of such \emph{fluid-kinetic systems}, which stem from the works of O'Rourke \cite{oro} and Williams \cite{will}, the dispersed phase and the fluid are related through a particular coupling (see \cite{desvill}). Here, we neglect the effect of collisions, coalescence and fragmentation between particles and we work under the assumption of \emph{thin sprays}, considering that the volume occupied by the droplets is negligeable compared to that occupied by the fluid. This modelling corresponds to a moderate Reynolds number for the particles and is for instance a prototype for the description of an aerosol in the air. The coupling is thus made of a drag term in the Vlasov equation (\ref{eq:Vlasov}), and of a source term in the Navier-Stokes equations (\ref{eq:NS}), called the \emph{Brinkman force}, which describes the exchange of momentum between the particles and the fluid. In this system called the Vlasov-Navier-Stokes system, where physical constants are all normalized, both unknowns $f$ and $u$ depend on each other.

The Vlasov-Navier-Stokes system and its mathematical analysis have received a lot of attention in the past twenty years. The question of the global existence of weak solutions (or local existence of strong solutions) to the Cauchy problem is now well-understood on a large class of spatial tridimensional domains, like the flat torus $\T^3$ in \cite{BDGM,ChKw}, a fixed bounded domain in \cite{ABdM} and even a time-dependent bounded domain in \cite{BGM,BMM}.

Concerning the rigorous derivation of theses equations from "first principles", little is known about the whole system. A first approach consists at looking at fluid-solid equations and their \textit{mean-field} limits. In a stationary framework, homogenization techniques have been used in \cite{DGR,H,HMS,CH} to derive the Brinkman force in the fluid equation. The dynamical problem has only been adressed when the inertia of the particles is neglected and in some dilute regime (see \cite{Hof,Mech}). An alternative, still partly formal, program has been proposed in \cite{BDGR1,BDGR2} where the starting point is a coupling between two Boltzmann equations and where both fluid and particles are considered as dispersed phases.

The  large time behavior for global solutions to the Vlasov-Navier-Stokes system appears as one of the next important steps in the understanding of such fluid kinetic model. A general setting has been highlighted by Jabin in \cite{Jab} where an asymptotic scenario with concentration in velocity for the particles (namely, a convergence of the distribution function towards a Dirac mass for the velocity part) has been described for some class a kinetic equations (but with a different coupling between the fluid and the particles).

The asymptotic dynamics of the Vlasov-Navier-Stokes has been studied in a recent series of papers \cite{GHKM,HKMM,HK}: a concentration phenomenon in velocity, leading to monokinetic large time behavior for the particles, is proven in \cite{HKMM} in the case of the torus $\T^3$ and in \cite{HK} in the case of the whole space $\R^3$. This asymptotics is obtained for data that are in some sense close to equilibrium
(we mention that a first contribution of Choi and Kwon has been made in the same direction on the torus in \cite{ChKw}, but with an \textit{a priori} assumption on the solutions which is not made in the two previous articles). 
We refer to Section \ref{tore-space-rec} of this article where we give more details about the strategy which is used. 
When a Fokker-Planck dissipation term (namely, a term of the form $- \Delta_v f$ ) is added in the kinetic equation (\ref{eq:Vlasov}), global classical solutions can be constructed for data close to Maxwellian equilibria and this non-singular steady states locally attract these solutions (see \cite{GLMZ}).

Fore more general and physical domains, very few articles deal with the question of the large time behavior of the Vlasov-Navier-Stokes system. A seminal article of Hamdache \cite{Ham} deals with the Vlasov-Stokes system with specular reflexion on a bounded domain of $\R^3$ and shows that, up to a subsequence, the sequence of solutions $(u_n(t),f_n(t)):=(u(t+t_n),f(t+t_n))$ converges in some sense to a solution of the same system as $t_n \overset{n \rightarrow + \infty}{\longrightarrow} + \infty$. 
In \cite{GHKM}, the authors consider the particular case of the Vlasov-Navier-Stokes system on a rectangle in $\R^2$, where the particular geometry of the domain allows to construct weak solutions around non-singular stationary equilibrium and where a geometric control condition helps to avoid concentration in velocity scenario for the particles (we also refer to Section \ref{tore-space-rec} below for more details).

\medskip

\textbf{Main contribution of this paper.} \ We study the Vlasov-Navier-Stokes system on a bounded domain of $\R^3$ with absorption boundary condition for the distribution function and homogeneous Dirichlet boundary condition for the fluid velocity. We aim at giving a proof of the large time monokinetic behavior to solutions to the system. This singular behavior is adressed for global weak solutions satisfying a natural energy-dissipation inequality and starting at initial data close to equilibrium. Furthermore, there is a competition between concentration and absorption which determine the final dynamics. A broad range of outcomes is possible for the asymptotic spatial profile as we illustrate by exhibiting examples of initial data, leading to a variety of behaviors, from total absorption of the particles to no absorption at all.

\subsection{Notations and definitions}\label{subsec:notations}
We denote the phase-space domain by
\begin{align*}
\mathcal{O}:=\Omega \times \R^3.
\end{align*}
In the following, the outer-pointing normal to the boundary at a point $x \in \partial \Omega$ will be denoted by $n(x)$. We denote by $\mathscr{D}_{\mathrm{div}}(\Omega)$ the set of smooth $\R^3$ valued divergence free vector-fields having compact support in $\Omega$. The closures of $\mathscr{D}_{\mathrm{div}}(\Omega)$ in $\Ld^2(\Omega)$ and in $\H^1(\Omega)$ are respectively denoted by $\Ld^2_{\mathrm{div}}(\Omega)$ and by $\H^1_{\mathrm{div}}(\Omega)$. We write $\H^{-1}_{\mathrm{div}}(\Omega)$ for the dual of the later. 

\medskip

We first define the class of admissible initial data for (\ref{eq:NS})-(\ref{eq:Vlasov}).
\begin{definition}[Initial condition]\label{CIadmissible}
We shall say that a couple $(u_0,f_0)$ is an admissible initial condition if 
\begin{align*}
u_0 &\in \Ld^2_{\mathrm{div}}(\Omega), \\
f_0 &\in \Ld^1 \cap \Ld^{\infty}(\Omega \times \R^3),\\
f_0 &\geq 0, \ \ \int_{\Omega \times \R^3} f_0 (x,v)\, \mathrm{d}x \, \mathrm{d}v =1,\\
(x,v)& \mapsto f_0 (x,v)\vert v \vert^2 \in \Ld^1(\Omega \times \R^3). 
\end{align*}
\end{definition}

The system (\ref{eq:NS})-(\ref{eq:Vlasov}) is supplemented with the following initial conditions for the fluid velocity $u$ and the distribution function $f$
\begin{align*}
u_{\mid t=0}&=u_0 \text{ in } \Omega,\\
f_{\mid t=0}&=f_0 \text{ in } \mathcal{O}.
\end{align*}

We prescribe the following homogeneous Dirichlet boundary condition for the fluid velocity
\begin{align}\label{bcond-fluid}
u(t,\cdot)=0, \text{ on } \partial \Omega.
\end{align}
We also need to introduce the following outgoing/incoming phase-space boundary for the dispersed phase:
\begin{align*}
\Sigma^{\pm}&:= \left\lbrace  (x,v) \in \partial \Omega \times \R^3   \mid \pm v \cdot n(x)>0 \right\rbrace,\\
\Sigma_0&:= \left\lbrace  (x,v) \in \partial \Omega \times \R^3   \mid v \cdot n(x)=0 \right\rbrace,\\
\Sigma &:= \Sigma^+ \sqcup \Sigma^- \sqcup \Sigma_0= \partial \Omega \times \R^3.
\end{align*}
Then, we prescribe the following absorption boundary condition for the distribution function $f$:
\begin{align}\label{bcond-f}
f(t,\cdot,\cdot)=0, \text{ on } \Sigma^{-},
\end{align}
meaning that all the particles reaching transversally the physical boundary are deposited.
\medskip

We then define the energy and the dissipation of the whole system.

\begin{definition}
\begin{enumerate}

\item The \textbf{kinetic energy} of the Vlasov-Navier-Stokes system is defined for all $t \geq 0$ as 
\begin{align}\label{eq:Energy}
\mathrm{E}(t):= \dfrac{1}{2}\int_{\Omega} \vert u(t,x)\vert^2  \, \mathrm{d} x
+\dfrac{1}{2}\int_{\Omega \times \R^3} f(t,x,v) \vert v \vert^2 \, \mathrm{d} x \, \mathrm{d} v.
\end{align}  

\item The \textbf{dissipation} of the Vlasov-Navier-Stokes system is defined for all $t \geq 0$ as
\begin{align}\label{eq:Dissipation}
\mathrm{D}(t):= \int_{\Omega} \vert \nabla u(t,x)\vert^2 \, \mathrm{d} x+ \int_{\Omega \times \R^3} f(t,x,v) \vert u(t,x)-v\vert^2 \, \mathrm{d} x \, \mathrm{d} v. 
\end{align}  
\end{enumerate}
\end{definition}

These two functionals naturally appear when looking for \textit{a priori} estimates satisfied by solutions to the Vlasov-Navier-Stokes system. One can indeed check that the following energy-dissipation identity formally holds
\begin{align}\label{introEnergy}
\dfrac{\mathrm{d}}{\mathrm{d}t}\mathrm{E}(t) + \mathrm{D}(t)=0.
\end{align}

We then introduce the notion of weak solution to the system.

\begin{definition}[Weak solution]\label{weak:sol}
Consider an admissible initial condition $(u_0,f_0)$  in the sense of Definition \ref{CIadmissible}. A global weak solution to the Vlasov-Navier-Stokes system with initial condition $(u_0,f_0)$ on $\Omega$ is a pair $(u,f)$ such that
\begin{align*}
u &\in \Ld^{\infty}_{\mathrm{loc}}(\R^+;\Ld^2_{\mathrm{div}}(\Omega))\cap  \Ld^{2}_{\mathrm{loc}}(\R^+;\H^1_{\mathrm{div}}(\Omega)),\\
f &\in \Ld^{\infty}_{\mathrm{loc}}(\R^+;\Ld^1 \cap \Ld^{\infty}(\Omega \times \R^d)),\\
j_f-\rho_f u &\in \Ld^{2}_{\mathrm{loc}}(\R^+;\H^{-1}_{\mathrm{div}}(\Omega)),\\
f(t,x,v) & \geq 0 \ \ \text{for almost all} \ (t,x,v) \in \R^+ \times \Omega \times \R^3,
\end{align*}
with $u$ being a Leray solution to the Navier-Stokes equations (\ref{eq:NS})-(\ref{eq:NS2}) with strong energy inequality (with initial condition $u_0$) and $f$ being a renormalized solution in the sense of DiPerna-Lions (see Appendix \ref{DiPernaLions}) to the Vlasov equation (\ref{eq:Vlasov}) (with initial condition $f_0$).
Furthermore, we require that the following energy estimate holds for almost all $s \geq 0$ (including $s=0$) and all $t \geq s$ 
\begin{align}\label{ineq-energy}
\mathrm{E}(t) + \int_s ^t \mathrm{D}(\sigma)\mathrm{d}\sigma \leq \mathrm{E}(s).
\end{align}
\end{definition}

\begin{remarque}
A weak solution $(u,f)$ to the Vlasov-Navier-Stokes system satisfies the following weak formulations.

\textbullet \ \ 
For all $T>0$, for all $\phi \in \mathscr{D}([0,T] \times \Omega)$ such that $\phi(T)=0$ and ${\rm div}_x \,  \phi=0$ 
\begin{multline*}\int_0^T\int_{\Omega}\left[ u \cdot \partial_t \phi +(u \otimes u) : \nabla_x \phi -\nabla_x u :\nabla_x \phi \right](t,x) \, \mathrm{d}x \, \mathrm{d}t\\ = - \int_0^T \langle j_f-\rho_f u, \phi \rangle(t) \, \mathrm{d}t  - \int_{\Omega} u_0(x) \cdot \phi(0,x)\, \mathrm{d}x.
\end{multline*}

\textbullet \ \ 
For all $T>0$, for all $\psi \in \mathscr{D}([0,T] \times\overline{ \Omega} \times \R^3)$ such that $\psi(T)=0$ and vanishing on $\R^+ \times (\Sigma^+ \sqcup \Sigma^0)$
\begin{multline*}
 \int_0^T \int  \! \! \! \! \int_{\Omega \times \R^3} f \left[ \partial_t \psi + v \cdot \nabla_x \psi +(u-v) \cdot \nabla_v \psi \right](t,x,v) \, \mathrm{d}x \, \mathrm{d}v \, \mathrm{d}t \\ =-\int  \! \! \! \!  \int_{\Omega \times \R^3} f_0(x,v) \psi(0,x,v) \, \mathrm{d}x \, \mathrm{d}v.
\end{multline*}
Furthermore, for such a weak solution to the Vlasov equation, we can define a trace on the phase space boundary $\Sigma$ in the DiPerna-Lions framework for transport equations $\Omega \times \R^3$: we refer to Section \ref{DiPernaLions} of the Appendix for further properties of the solution to this initial boundary value problem (see \cite[Section 1 - Chap 6]{BF} and \cite{Mi}).
\end{remarque}

\medskip

Finally, we introduce the following definitions that will be useful later.
\begin{definition}\label{notation:moments}
For any $\alpha >0$, and any measurable function $f:\R^+ \times \Omega \times \R^3 \rightarrow \R^+$, we set
\begin{align*}
m_{\alpha}f(t,x)&:=\int_{\R^3} \vert v \vert^{\alpha} f(t,x,v)\,\mathrm{d}v, \\
M_{\alpha}f(t,x)&:=\int_{\Omega \times \R^3} \vert v \vert^{\alpha} f(t,x,v)\,\mathrm{d}x \, \mathrm{d}v=\int_{\Omega} m_{\alpha} f(t,x) \,\mathrm{d}x.
\end{align*} 
\end{definition}

\begin{definition}\label{def:pointdecay}
We say that an initial kinetic condition $f_0$ satisfies the pointwise decay assumption of order $q>0$ if 
\begin{align}\label{def:Nq}
N_q(f_0):= \underset{\substack{x \in \Omega \\ v \in \R^3}}{\sup} (1+\vert v \vert^q) f^0(x,v) < \infty.
\end{align}
\end{definition}

\subsection{The torus, the whole space and the rectangle}\label{tore-space-rec}
As already said, the large-time behavior of the Vlasov-Navier-Stokes has been tackled in the case of the torus $\T^3$ in \cite{HKMM} or the whole space $\R^3$ in \cite{HK}, and for which concentration in velocity happens in some small data regime (that is to say, a convergence towards a Dirac mass in velocity). We refer to Section 1.2 of the introduction of \cite{HKMM} for heuristics about this monokinetic asymptotic phenomenon when there are no boundaries, relying on an explicit formula for the solutions to the linearized equations around states of the form $(\overline{U},\overline{f}=0)$ with $\overline{U} \in \R^3$.

On the contrary, non-singular equilibra have been constructed in the case of a bidimensional rectangle in \cite{GHKM}. Taking advantage of the specific geometry of the domain and of absorbing boundary conditions for the particles, it has been proven that such stationary solutions are locally asymptotically stable relatively to compact perturbations

Let us describe more specifically the main strategy that has been employed in each of these situations.

\textbullet \ \ On the torus, Choi and Kwon have introduced in \cite{ChKw} a version of the following so-called \emph{modulated energy} 
\begin{multline}\label{mod-energy-torus}
\mathcal{E}_{\T^3}(t):= 
\dfrac{1}{2}\int_{\T^3 \times \R^3} f(t,x,v) \vert v-\langle j_f(t) \rangle\vert^2 \, \mathrm{d} v \, \mathrm{d} x \\
+  \dfrac{1}{2}\int_{\T^3} \vert u(t,x)- \langle  u(t) \rangle \vert^2 \, \mathrm{d} x
+ \dfrac{1}{4} \vert \langle j_f(t) \rangle - \langle u(t) \rangle \vert ^2, 
\end{multline}
where $\langle \cdot \rangle$ stands for the average on $\T^3$. They have used this functional to describe the large-time dynamics of the system but under an \textit{a priori} assumption on the solutions. In \cite{HKMM}, the authors have worked  in a small data regime and have provided the first complete description of the asymptotic behavior of the system. Loosely speaking, under a condition of the type
\begin{align}\label{Intro:datasmall}
\mathcal{E}_{\T^3}(0) + \Vert u_0 \Vert_{\dot\H^{1/2}(\T^3)} \ll 1,
\end{align} 
the authors have showed that the fluid velocity $u(t)$ homogenizes when $t \rightarrow + \infty$ to the constant velocity $U_0 :=\langle u_0 + j_{f_0} \rangle/2$, while the kinetic distribution function $f(t)$ converges in velocity to the Dirac mass supported at $U_0$. Moreover, this convergence is exponentially fast and can be measured thanks to the 1-Wasserstein distance on $\T^3 \times \R^3$.

The previous modulated energy is at the heart of the proof of this monokinetic large time behavior for the system on the torus and is linked to the dissipation thanks to the formal identity 
\begin{align*}
\dfrac{\mathrm{d}}{\mathrm{d}t}\mathcal{E}_{\T^3}(t) + \mathrm{D}(t)=0.
\end{align*}

As a matter of fact, this modulated energy essentially captures concentration phenomena so that controlling such a quantity is the main key to understand the large-time dynamics of the system. In short, under the assumption that $\rho_f \in \Ld^{\infty}(\R^+;\Ld^{3/2}(\T^3))$, Choi and Kwon proved in \cite{ChKw} that
\begin{align}\label{decayT3intro}
\forall t \geq 0, \ \ \ \mathcal{E}_{\T^3}(t)\lesssim e^{- \lambda t} \mathcal{E}_{\T^3}(0),
\end{align}
for some $\lambda>0$. Thanks to this exponential decay, they deduced that the asymptotics we have mentioned above hold for the fluid velocity and for the kinetic distribution.

The main strategy of \cite{HKMM} is based on a bootstrap analysis whose aim is to ensure that $\Vert \rho_f \Vert_{\Ld^{\infty}(\R^+;\Ld^{\infty}(\T^3))} < \infty$. Thanks to a straigthening change of variable in velocity, it is shown in \cite{HKMM} that this condition is actually implied by an estimate of the type 
\begin{align}\label{eq:smallnauT3intro}
\int_0^{\infty} \Vert \nabla u(s) \Vert_{\Ld^{\infty}(\T^3)} \, \mathrm{d}s \ll 1.
\end{align}
Then, the framework of the Fujita-Kato solutions to the Navier-Stokes system has been leveraged in \cite{HKMM} to ensure that such a control holds under the condition (\ref{Intro:datasmall}).

\textbullet \ \ In the case of the euclidean space $\R^3$, the large-time behavior of the system has been investigated in \cite{HK} where it has been shown that concentration in velocity towards a Dirac mass supported at 0 occurs also for the kinetic distribution, while the fluid velocity converges to 0. Because of the unbounded nature of the spatial domain, this convergence is only at a polynomial rate. More precisely, the good functional to look at to understand the concentration phenomenon turns out to be the kinetic energy $\mathrm{E}$. One of the the main results of \cite{HK} is a conditional decay of the form
\begin{align*}
\forall t \geq 0, \ \ \ \mathrm{E}(t) \leq \dfrac{\varphi_{\alpha}(\mathrm{E(0)})}{(1+t)^{\alpha}}, \ \ \text{for all} \ \ \alpha \in ]0,3/2[,
\end{align*}
for some function $\varphi_{\alpha}$, up to an \textit{a priori} control on the moment $\rho_f$. As in the torus case, a bootstrap analysis is used to obtain such a control for small data solutions. But in this unbounded context where the decay of the energy is only polynomial, the Brinkman force requires a more careful treatment which has led to the derivation of a new family of identities for higher order dissipation functionals.

\textbullet \ \ In \cite{GHKM}, the authors have dealt with the particular case of a bidimensional rectangle $\Omega:=(-L,L) \times (-1,1)$. Here, the boundary condition for the fluid part is a Dirichlet boundary condition matching a Poiseuille flow $u_p$ (which is a stationary solution to the Navier-Stokes equations). For the particles, partly absorbing boundary conditions are used with absorption boundary conditions on the horizontal parts and an injection boundary condition on the vertical left part.
Thanks to a geometric control condition (referred to as the \textit{exit geometric condition}),
compelling the particles to be absorbed by the boundary before a fixed finite time, one can construct non-trivial smooth equilibria $(\overline{u},\overline{f})$ for the system, so that the concentration in velocity scenario does not occur. Furthermore, if $(\overline{u},\overline{f})$ is a smooth stationary solution close to $(u_p,0)$ (for some small injection term and a small Poiseuille flow $u_p$), any $(u_0,f_0)$ which is a small perturbation of $(\overline{u},\overline{f})$ gives birth to a weak solution $(u,f)$ to the system satisfying
\begin{align*}
 \forall t \in \R^+, \ \ \Vert f(t)-\overline{f} \Vert_{\Ld^{2}_{x,v}} +  \Vert u(t)-\overline{u} \Vert_{\Ld^{2}_x} \lesssim e^{-t}.
\end{align*}
In short, the previous equilibria are asymptotically stable in $\Ld^2$, under small localized perturbations. 

%
%
%

\medskip
In the following, we will consider solutions starting at initial data close to the equilibrium $(0,0)$ and prove the existence of an asymptotic profile $\rho^{\infty} \in \Ld^{\infty}(\Omega)$ such that the following weak convergence holds 
\begin{align*}
f(t) \xrightharpoonup[t \to +\infty]{} \rho^{\infty} \otimes \delta_{v=0},
\end{align*}
where the tensor product is in $(x,v)$. As said before, this kind of singular limit was already present in the case of the torus and of the whole space while geometric control was preventing such monokinetic behavior in the case of the rectangle. Once this monokinetic large time behavior is established, we will be interested in a possible further study of the spatial profile $\rho^{\infty}$: because of the boundary effects, it may indeed exist some particle trajectories which are leaving the domain during the evolution. Therefore, we would like to study different possible scenarios, where absorption of the particles can prevail or not, in order to describe the asymptotic local density.

\subsection{Main results}
The main result of this article is stated in the following Theorem and Corollary.
\begin{theoreme}\label{theoreme1}
There exists a universal constant $\varepsilon_0>0$ and a nondecreasing function $\varphi : \R^+ \rightarrow \R^+$ such that the following holds. Let $(u_0,f_0)$ be an admissible initial condition in the sense of Definition \ref{CIadmissible} satisfying 
\begin{align}\label{data:hyp}
\begin{split}
& u_0 \in \H^1_{\mathrm{div}}(\Omega), \\
& M_6 f_0 + N_q f_0<\infty, \ \  \text{for some } q>3.
\end{split}
\end{align}
If
\begin{align}\label{smallness:condition}
\varphi \Big(1+ N_q f_0 \Big)\E(0)<\varepsilon_0, \ \ \Vert \nabla u_0 \Vert_{\Ld^{2}(\Omega)}<\varepsilon_0,
\end{align}
where $\E$ is defined in (\ref{eq:Energy}), then for any weak solution $(u,f)$ to the Vlasov-Navier-Stokes system with initial data $(u_0,f_0)$, there exist constants $\lambda, C_{\lambda} >0$ such that
\begin{align}
& \rho_f \in \Ld^{\infty}(\R^+;\Ld^{\infty}(\Omega)), \\[2mm]
& \E(t) \leq \E(0) C_{\lambda} \exp(- \lambda t), \ \ t \geq 0.
\end{align}
\end{theoreme}

\begin{remarque}
In the case of the torus \cite{HKMM}, the required assumption on the initial fluid velocity $u_0$ was $u_0 \in \dot\H^{1/2}(\T^3)$ and allowed to rely on some parabolic smoothing for the solution $u$. Here, we have preferred to state the result with the assumption (\ref{data:hyp}) in order to avoid unnecessary technical developments.
\end{remarque}


\begin{corollaire}\label{coro:result}
Under the same assumptions of Theorem (\ref{theoreme1}), for any weak solution $(u,f)$ to the Vlasov-Navier-Stokes system with admissible initial  data $(u_0,f_0)$, there exist constants $\lambda, C_{\lambda} >0$ and $\rho^{\infty} \in \Ld^{\infty}(\Omega)$ such that for all $t \geq 0$
\begin{align}
&\mathrm{W}_1\Big( f(t), \rho_f(t)\otimes \delta_{v=0} \Big) + \left\Vert u(t) \right\Vert_{\Ld^{2}(\Omega)} \leq \E(0)^{1/2} C_{\lambda} \exp(- \lambda t), \\[2mm]
&\rho_f(t) \overset{t \rightarrow +\infty}{\longrightarrow} \rho^{\infty} \ \ in \ \H^{-1}(\Omega),
\end{align}
where $\mathrm{W}_1$ stands for the Wasserstein distance on $\overline{\Omega}\times \R^3$. Moreover, the last convergence also occurs with an exponential rate.
\end{corollaire}

As in \cite{HKMM}, it is possible to provide a further description of the limit profile $\rho^{\infty}$, using a notion of asymptotic characteristics.
\begin{theoreme}
\label{thm-infini}
For $\delta$ small enough, under the assumptions of Theorem \ref{theoreme1}, and if
\begin{align}\label{ineq:asympnab}
\int_0^{+\infty} \|\nabla u(\tau)\|_{\Ld^{\infty}(\Omega)} \, \mathrm{d} \tau\leq \delta,
\end{align}
then there exists a vector field
\begin{align*}
  \R^3\times\R^3 &\longrightarrow \R^3\\
  (x,v) &\longmapsto \X_{\infty}(x,v),
\end{align*}
belonging to $\mathscr{C}^1(\R^3\times\R^3)$ and such that the following holds. For all $v \in \R^3$, the mapping $\X_{\infty,v} : x \mapsto \X_{\infty}(x,v)$ is a  $\mathscr{C}^1$-diffeomorphism from $\R^3$ to itself and we have for almost every $ x \in \Omega$ 
\begin{align}
\label{eq-rhoinfini}
\rho^{\infty}(x)=\int_{\R^3}\mathbf{1}_{\mathcal{U}^{\infty}}(x,v)f_0\left(\X^{-1}_{\infty,v}(x),v \right) \vert \det \, \mathrm{D}_x \X^{-1}_{\infty,v}(x) \vert \, \mathrm{d}v,
\end{align} 
where the set $\mathcal{U}^{\infty}$ is defined as follows: for $(x,v) \in \Omega \times \R^3$,  
\begin{align*}
(x,v) \in \mathcal{U}^{\infty}
\Longleftrightarrow \, & \exists ! y \in \Omega , \ \ x=\X_{\infty}(y,v) \ \ \text{and} \ \  \forall t \geq 0, \ \  \X^t_0(y,v) \in \Omega,
\end{align*}
where
\begin{align}
\mathrm{X}^t_0(y,v)&=x+(1-e^{-t})v+\int_0^t (1-e^{\tau-t}) u(\tau,\mathrm{X}^{\tau}_0(x,v))  \, \mathrm{d}\tau,\label{eq:X} \\ 
\X_{\infty}(y,v) &=x+v+\int_0^{\infty} u(\tau,\mathrm{X}^{\tau}_{0}(y,v))  \, \mathrm{d}\tau.\label{eq:Xinfini}
\end{align}
\end{theoreme}

\begin{remarque}
In the previous statement, we can actually get rid of the assumption (\ref{ineq:asympnab}) by only considering the evolution of the system from time $t=1$. Indeed, the proof of Theorem \ref{theoreme1} and Corollary \ref{coro:result} will ensure (see the bootstrap procedure in Section \ref{section:bootstrap}) that for all $\epsilon>0$, the quantity $\Vert \nabla u \Vert_{\Ld^1(\varepsilon,+\infty;\Ld^{\infty}(\Omega))}$ can be ensured as small as required up to imposing a relevant smallness assumption (\ref{smallness:condition}). Nevertheless, for the sake of clarity, we have decided to state the result from time $t=0$. If one wants to get the result without further assumption near the time $0$, one has to replace $f_0$ by $f_{\mid t=1}$ together with integrals starting at $t=1$ and for a set $\mathcal{U}^{\infty}$ defined with the function
\begin{align*}
\mathrm{X}^t_1(y,v)&=x+(1-e^{1-t})v+\int_1^t (1-e^{\tau-t}) u(\tau,\mathrm{X}^{\tau}_1(x,v))  \, \mathrm{d}\tau.
\end{align*}
\end{remarque}

As mentioned before, an important difference with the case of the torus is that the particle trajectory may possibly escape the domain $\Omega$ because we prescribe absorption boundary conditions for the Vlasov equation. 
This means that, in Corollary \ref{coro:result} and Theorem \ref{thm-infini}, we obtain an asymptotic spatial profile $\rho^{\infty}$ whose total mass is unknown: indeed, part of the initial mass of the system may have disappeared throughout the evolution. 

In the following result, we show that any fixed mass which is less than or equal to the initial mass can be reached by the system for some well-chosen data. We recall that we consider initial distribution functions $f_0$ such that $\int_{\Omega \times \R^3} f_0(x,v) \, \mathrm{d}x \, \mathrm{d}v=1$.
\begin{proposition}\label{prop:massalpha}
Let $\alpha \in [0,1]$. There exists $(u_0,f_0)$ an admissible initial condition in the sense of Definition \ref{CIadmissible} such that
for any weak solution $(u,f)$ to the Vlasov-Navier-Stokes system starting at $(u_0,f_0)$, their exists $\rho^{\infty} \in \Ld^{\infty}(\Omega)$ which satisfies
\begin{align}
&\rho_f(t) \xrightharpoonup[t \to +\infty]{} \rho^{\infty} \ \ in \ \ \mathscr{C}(\overline{\Omega})',\\[2mm]
&\int_{\Omega} \rho^{\infty}(x) \, \mathrm{d}x=\alpha.
\end{align}
\end{proposition}

\bigskip

\textbf{Strategy and outline of the paper}. \ Let us describe the main strategy that we use in this paper and how it is organized. Our method is reminiscent of the work of \cite{HKMM}. 

In Section \ref{section:energy}, we study the conditional decay of the kinetic energy and its consequences. As explained before, such a decay will essentially be enough for the concentration in velocity to happen. Here, we strongly rely on the Poincaré inequality which holds on the domain $\Omega$ for the fluid velocity. In short, we can hope for an exponential decay of the energy provided that the local density $\rho_f$ is controlled in $\Ld^{\infty} \Ld^{\infty}$. This result is stated in Proposition \ref{decrexpo:energy}. We then show in Proposition \ref{asymp-profil} that it also provides the existence of an asymptotic spatial profile $\rho^{\infty}$.

In Section \ref{Section:particle:traj}, we investigate a way to get the desired control on this moment $\rho_f$. We mainly rely on trajectorial estimates, by introducing the characteristic curves $(\X,\V)$ for the Vlasov equation as the solution of the differential system
\begin{equation*}
\left\{
      \begin{aligned}
        \dot{\mathrm{X}}(s;t,x,v) &=\mathrm{V}(s;t,x,v),\\
\dot{\mathrm{V}} (s;t,x,v)&= u(s,\mathrm{X}(s;t,x,v))-\mathrm{V}(s;t,x,v),
      \end{aligned}
    \right.
\end{equation*}
with $(\mathrm{X}(t;t,x,v),\mathrm{V}(t;t,x,v))=(x,v)$ and where we extend $u$ by $0$ outside of $\Omega$. In Proposition \ref{Prop:formulerep}, we then derive a representation formula for weak solutions to the Vlasov equation that is valid on a bounded domain and that takes into account the boundary effects. This leads to 
\begin{align}\label{formulerepINTRO}
\rho_f(t,x)=e^{3t} \int_{\R^3} \mathbf{1}_{\tau^-(t,x,v)<0} \, f_0(\mathrm{X}(0;t,x,v),\mathrm{V}(0;t,x,v)) \, \mathrm{d}v,
\end{align}
where 
\begin{align*}
\tau^{-}(t,x,v)&:=\inf \left\lbrace s \leq t \ \mid \forall \sigma \in [s,t], \  \mathrm{X}(\sigma;t,x,v) \in \Omega \right \rbrace,
\end{align*}
stands for the first backward exit time of the particle trajectory. Thanks to the change of variable in velocity $v \mapsto \V(0;t,x,v)$, obtaining a uniform control on $\rho_f$ until time $t$ actually reduces to an inequality of the type $\Vert \nabla u \Vert_{\Ld^1(0,t;\Ld^{\infty}(\Omega))} \ll 1$. This fact is stated in Proposition \ref{Propo:gradFirst}.

A third step is made of a bootstrap analysis and aims at obtaining the previous control on $\Vert \nabla u \Vert_{\Ld^{\infty}(\Omega)}$. To do so, we interpolate this quantity by second order derivatives of $u$ together with the kinetic energy of the system, in order to get integrability in time thanks to the exponential decay of this energy. This procedure requires some integrability and higher regularity estimates for the solutions to the Navier-Stokes equations. Section \ref{Section:prepbootstrap} yields preliminaries for this subsequent bootstrap argument. In Proposition \ref{controlLinfiniLOC}, we obtain local in time controls for the velocity field $u$, as well as for the moments $\rho_f$ and $j_f$. We also provide $\H^1$ energy estimates for the velocity field $u$, using the framework of strong solutions to the Navier-Stokes equations. This one is allowed if the initial condition and the Brinkman force are small enough, as seen in Proposition \ref{Prop:H1estimates}. In order to propagate this condition, we introduce the notion of strong existence time.

In Section \ref{Section:estimateforces}, we start to estimate the second order derivatives $\mathrm{D}^2 u$ in some $\Ld^q\Ld^p$ space. To do so, we rely on maximal regularity for the Stokes system. This requires to control the Brinkman force $j_f-\rho_f u$ and the convection term $ (u \cdot \nabla) u$ in some $\Ld^a \Ld^b$ space, which is possible thanks to the higher order energy estimates satisfied by the fluid velocity.

The bootstrap analysis then takes place in Section \ref{section:bootstrap} where the smallness condition (\ref{smallness:condition}) we have used in our statements plays a key role. We essentially show that this one can ensure the global control $\Vert \nabla u \Vert_{\Ld^1(1,+\infty;\Ld^{\infty}(\Omega))} \ll 1$, ending the proof of Theorem \ref{theoreme1} and Corollary \ref{coro:result}.

In Section \ref{Section:asympprofile}, we give a further description of the asymptotic spatial density profile $\rho^{\infty}$ of the particles. In order to prove Theorem \ref{thm-infini}, we rely on the Lagrangian structure of the Vlasov equation. The main guiding line is to pass to the limit in the representation formula \eqref{formulerepINTRO} when $t\rightarrow + \infty$, by considering asymptotic characteristic curves. Here, our analysis is based upon a fine analysis of the flow $(\X,\V)$ where we have to deal again carefully with the possible exit of these trajectories from the domain $\Omega$.

As the structure of the asymptotic profile remains quite implicit, the question of its total mass is discussed in Section \ref{section:prescribe:mass}. We first study sufficient conditions on the support of the initial distribution $f_0$ and on the initial velocity $u_0$ to ensure the possible following scenarios: the total mass of the system can be preserved throughout the evolution or, on the contrary, can vanish after a certain finite time. These two opposite situations are stated in Propositions \ref{prop:inside}-\ref{prop:outside} and strongly depend on the support in velocity of $f_0$. Their analysis are strongly based on the consequences of the previous bootstrap. Combining these two examples eventually leads to Proposition \ref{prop:massalpha}.

\section{Energy dissipation: towards concentration}\label{section:energy}

We start with a Lemma which states that the total mass of the system is nonincreasing along the evolution. The proof of this result is postponed to Subsection \ref{Subsec:repform} because it requires the use of a representation formula for the distribution function $f$.

\begin{lemme}\label{loose:mass}
For any weak solution $(u,f)$ to the Vlasov-Navier-Stokes system and for all $t \geq 0$, we have
\begin{align*}
\int_{\Omega \times \R^3} f(t,x,v) \, \mathrm{d}x \, \mathrm{d}v \leq \int_{\Omega \times \R^3} f_0(x,v) \, \mathrm{d}x \, \mathrm{d}v.
\end{align*}
\end{lemme}

We can now state the following inequality which highlights the role of the kinetic energy for the study of the asymptotics.
\begin{lemme}\label{lem:W1}
For any weak solution $(u,f)$ to the Vlasov-Navier-Stokes system and for all $t \geq 0$, we have the following inequality
\begin{align}
\mathrm{W}_1\left( f(t), \rho_f(t)\otimes \delta_{v=0 } \right) + \left\Vert u(t) \right\Vert_{\Ld^{2}(\Omega)} \lesssim \E(t)^{1/2},
\end{align}
where $\mathrm{W}_1$ stands for the Wasserstein distance on $\overline{\Omega} \times \R^3$.
\end{lemme}
\begin{proof}
First, we observe that for all times $t\geq 0$, the measures $f(t)\,\mathrm{d}x \,\mathrm{d}v$ and $(\rho_f(t)\mathrm{d}x )\otimes \delta_{v=0}$ have the same mass on $\Omega \times \R^3$. Therefore we use the Monge-Kantorovich duality for $\mathrm{W}_1$ on $\overline{\Omega} \times \R^3$ (see Appendix \ref{Appendix:Wasserstein})  to get
\begin{align*}
\mathrm{W}_1\left( f(t), \rho_f(t)\otimes \delta_{v=0} \right) &=\underset{\Vert\nabla_{x,v}\phi\Vert_{\infty} \leq 1 }{\mathrm{sup}} \left\lbrace \left\vert \int_{\Omega \times \R^3} f(t,x,v)(\phi(x,v)-\phi(x,0)) \, \mathrm{d}x \, \mathrm{d}v \right\vert \right\rbrace \\
& \leq \int_{\Omega \times \R^3} f(t,x,v)\vert v \vert \, \mathrm{d}x \, \mathrm{d}v.
\end{align*}
The Cauchy-Schwarz inequality and the previous Lemma \ref{loose:mass} together with the normalisation for the density $f_0$ give us
\begin{multline*}
\mathrm{W}_1\left( f(t), \rho_f(t)\otimes \delta_{v=0} \right) \\
\leq \left(\int_{\Omega \times \R^3} f(t,x,v)\vert v \vert^2 \, \mathrm{d}x \, \mathrm{d}v\right)^{1/2} \left(\int_{\Omega \times \R^3} f(t,x,v) \, \mathrm{d}x \, \mathrm{d}v\right)^{1/2} 
\leq \sqrt{2} \E(t)^{1/2},
\end{multline*}
by definition of the kinetic energy. The inequality $ \left \Vert u(t) \right\Vert_{\Ld^{2}(\Omega)} \lesssim \E(t)^{1/2}$ is also straigthforward.
\end{proof}

We now state the following Proposition relating the dissipation and the kinetic energy.

\begin{proposition}\label{decrexpo:energy}
The exists a continuous nondecreasing function $\psi : \R^{+} \rightarrow \R^{+}$ such that the following holds. Let $(u,f)$ be a weak solution to the Vlasov-Navier-Stokes system such that $\rho_f \in \Ld^{\infty}_{\mathrm{loc}}(\R^+;\Ld^{\infty}(\Omega))$. Fix $T \in (0,+\infty]$ and set
\begin{align*}
\lambda :=\psi\left( \underset{[0,T]}{\mathrm{sup}} \ \  \Vert \rho_f (t) \Vert_{\Ld^{\infty}(\Omega)} \right),
\end{align*}  
then
\begin{align}\label{control:lDissipenergy}
\forall t \in [0,T], \ \ \lambda \E(t) \leq \mathrm{D}(t),
\end{align}
and we have the following exponential decay
\begin{align}\label{decay:expo}
\forall t \in [0,T], \ \  \E(t) \lesssim_{\lambda}\E(0) e^{- \lambda t},
\end{align}
where $\lesssim_{\lambda}$ refers to a constant only depending on $\lambda$.
\end{proposition}
\begin{proof}
The proof follows the same arguments as in \cite[Lemma 3.4]{HKMM} and was originally obtained in \cite[Theorem 1.2]{ChKw} in the case of the torus. It mainly relies on the fact that the Poincaré inequality is valid for the fluid velocity on $\Omega$.

\textbullet \ \   First, we note that (\ref{control:lDissipenergy}) implies (\ref{decay:expo}). Indeed, if (\ref{control:lDissipenergy}) holds, the energy-dissipation inequality (\ref{ineq-energy}) shows that for almost all $0 \leq s \leq t \leq T$,
\begin{align*}
\E(t) + \lambda \int_s ^t \E(\tau)  \, \mathrm{d}\tau\leq \E(s),
\end{align*}
from which we get
\begin{align*}
\lambda \int_t ^{T} \E(\tau)  \, \mathrm{d}\tau \leq \E(t), \ \ \ \  \E(t) \leq \E(s), \ \ \ \ s \leq t.
\end{align*}
We get the desired exponential decay (\ref{decay:expo}) thanks to Lemma \ref{variante:gronwall} (see Appendix \ref{Annexe:gronwall}).

\textbullet \ \ It remains to get a $\lambda >0$ such that $\lambda \E \leq \mathrm{D}$ on $[0,T]$. 
Thanks to the Poincaré inequality in $\H^1_0(\Omega)$, the exists a constant $c_P>0$ such that
\begin{align*}
\mathrm{D}(t) \geq \dfrac{1}{2}\int_{\Omega \times \R^3} f(t) \vert u(t)-v\vert^2   \, \mathrm{d} v \, \mathrm{d} x 
+  \dfrac{1}{2} c_P \Vert u(t)\Vert_{\Ld^2(\Omega)}^2.
\end{align*}
Let us denote
\begin{align*}
\widetilde{\E}(t):=\E(t)- \dfrac{1}{2} \Vert u(t)\Vert_{\Ld^2(\Omega)}^2=\dfrac{1}{2}\int_{\Omega \times \R^3} f (t,x,v)\vert v\vert^2  \, \mathrm{d} v \,  \mathrm{d} x .
\end{align*}
We see that it is sufficient to find some  $\gamma >0$ and $\beta \in (0,c_P)$ such that
\begin{align}\label{inegtechnique:energiemod}
\int_{\Omega \times \R^3} f(t) \vert u(t)-v\vert^2   \, \mathrm{d} v \, \mathrm{d} x  \geq \gamma \widetilde{\E}(t) -\beta \Vert u(t)\Vert_{\Ld^2(\Omega)}^2.
\end{align}
Indeed, if this holds, we can define $\lambda := \min(\gamma/2, c_P-\beta)>0$ and the previous inequalities yield
\begin{align*}
\mathrm{D}(t) \geq \lambda \E(t)+\dfrac{1}{2} \left( c_P-\beta-\dfrac{\lambda}{2}\right)\Vert u(t) \Vert_{\Ld^2(\Omega)}^2  \geq \lambda \E(t).
\end{align*}
In order to get the inequality (\ref{inegtechnique:energiemod}), we first write
\begin{align*}
\int_{\Omega \times \R^3} f(t) \vert u(t)-v\vert^2   \, \mathrm{d} v \, \mathrm{d} x =\int_{\Omega \times \R^3} f \vert v\vert^2  \, \mathrm{d} v \,  \mathrm{d} x 
+\int_{\Omega} \rho_f(t) \vert u(t) \vert^2   \, \mathrm{d} x  - 2  \int_{\Omega \times \R^3} f(t) \, v\cdot u(t)  \, \mathrm{d} v \, \mathrm{d} x.
\end{align*}
For the last term, we use Cauchy-Schwarz and Young inequalities to get
\begin{align*}
- 2  \int_{\Omega \times \R^3} f(t) \, v\cdot u(t)  \, \mathrm{d} v \, \mathrm{d} x \geq -\alpha \int_{\Omega \times \R^3} f(t) \vert v \vert^2   \, \mathrm{d} v \, \mathrm{d} x    - \alpha^{-1}\int_{\Omega} \rho_f(t) \vert u(t)\vert^2    \, \mathrm{d} x ,
\end{align*}
for some $\alpha \in (0,1)$. We infer that
\begin{align*}
\int_{\Omega \times \R^3} f(t) \vert v-u(t)\vert^2   \, \mathrm{d} v \, \mathrm{d} x \geq (1-\alpha) \int_{\Omega \times \R^3} f (t)\vert v\vert^2   \, \mathrm{d} v \, \mathrm{d} x  -(\alpha^{-1}-1) \int_{\Omega} \rho_f(t) \vert u(t) \vert^2    \, \mathrm{d} x.
\end{align*}
From the definition of $\widetilde{\E}(t)$, we can deduce that
\begin{align*}
\int_{\Omega \times \R^3} f(t) \vert u(t)-v\vert^2   \, \mathrm{d} v \,\mathrm{d} x &\geq (1-\alpha) \widetilde{\E}(t)-(\alpha^{-1}-1) \int_{\Omega} \rho_f \vert  u(t) \vert^2    \, \mathrm{d} x    \\ &\geq (1-\alpha) \widetilde{\E}(t)-(\alpha^{-1}-1)\underset{ s \in [0,T]}{\mathrm{sup}} \ \  \Vert \rho_f (s) \Vert_{\Ld^{\infty}(\Omega)}\Vert u(t) \Vert_{\Ld^2(\Omega)}^2.
\end{align*}
By our assumptions, the quantity $\underset{ s \in [0,T]}{\mathrm{sup}} \ \  \Vert \rho_f (s) \Vert_{\Ld^{\infty}(\Omega)}$ is finite and we can also assume that it is not equal to $0$. In order to get the inequality (\ref{inegtechnique:energiemod}), we eventually choose \begin{align*}
\alpha :=\left( 1+\dfrac{c_P}{2 \underset{ s \in [0,T]}{\mathrm{sup}} \ \  \Vert \rho_f (s) \Vert_{\Ld^{\infty}(\Omega)}  }\right)^{-1} \in \ (0,1),
\end{align*}
and then $\beta:=(\alpha^{-1}-1)\underset{ s \in [0,T]}{\mathrm{sup}} \ \  \Vert \rho_f (s) \Vert_{\Ld^{\infty}(\Omega)}=c_P/2$ and $\gamma:=1-\alpha$. Such an $\alpha$ is a continuous nondecreasing function of the variable $\underset{ s \in [0,T]}{\mathrm{sup}} \ \  \Vert \rho_f (s) \Vert_{\Ld^{\infty}(\Omega)}$ and does not vanish. Thus, we get $$\lambda := \min((1-\alpha)/2, c_P/2),$$ which is of the desired form. The proof is therefore complete.
\end{proof}

The following Proposition highlights the fact that a global bound on the local density $\rho_f$ is enough to obtain the convergence towards an asymptotic profile.

\begin{proposition}\label{asymp-profil}
For any weak solution $(u,f)$ to the Vlasov-Navier-Stokes system for which
\begin{align}
&\underset{t\geq 0}{\mathrm{sup}} \,  \Vert \rho_f (t) \Vert_{\Ld^{\infty}(\Omega)} < \infty, \label{hyprho:boundunif} 
\end{align}
there exists a profile $\rho^{\infty} \in \Ld^{\infty}(\Omega)$ such that
\begin{align}
\rho_f(t) \overset{t \rightarrow +\infty}{\longrightarrow} \rho^{\infty} \ \ in \ \H^{-1}(\Omega), 
\end{align}
exponentially fast.
\end{proposition}

\begin{proof}
It relies strongly on the exponential decay of the kinetic energy: indeed, as a consequence of Proposition \ref{decrexpo:energy}, we have $\E(t) \longrightarrow 0$ exponentially fast when $t \rightarrow +\infty$.  We will use Cauchy's criterion in $\H^{-1}(\Omega)$ in order to obtain the existence of an asymptotic profile.

To derive an estimate in $\H^{-1}$ norm, we fix $\psi \in \mathscr{C}_c^{\infty}(\Omega)$ and we use the weak formulation for the Vlasov equation (see Appendix \ref{DiPernaLions}) to write
\begin{align*}
\int_{\Omega}\psi \rho_f (t)-\int_{\Omega}\psi \rho_f (s)=\int_s^t \int_{\Omega} \nabla\psi \cdot j_f(\tau)   \, \mathrm{d}\tau
-\int_s^t \int_{\partial \Omega \times \R^3} (\gamma f) \psi(x) v \cdot n(x) \, \mathrm{d}\sigma(x) \, \mathrm{d}v  \, \mathrm{d}\tau,
\end{align*}
for all $0 \leq s \leq t$.
Since $\psi$ is compactly supported in $\Omega$, the last term vanishes. For the other term, we then use the Cauchy-Schwarz inequality to get, omitting the variable to simplify the formula
\begin{align*}
\left\vert \int_{\Omega} \psi \rho_f(t) \, \mathrm{d}x -\int_{\Omega}\psi \rho_f (s) \, \mathrm{d}x  \right\vert 
 \leq \int_s^t \left(\int_{\Omega } \rho_f  \vert \nabla \psi \vert^2  \, \mathrm{d}x \right)^{1/2} \left(\int_{\Omega \times \R^3} f  \vert v \vert^2 \, \mathrm{d}x  \, \mathrm{d}v \right)^{1/2} \, \mathrm{d}\tau,
\end{align*}
so that we eventually obtain
\begin{align*}
\left\vert \int_{\Omega} \psi \rho_f(t) \, \mathrm{d}x -\int_{\Omega}\psi \rho_f (s) \, \mathrm{d}x \right\vert 
 \leq \Vert \rho_f \Vert_{\Ld^{\infty}(\R^+;\Ld^{\infty}(\Omega))}^{1/2} \Vert \nabla \psi \Vert_{\Ld^2(\Omega)} \int_s^t \E(\tau)^{1/2}   \, \mathrm{d}\tau .
\end{align*}
Thanks to the assumption (\ref{hyprho:boundunif}),
we infer that 
\begin{align}\label{estimate:H-1}
\Vert \rho_f(t) - \rho(s) \Vert_{\H^{-1}(\Omega)} \lesssim \Vert \rho_f \Vert_{\Ld^{\infty}(\R^+;\Ld^{\infty}(\Omega))}^{1/2} \int_s^t \E(\tau)^{1/2}   \, \mathrm{d}\tau.
\end{align}
By the exponential decay of the energy on $\R^+$, the function $t \mapsto \E^{1/2}(t)$ in integrable on $\R^+$. Therefore, the Cauchy criterion for the $\Vert \cdot \Vert_{\H^{-1}(\Omega)}$ norm applied to the function $\rho_f(t)$ when $t \rightarrow + \infty$ gives us the convergence $\rho_f(t) \rightarrow \rho^{\infty} $ in $\H^{-1}(\Omega)$ when $t \rightarrow + \infty$, for some $\rho^{\infty} \in \H^{-1}(\Omega)$.

But again thanks to (\ref{hyprho:boundunif}), any sequence $(t_n)_n \rightarrow + \infty$ produces a sequence of functions $(\rho_f(t_n))_n$ bounded in $\Ld^{\infty}(\Omega)$. We deduce that $\rho^{\infty} \in \Ld^{\infty}(\Omega)$ and it is easy to see that the convergence towards this asymptotic profile is exponentially fast by looking at the estimate (\ref{estimate:H-1}) and letting $t \rightarrow + \infty$.
%
%
%
\end{proof}

\section{The particle trajectory }\label{Section:particle:traj}

In this section, we first derive a representation formula for the kinetic distribution thanks to the method of characteristics for the Vlasov equation, taking into account the absorption boundary condition. Then, we explain how to use a straightening change of variable in velocity in the integral formula giving $\rho_f$ to get a bound of the form
$$ \underset{0 \leq s \leq t}{\mathrm{sup}} \,  \Vert \rho_f (s) \Vert_{\Ld^{\infty}(\Omega)}\lesssim 1.$$
This idea has already been used in \cite{HKMM, HK} and can be applied if the change of variable is valid, which is the case if the semi norm $\Vert \nabla u \Vert_{\Ld^1(0,t;\Ld^{\infty}(\Omega))}$ is small enough.
\subsection{Characteristic curves for the system}
Given a time-dependent vector field $u$ on $\R^+ \times \Omega$, a time $t \in \R^+$ and a point $(x,v) \in \R^3 \times \R^3$, we define the characteristic curves $s \in \R^+ \mapsto (\mathrm{X}(s;t,x,v), \mathrm{V}(s;t,x,v)) \in \R^3 \times \R^3$ for the Vlasov equation (associated to $u$) as the solution of the following system of ordinary differential equations
\begin{equation}\label{EDO-charac}
\left\{
      \begin{aligned}
        \dot{\mathrm{X}}(s;t,x,v) &=\mathrm{V}(s;t,x,v),\\
\dot{\mathrm{V}} (s;t,x,v)&= (Pu)(s,\mathrm{X}(s;t,x,v))-\mathrm{V}(s;t,x,v),\\
	\mathrm{X}(t;t,x,v)&=x,\\
	\mathrm{V}(t;t,x,v)&=v,
      \end{aligned}
    \right.
\end{equation}
where the dot means derivative along the first
variable. Here, $P$ is the linear extension operator continuous from $\Ld^{\infty}(\Omega)$ to $\Ld^{\infty}(\R^3)$ and from $\W^{1,\infty}_0(\Omega)$ to $\W^{1,\infty}(\R^3)$ defined by 
\begin{equation}\label{opP:def}
\forall x \in \R^d, \ \ (Pw)(x):=\left\{
      \begin{aligned}
        & w(x), \ &&\text{if} \ x \in \Omega, \\
	& 0\ &&\text{if} \ x \in \R^3 \setminus \Omega,
      \end{aligned}
    \right.
\end{equation} and which satisfies
\begin{align}
& \forall w \in \Ld^{\infty}(\Omega), \ \Vert Pw\Vert_{\Ld^{\infty}(\R^3)} = \Vert w \Vert_{\Ld^{\infty}(\Omega)}, \label{opP:Linfini} \\[2mm]
& \forall w \in \W^{1,\infty}_0(\Omega), \ \ \Vert \nabla(Pw)\Vert_{\Ld^{\infty}(\R^3)} \leq \Vert \nabla w \Vert_{\Ld^{\infty}(\Omega)}. \label{opP:grad}
\end{align}
For the sake of completeness, we refer to Section \ref{Appendix:ExtensionLip} of the Appendix for a simple justification of these facts. Also, we will use the convention
\begin{align*}
(Pu)(t,\cdot)=P(u(t,\cdot)).
\end{align*}

\bigskip

Now, let $T>0$ be fixed and take $$u \in \Ld^2(0,T;\H^1_0(\Omega)) \cap \Ld^1(0,T;\W^{1,\infty}(\Omega)).$$ 
We can apply the Cauchy-Lipschitz theorem to show the following Proposition.

\begin{proposition}\label{Propo:diffeoZ}
Given $(x,v) \in \R^3 \times \R^3$ and a time $t \in [0,T]$, the system (\ref{EDO-charac}) admits a unique solution $s \mapsto \mathrm{Z}_{s,t}(x,v) \in \R^3 \times \R^3$ on $[0,T]$ and
\begin{equation*}
\mathrm{Z}_{s,t} : \left\{
   \begin{aligned}
      & \R^3 \times \R^3 &&\longrightarrow \R^3 \times \R^3\\
      & \ (x,v) &&\longmapsto \mathrm{Z}_{s,t}(x,v):= (\mathrm{X}(s;t,x,v),\mathrm{V}(s;t,x,v))
      \end{aligned}
    \right.
\end{equation*}
is a diffeomorphism of $\R^3 \times \R^3$ whose inverse is given by $\mathrm{Z}_{s,t}^{-1}=\mathrm{Z}_{t,s}$ and whose Jacobian determinant is $e^{3(t-s)}$.
\end{proposition}
In what follows, for all $(x,v) \in \R^3 \times \R^3$, we will sometimes use the notation
\begin{align}\label{notation:Z}
\mathrm{Z}_{s,t}(x,v):=(\mathrm{X}_{s,t}(x,v),\mathrm{V}_{s,t}(x,v)):=(\mathrm{X}(s;t,x,v),\mathrm{V}(s;t,x,v)) .
\end{align}

\subsection{A representation formula for the kinetic distribution}\label{Subsec:repform}
We keep the same assumptions and notations as in the previous subsection.
For $(x,v) \in \Omega \times \R^3$, we now define for any $t \geq 0$
\begin{align}\label{def:tau-}
\tau^{-}(t,x,v)&:=\inf \left\lbrace s \leq t \ \mid \forall \sigma \in [s,t], \  \mathrm{X}(\sigma;t,x,v) \in \Omega \right \rbrace,
\end{align}
where we use the harmless convention $\X(\sigma;t,x,v) =\X(0,t,x,v)$ for all $\sigma<0$ and $t \geq 0$.

If $t \geq 0$ is fixed, we set 
\begin{align}
\mathcal{O}^t &:= \left\lbrace  (x,v) \in \Omega \times \R^3 \ \mid \ \tau^{-}(t,x,v)<0  \right\rbrace. \label{equivOt}
\end{align}
We now derive an important representation formula for the distribution function which solves the Vlasov equation in the weak sense.
\begin{proposition}\label{Prop:formulerep}
Let $f$ be the weak solution to the initial boundary value problem for the Vlasov equation, associated to the velocity field $u \in  \Ld^2(0,T;\H^1_0(\Omega)) \cap \Ld^1(0,T;\W^{1,\infty}(\Omega))$ with initial condition $f_0$ and with absorption boundary condition. There holds
\begin{align}\label{formule-rep}
 f(t,x,v)= e^{3t} \mathbf{1}_{\mathcal{O}^t}(x,v) \, f_0(\mathrm{Z}_{0,t}(x,v)) \ \ \text{a.e.}
\end{align}
\end{proposition}
This formula seems to be folklore but we have not been able to find a proof in the literature. For the sake of completeness, and because this will be useful in the later Sections \ref{Section:asympprofile}-\ref{section:prescribe:mass}, we give a complete proof of Proposition \ref{Prop:formulerep} in Section \ref{Annexe:proofrep} of the Appendix.

A first application of this formula is the proof of Lemma \ref{loose:mass} that we had stated in Section \ref{section:energy}.

\begin{proof}[Proof of Lemma \ref{loose:mass}]
First suppose that $u$ and $f_0$ are smooth enough. By Proposition \ref{Prop:formulerep}, we have
\begin{align*}
 f(t,x,v)= e^{3t} \mathbf{1}_{\mathcal{O}^t}(x,v) \, f_0(\mathrm{Z}_{0,t}(x,v)) \ \ \text{a.e.}
\end{align*}
therefore 
\begin{align*}
\int_{\Omega \times \R^3} f(t,x,v) \, \mathrm{d}x \, \mathrm{d}v &=\int_{\Omega \times \R^3} e^{3t} \mathbf{1}_{\mathcal{O}^t}(x,v) \, f_0(\mathrm{Z}_{0,t}(x,v)) \, \mathrm{d}x \, \mathrm{d}v \\
&=\int_{\Omega \times \R^3}  \mathbf{1}_{\tau^{+}(0,x,v)>t} \, f_0(x,v) \, \mathrm{d}x \, \mathrm{d}v,
\end{align*}
thanks to the change of variable $(x,v) \mapsto \mathrm{Z}_{t,0}(x,v)$ (see Proposition \ref{Propo:diffeoZ} and the proof of Proposition \ref{Prop:formulerep} in Section \ref{Annexe:proofrep}). We thus deduce that
\begin{align*}
\int_{\Omega \times \R^3} f(t,x,v) \, \mathrm{d}x \, \mathrm{d}v \leq \int_{\Omega \times \R^3} f_0(x,v) \, \mathrm{d}x \, \mathrm{d}v.
\end{align*}

In the general case, we use a stability result coming from the DiPerna-Lions theory for transport equations  (see Section \ref{DiPernaLions} of the Appendix): we consider a sequence of nonnegative distribution functions $(f_n)_n$ solutions to the Vlasov equation with the absorption boundary condition, associated to regularized velocity fields $(u_n)_n$ and regularized and truncated initial conditions $(f_{0,n})_n$, converging respectively to $u$ and $f_0$. The associated characteristic curves $\mathrm{Z}^n $ for the Vlasov equation are classicaly defined and we have for all $n \in \N$
\begin{align*}
 f_n(t,x,v)= e^{3t} \mathbf{1}_{\mathcal{O}_n^t}(x,v) \, f_{0,n}(\mathrm{Z}^n_{0,t}(x,v)),
\end{align*}
where $\mathcal{O}_n^t$ is defined with respect to $\mathrm{Z}^n $ according to \eqref{equivOt}. 

We also know that $f$ is the strong limit in $\Ld^{\infty}_{\mathrm{loc}}(\R^+;\Ld^p(\Omega \times \R^3))$ for $1\leq p<\infty$ and the weak-$\star$ limit in $\Ld^{\infty}([0,T] \times \Omega \times \R^3)$ of the sequence $(f_n)_n$. This is thus enough to pass to the limit in the inequality 
\begin{align*}
\int_{\Omega \times \R^3} f_n(t,x,v) \, \mathrm{d}x \, \mathrm{d}v \leq \int_{\Omega \times \R^3} f_{0,n}(x,v) \, \mathrm{d}x \, \mathrm{d}v,
\end{align*}
which holds true thanks to the first part of the analysis. This concludes the proof because $f \in \mathscr{C}(\R^+;\Ld^1_{\mathrm{loc}}(\overline{\Omega} \times \R^3))$ (see \ref{DiPernaLions} in Appendix).
\end{proof}

\subsection{Change of variable in velocity and bounds on moments}
In order to get global bounds on the moments $\rho_f$ and $j_f$, we rely on a change of variable in velocity (inspired by \cite{BD}). As in \cite{HKMM,HK}, such a strategy can be used if the quantity $\Vert \nabla u \Vert_{\Ld^1(\R^+;\Ld^{\infty}(\Omega))}$ is small enough.
We use the same notations as before.
\begin{lemme}\label{chgmt-var-prepa}
Suppose $u \in \Ld^2_{\mathrm{loc}}(\R^+,\H^1_0(\Omega)) \cap\Ld^1_{\mathrm{loc}}(\R^+;\Ld^{\infty}(\Omega))$. Fix $\delta >0$ satisfying $\delta e^{\delta} \leq 1/9$. Then, for all times $t \in \R^{+}$ satisfiying
\begin{align}\label{borne-gradient0}
\int_0 ^t \Vert \nabla u(s) \Vert_{\Ld^{\infty}(\Omega)}  \, \mathrm{d}s < \delta,
\end{align}
and for all $x \in \Omega$, the map  
$$ \Gamma_{t,x} : v \mapsto \mathrm{V}(0;t,x,v),$$ is a global $\mathscr{C}^{1}$-diffeomorphism from $\R^3$ to itself satisfying furthermore 
\begin{align}\label{control:jacob}
\forall v \in \R^3, \ \ \vert \det \, \mathrm{D}_v \Gamma_{t,x} (v)  \vert \geq \dfrac{e^{3t}}{2}.
\end{align} 
\end{lemme}

\begin{proof}
The proof is a straightforward adaptation of that of \cite[Lemma 4.4]{HKMM}: here, the characteristic curves $\mathrm{Z}=(\X,\V)$ are defined in $\R^3 \times \R^3$ but we can use the inequality
$$\int_0^t \Vert \nabla (P u)(s)\Vert_{\Ld^{\infty}(\R^3)}   \, \mathrm{d}s \leq \int_0^t \Vert \nabla u(s)\Vert_{\Ld^{\infty}(\Omega)}   \, \mathrm{d}s.$$
which follows from (\ref{opP:grad}). The rest of the proof is then similar.
\end{proof}
As a consequence, we obtain the following proposition which is based on the representation formula for the  distribution function.

\begin{proposition}\label{Propo:gradFirst}
Suppose $u \in \Ld^2_{\mathrm{loc}}(\R^+,\H^1_0(\Omega)) \cap\Ld^1_{\mathrm{loc}}(\R^+;\Ld^{\infty}(\Omega))$. If the assumption (\ref{borne-gradient0}) is satisfied at a time $t\geq 0$, then we have
\begin{align*}
\Vert \rho_f(t) \Vert_{\Ld^{\infty}(\Omega)}  & \leq 2 I_q N_q(f_0), \\
\Vert j_f(t) \Vert_{\Ld^{\infty}(\Omega)} & \leq 2 I_q e^{-t} \left(\int_0 ^t e^s \Vert u(s) \Vert_{\Ld^{\infty}(\Omega)}   \, \mathrm{d}s +1 \right) N_q(f_0),
\end{align*}
where we recall that $N_q(f_0)$ is defined is (\ref{def:Nq})
and where
$$I_q:= \int_{\R^3} \dfrac{1+\vert v\vert}{1+\vert v \vert ^q}  \, \mathrm{d}v.$$
\end{proposition}
\begin{proof}
Note that because of the assumption (\ref{borne-gradient0}) on $u$, we can define the characteristics curves for the Vlasov-Navier-Stokes system in a classical sense as in the previous subsection. We use the representation formula (\ref{formule-rep}) to write that for all $x \in \Omega$
\begin{align*}
\rho_f(t,x)=e^{3t} \int_{\R^3} \mathbf{1}_{\mathcal{O}^t}(x,v) f_0(\mathrm{X}_{0,t}(x,v),\mathrm{V}_{0,t}(x,v)) \, \mathrm{d}v.
\end{align*}
Now, using the change of variable $v \mapsto \Gamma_{t,x}(v)= \mathrm{V}(0;t,x,v)$, we get
\begin{align*}
\rho_f(t,x)=e^{3t} \int_{\R^3} \mathbf{1}_{\mathcal{O}^t}(x,\Gamma_{t,x}^{-1}(w)) f_0(\mathrm{X}_{0,t}(x,\Gamma_{t,x}^{-1}(w)),w) \vert \det \, \mathrm{D}_w ( \Gamma_{t,x})^{-1}(w) \vert \, \mathrm{d}w,
\end{align*}
and thanks to the bound (\ref{control:jacob}), we obtain 
\begin{align*}
\rho_f(t,x) \leq 2 \int_{\R^3} \mathbf{1}_{\mathcal{O}^t}(x,\Gamma_{t,x}^{-1}(w)) f_0(\mathrm{X}_{0,t}(x,\Gamma_{t,x}^{-1}(w)),w)   \, \mathrm{d}w.
\end{align*}
Now, thanks to the definition (\ref{equivOt}) of $\mathcal{O}^t$,
$(x,\Gamma_{t,x}^{-1}(w)) \in \mathcal{O}^t$ if and only if $\tau^{-}(t,x,\Gamma_{t,x}^{-1}(w))<0$,
therefore 
\begin{align*}
(1+\vert w \vert ^q)\mathbf{1}_{\mathcal{O}^t}(x,\Gamma_{t,x}^{-1}(w)) f_0(\mathrm{X}_{0,t}(x,\Gamma_{t,x}^{-1}(w)),w) \leq  \underset{(x,v) \in \Omega \times \R^3}{\sup}(1+\vert v \vert ^q)f_0(x,v)=N_q(f_0).
\end{align*}
We thus deduce the desired inequality on $\rho_f$. We proceed in the same way for the bound on $j_f$, namely we start from the following representation formula 
\begin{align*}
j_f(t,x)=e^{3t} \int_{\R^3} \Gamma_{t,x}^{-1}(w) \mathbf{1}_{\mathcal{O}^t}(x,\Gamma_{t,x}^{-1}(w)) f_0(\mathrm{X}_{0,t}(x,\Gamma_{t,x}^{-1}(w)),w) \vert \det \, \mathrm{D}_w ( \Gamma_{t,x})^{-1}(w) \vert \, \mathrm{d}w.
\end{align*}
Using the formula 
\begin{align*}
 \mathrm{V}(0;t,x,v) = e^{t}v-\int_0^t e^{\tau}(P u)(\tau,\mathrm{X}(\tau;t,x,v)) \,  \mathrm{d}\tau,
\end{align*}
we infer that
\begin{align*}
w= e^{t}\Gamma_{t,x}^{-1}(w)-\int_0^t e^{\tau}(P u)(\tau,\mathrm{X}(\tau;t,x,\Gamma_{t,x}^{-1}(w))) \,  \mathrm{d}\tau,
\end{align*}
which becomes
\begin{align*}
\vert \Gamma_{t,x}^{-1}(w) \vert & \leq e^{-t} \left( \vert w \vert + \int_0^t e^{\tau} \Vert (P u)(\tau) \Vert_{\Ld^{\infty}(\R^3)} \,  \mathrm{d}\tau \right) \\
 &\leq e^{-t} (1 + \vert w \vert ) \left( 1 +\int_0^t e^{\tau} \Vert u(\tau) \Vert_{\Ld^{\infty}(\Omega)} \,  \mathrm{d}\tau \right).
\end{align*}
Coming back to the representation formula for $j_f$, we can conclude exactly as in the previous case.
\end{proof}
%
In the following Lemma, we study how the pointwise decay condition of Definition \ref{def:pointdecay} can be propagated after times $t=0$.
\begin{lemme}\label{inegdecal}
Let $t_0 >0$. If $N_q(f_0) < \infty$ and if $u \in \mathrm{L}_{\mathrm{loc}}^1(\R_+ ; \mathrm{H}^1_0 \cap \mathrm{L}^{\infty}(\Omega))$ then $N_q(f(t_0)) < \infty $ with 
\begin{align*}
N_q(f(t_0))  \lesssim e^{3 t_0} (1+ \Vert u \Vert_{\mathrm{L}^1 (0,t_0 ; \mathrm{L}^{\infty}(\Omega))} ^q) N_q(f_0).
\end{align*}
\end{lemme}
\begin{proof}
In this proof, we apply the stability results coming from the DiPerna-Lions theory (see Appendix \ref{DiPernaLions}) since the vector field $u$ has not enough regularity to define the characteristic curves in a classical sense. We omit the details and write the proof of the desired inequality with only sufficiently regular solutions and data (we refer to the proof of Lemma \ref{loose:mass} where we have already explained this regularization argument).

We use the representation formula
$$f(t_0,x,v)=e^{3t_0} \mathbf{1}_{\mathcal{O}^{t_0}}(x,v) f_0(\mathrm{Z}_{0,t_0}(x,v)).
 $$
Furthermore, we write
\begin{align*}
\mathrm{V}(0;t_0,x,v) &= e^{t_0}v -\int_0 ^{t_0} e^s (Pu)(s,\mathrm{X}(s;t_0,x,v))  \, \mathrm{d}s,
\end{align*}
to get 
$$v =e^{-t_0}\mathrm{V}(0;t_0,x,v)  + \int_0 ^{t_0} e^{s-t_0} (Pu)(s,\mathrm{X}(s;t_0,x,v))  \, \mathrm{d}s,  $$
and therefore obtain
$$ \vert v \vert \leq \vert \mathrm{V}(0;t_0,x,v) \vert + \int_0^{t_0} \Vert u(s) \Vert_{\Ld^{\infty}(\Omega)}   \, \mathrm{d}s.$$ Thus, we can estimate
\begin{multline*}
(1+\vert v \vert ^q)f(t_0,x,v)  \lesssim \left(1+\vert \mathrm{V}(0;t_0,x,v) \vert^q + \Vert u \Vert_{\mathrm{L}^1 (0,t_0 ; \mathrm{L}^{\infty}(\Omega))}^q \right)  e^{3 t_0}\mathbf{1}_{\mathcal{O}^{t_0}}(x,v) f_0(\mathrm{Z}(0;t_0,x,v)),
\end{multline*}
and the very definition of $N_q(f_0)$ leads to
\begin{align*}
(1+\vert v \vert ^q)f(t_0,x,v) \lesssim e^{3 t_0} (1+ \Vert u \Vert_{\mathrm{L}^1 (0,t_0 ; \mathrm{L}^{\infty}(\Omega))} ^q) N_q(f_0),
\end{align*}
which concludes the proof.

\end{proof}
By performing the same analysis as in the three previous results and by replacing the initial time $t=0$ by $t_0$ and $f_0$ by $f(t_0)$, we get the following statement.
\begin{lemme}\label{lastLMgradient}
Let $u \in \Ld^2_{\mathrm{loc}}(\R^+,\H^1_0(\Omega)) \cap\Ld^1_{\mathrm{loc}}(\R^+;\Ld^{\infty}(\Omega))$. Let $\delta$ be fixed such that $\delta e^{\delta} \leq 1/9$. For all times $t \geq t_0 \geq 0$ such that
\begin{align}\label{borne-gradient}
\int_{t_0} ^t \Vert \nabla u(s) \Vert_{\Ld^{\infty}(\Omega)}  \, \mathrm{d}s \leq \delta,
\end{align}
we have
\begin{align*}
\Vert \rho_f(t) \Vert_{\Ld^{\infty}(\Omega)}  & \lesssim e^{3t_0}N_q(f_0) (1+ \Vert u \Vert_{\mathrm{L}^1 (0,t_0 ; \mathrm{L}^{\infty}(\Omega))} ^q), \\
\Vert j_f(t) \Vert_{\Ld^{\infty}(\Omega)} & \lesssim  e^{3t_0-t} \left(\int_{t_0} ^t e^s \Vert u(s) \Vert_{\Ld^{\infty}(\Omega)}   \, \mathrm{d}s +1 \right) N_q(f_0) (1+ \Vert u \Vert_{\mathrm{L}^1 (0,t_0 ; \mathrm{L}^{\infty}(\Omega))} ^q).
\end{align*}
\end{lemme}

\begin{remarque}
We will actually see later that for any weak solution $(u,f)$ to the Vlasov-Navier-Stokes system, we have $u \in \Ld^1_{\mathrm{loc}}(\R^+;\Ld^{\infty}(\Omega))$: this comes from the local estimates of Proposition \ref{controlLinfiniLOC} which are independent of this section.
\end{remarque}

\section{Preparation for the bootstrap}\label{Section:prepbootstrap}
This section aims at providing local in time controls on $\rho_f$ and $j_f$ by using the results of the previous section and at deriving a $\Ld^{\infty}\H^1 \cap \Ld^2\H^2$ estimate for the fluid velocity. These two types of estimates will be crucial in order to show Theorem \ref{theoreme1}.

\bigskip

We first introduce the following useful notations, which allows us to track the dependency on the initial data in the later estimates. This will be useful at the very end of the our proof of Theorem \ref{theoreme1} because of the smallness condition (\ref{smallness:condition}).
\begin{notation}\label{nota:lesssim_0}
The notation $A \lesssim_0 B$ means
\begin{align}\label{Lesssim0}
A \leq \varphi \Big( 1+\Vert \nabla u_0 \Vert_{\Ld^2(\Omega)}+ N_q f_0 + \mathrm{E}(0) \Big)B,
\end{align}
where $\varphi : \R^+ \rightarrow \R^+$ is onto, continuous and nondecreasing, and $q>3$ is the exponent given in the statement of Theorem \ref{theoreme1}.
\end{notation}

From now on, we fix $(u,f)$ a weak solution to the Vlasov-Navier-Stokes system in the sense of Definition \ref{weak:sol} with admissible initial data $(u_0,f_0)$ in the sense of Definition \ref{CIadmissible}. Such a weak solution has for instance been constructed in \cite{BGM}.

\begin{notation}
We set
\begin{align*}
F&:=j_f-\rho_f u,\\
S&:=F- (u \cdot\nabla)u.
\end{align*}
\end{notation}

First, we have the following estimate on the Brinkman force $F$.
\begin{lemme}\label{ineg-Brinkman}
For all $t \geq 0$, we have
\begin{align*}
\int_0 ^t \Vert F(s) \Vert_{\Ld^2(\Omega)}^2   \, \mathrm{d}s \leq \mathrm{E}(0) \underset{s \in [0,t]}{\sup} \Vert \rho_f (s) \Vert_{\Ld^{\infty}(\Omega)}.
\end{align*}
\end{lemme}

\begin{proof}
Using Cauchy-Schwarz inequality, we get by dropping the time variable,
$$ \vert F\vert \ = \left\vert \int_{\R^3} f(v-u)  \, \mathrm{d}v \right\vert\leq \rho_f^{1/2} \left(\int_{\R^3} f\vert v-u\vert^2   \, \mathrm{d}v \right)^{1/2},$$
therefore for almost every $s \in [0,t]$
$$ \Vert F(s) \Vert_{\Ld^2(\Omega)}^2 \leq \Vert \rho_f (s) \Vert_{\Ld^{\infty}(\Omega)}\mathrm{D}(s) \leq\underset{s \in [0,t]}{\sup} \Vert \rho_f (s) \Vert_{\Ld^{\infty}(\Omega)} \mathrm{D}(s),$$
where $\mathrm{D}$ is defined in (\ref{eq:Dissipation}).
Then, we integrate the last inequality between $0$ and $t$ and conclude thanks to the energy inequality (\ref{ineq-energy}).
\end{proof}
\subsection{Local estimates}

We first recall standard interpolation estimates on the moments (see for instance \cite{Ham}) for which we refer to the notations of Definition \ref{notation:moments}.

\begin{proposition}\label{interpo-moment}
Let $k >0$ and let $g$ be a nonnegative function in $\Ld^{\infty}(\R^+ \times \Omega \times \R^3)$. Then the following estimates hold for any $\ell \in [0,k]$ and a.e $(t,x) \in \R^+ \times \Omega$
\begin{align*}
m_{\ell}g(t,x) &\lesssim \left( \Vert g\Vert_{\Ld^{\infty}(\R^+ \times \Omega \times \R^3)} +1 \right) m_{k}g(t,x)^{\frac{\ell +3}{k+3}},\\
\Vert m_{\ell}g(t) \Vert_{\Ld^{\frac{k+3}{\ell+3}}(\Omega)} &\leq C_{k,\ell} \Vert g(t) \Vert_{\Ld^{\infty}(\Omega \times \R^3)}^{\frac{k-\ell}{k+3}}  M_{k}g(t)^{\frac{\ell +3}{k+3}}.
\end{align*}
\end{proposition}

\begin{lemme}\label{interpo-estimate}
Consider $\alpha \geq 2$ such that $u \in \Ld^1_{\mathrm{loc}}(\R^+; \Ld^{\alpha+3} \cap \W^{1,1}(\Omega))$ and $M_{\alpha} f_0 < \infty$. Then $M_{\alpha} f(t) < \infty$ for all $t>0$ and we have
\begin{align}\label{interpo-moment-vlasov}
M_{\alpha} f(t) \lesssim \left( {M_{\alpha} f_0}^{\frac{1}{\alpha + 3}} + e^{\frac{3t}{\alpha+3}} \int_0^t \Vert u(s) \Vert_{\Ld^{\alpha+3}(\Omega)} \mathrm{d}s \right)^{\alpha+3}.
\end{align}
\end{lemme}
\begin{proof}
The proof basically follows the same arguments as in \cite[Lemma 4.2]{HKMM}. As in the proof of Lemma \ref{loose:mass}, we rely on the stability results of DiPerna-Lions theory on $\Omega \times \R^3$ (see Appendix \ref{DiPernaLions}). We do not detail the argument and we write the proof as if $u$ and $f_0$ were smooth, so that the characteristic curves for the Vlasov equation are classicaly defined. Furthermore, by Proposition \ref{Prop:formulerep}, we have
\begin{align*}
 f(t,x,v)= e^{3t} \mathbf{1}_{\mathcal{O}^t}(x,v) \,  f_0(\mathrm{Z}_{0,t}(x,v)),
\end{align*}
therefore
\begin{align*}
M_{\alpha} f(t)&= \int_{\Omega \times \R^3}  \vert v \vert^{\alpha} e^{3t} \mathbf{1}_{\mathcal{O}^t}(x,v) \,  f_0(\mathrm{Z}_{0,t}(x,v)) \, \mathrm{d}x \, \mathrm{d}v \\
&=\int_{\Omega \times \R^3}  \vert \mathrm{V}_{t,0}(x,v) \vert^{\alpha} \mathbf{1}_{\tau^{+}(0,x,v)>t} \, f_0(x,v) \, \mathrm{d}x \, \mathrm{d}v,
\end{align*}
thanks to the change of variable $(x,v) \mapsto \mathrm{Z}_{t,0}(x,v)$ (see Proposition \ref{Propo:diffeoZ} and the proof of Proposition \ref{Prop:formulerep} in Section \ref{Annexe:proofrep}). 

Since $\alpha \geq 2$, we also have
\begin{align*}
\frac{\mathrm{d}}{\mathrm{d} s} |\V_{s,0}(x,v))|^{\alpha} &= \alpha \frac{\mathrm{d}}{\mathrm{d} s} [\V_{s,0}(x,v))] \cdot \V_{s,0}(x,v)  |
\V(s;0,x,v) |^{\alpha-2} \\[2mm]
&= \alpha \left[ u(s,\X_{s,0}(x,v)) - \V_{s,0}(x,v) \right] \cdot \V_{s,0}(x,v)  |\V_{s,0}(x,v)|^{\alpha-2},
\end{align*}
from which we infer that
\begin{align*}
 |\V_{t,0}(x,v))|^{\alpha} \leq \vert v \vert^{\alpha}+\alpha\int_0^t \left[u(s,\X_{s,0}(x,v)) \right] \cdot \V_{s,0}(x,v)  \vert \V_{s,0}(x,v) \vert^{\alpha-2}\, \mathrm{d}s.
\end{align*}
By Fubini Theorem, we get
\begin{multline*}
M_{\alpha} f(t) \leq \int_{\Omega \times \R^3}  \vert v \vert^{\alpha} \mathbf{1}_{\tau^{+}(0,x,v)>t} \, f_0(x,v) \, \mathrm{d}x \, \mathrm{d}v  
\\ + \alpha \int_0^t \int_{\Omega \times \R^3}  \vert u(s,\X_{s,0}(x,v)) \vert \vert \V_{s,0}(x,v)  \vert \V_{s,0}(x,v) \vert^{\alpha-2} \vert \mathbf{1}_{\tau^{+}(0,x,v)>t} \, f_0(x,v) \, \mathrm{d}x \, \mathrm{d}v.
\end{multline*}
Using the reverse change of variable in the last integral, we get 
\begin{align*}
M_{\alpha} f(t) & \leq M_{\alpha} f_0+ \alpha \int_0^t \int_{\Omega} \vert u(s,x) \vert  m_{\alpha-1}f(s,x) \, \mathrm{d}x \, \mathrm{d}s \\
\quad & \leq M_{\alpha} f_0+\alpha \int_0^t \Vert u(s) \Vert_{\Ld^{\alpha +3}(\Omega)} \Vert m_{\alpha-1}f(s) \Vert_{\Ld^{\frac{\alpha+3}{\alpha+2}}(\Omega)} \, \mathrm{d}s,
\end{align*}
thanks to Hölder inequality.
The rest of the proof is then similar to that of \cite[Lemma 4.2]{HKMM}, by using the interpolation estimate on the moments (Proposition \ref{interpo-moment}) with $(\ell,k)=(\alpha-1,\alpha)$, the rough control
\begin{align*}
\Vert f(t) \Vert_{\Ld^{\infty}(\Omega \times \R^3)} \leq e^{3t} \Vert f_0 \Vert_{\Ld^{\infty}(\Omega \times \R^3)},
\end{align*}
and a Grönwall lemma in time.
\end{proof}

\begin{lemme}\label{rough:bounds}
If $M_3 f_0 <+ \infty$, we have $M_3 f\in \Ld^{\infty}_{\mathrm{loc}}(\R^+)$. Moreover, there exists a continuous non-negative and nondecreasing function $\eta$ such that
\begin{align*}
\Vert \rho_f (t) \Vert_{\Ld^{2}(\Omega)} + \Vert j_f (t) \Vert_{\Ld^{3/2}(\Omega)} \lesssim_0 \eta(t). 
\end{align*}
\end{lemme}
\begin{proof}
Thanks to the estimate (\ref{interpo-moment-vlasov}) with $\alpha=3$, we obtain
\begin{align*}
M_{3} f(t) \lesssim \left( {M_{3} f_0}^{1/6} + e^{t/2} \int_0^t \Vert u(s) \Vert_{\Ld^6(\Omega)} \mathrm{d}s \right)^6,
\end{align*}
as in \cite[Lemma 4.2]{HKMM}. Then, the rest of the proof is similar to that of \cite[Lemma 4.3]{HKMM}: we use Sobolev embedding and the Poincaré inequality on $\H^1_0(\Omega)$ to write $\Vert u(s) \Vert_{\Ld^6(\Omega)} \lesssim \Vert \nabla u(s) \Vert_{\Ld^2(\Omega)}$ so that the Cauchy-Schwarz inequality yields
\begin{align*}
\int_0^t \Vert u(s) \Vert_{\Ld^6(\Omega)} \mathrm{d}s \lesssim t^{1/2} \left( \int_0^t \Vert \nabla u(s) \Vert_{\Ld^2(\Omega)}^2 \mathrm{d}s\right)^{1/2} \leq t^{1/2} \mathrm{E}(0)^{1/2},
\end{align*}
where we have used the energy inequality (\ref{ineq-energy}).
\end{proof}

In the following lemma, we slightly improve the rough bound on $\rho_f$ and $j_f$ which was given in Lemma \ref{rough:bounds}.

\begin{lemme}\label{integBrinkman-1}
There exists a continuous nondecreasing function $\eta : \R^+\rightarrow \R^+$ such that for all $s \geq 0$ 
\begin{align}\label{ineg:momentexpo}
\Vert \rho_f (s) \Vert_{\Ld^{3}(\Omega)} + \Vert j_f (s) \Vert_{\Ld^{9/4}(\Omega)} \lesssim \eta(s).
\end{align}
Moreover, 
\begin{align*}
j_f - \rho_f u \in \Ld^2_{\mathrm{loc}}(\R^+;\Ld^2(\Omega)).
\end{align*}
\end{lemme}
\begin{proof}
In the proof we denote by $\eta$ a generic positive nondecreasing continuous function which may vary from line to line.

\medskip
Recall that in (\ref{data:hyp}), we suppose $M_2 f_0<\infty$ and $M_6 f_0<\infty$ so that $M_3 f_0 \lesssim M_2 f_0 +M_6 f_0 < \infty$. We first write that the velocity field $u$ is solution of the following Stokes system 
\begin{equation}
\left\{
      \begin{aligned}
        \partial_t u +Au &= \mathbb{P} S,\\
        \mathrm{div}_x \, u =0,\\
        u_{\mid \partial \Omega}=0,\\
        u(0)=u_0,
      \end{aligned}
    \right.
\end{equation}
with
$$S:=j_f -\rho_f u -(u\cdot \nabla)u,$$
where $\mathbb{P}$ stands for the Leray projection on vector fields and where $A$ is the Stokes operator (see Section \ref{AnnexeMaxregStokes} of the Appendix).
To estimate this source term, we first infer from Lemma \ref{rough:bounds} that 
\begin{align*}
\Vert \rho_f (t) \Vert_{\Ld^{2}(\Omega)} + \Vert j_f (t) \Vert_{\Ld^{3/2}(\Omega)} \lesssim_0 \eta(t). 
\end{align*}
The previous estimate, with the Hölder inequality and the energy inequality (\ref{ineq-energy}), allows us to obtain the following inequalities for all $t \geq 0$
\begin{align*}
\int_0^t \Vert j_f (s) \Vert_{\Ld^{3/2}(\Omega)}^2\mathrm{d}s &\lesssim_0 \eta(t),\\
\int_0^t \Vert \rho_f (s) u(s)\Vert_{\Ld^{3/2}(\Omega)}^2\mathrm{d}s \leq \int_0^t \Vert u(s) \Vert_{\Ld^6(\Omega)}^2 \Vert \rho_f(s) \Vert_{\Ld^2(\Omega)}^2 \mathrm{d}s    &\lesssim_0 \eta(t),
\end{align*}
where we have again used the Sobolev inequality $\Vert u(s) \Vert_{\Ld^6(\Omega)} \lesssim \Vert \nabla u(s) \Vert_{\Ld^2(\Omega)}$. Thus, we get $$ j_f -\rho_f u \in \Ld^{2}_{\mathrm{loc}}(\R^+;\Ld^{3/2}(\Omega)).$$
Then, since $u\in \Ld^{\infty}_{\mathrm{loc}}(\R^+;\Ld^2(\Omega))\cap  \Ld^{2}_{\mathrm{loc}}(\R^+;\H^1(\Omega))$ and again by Sobolev embedding, we get 
$$(u\cdot \nabla)u \in \Ld^{2}_{\mathrm{loc}}(\R^+;\Ld^1(\Omega)) \ \ \text{and} \ \ (u\cdot \nabla)u \in \Ld^{1}_{\mathrm{loc}}(\R^+;\Ld^{3/2}(\Omega)),$$
so that $(u\cdot \nabla)u \in \Ld^{p}_{\mathrm{loc}}(\R^+;\Ld^{q}(\Omega))$ with
\begin{align*}
\dfrac{1}{p}=\dfrac{\theta+1}{2}, \ \ \ \ \dfrac{1}{q}=\dfrac{3-\theta}{3}, \ \ \ \ 0 \leq \theta \leq 1.
\end{align*}
We choose $\theta =2/3$ which yields $(u\cdot \nabla)u \in \Ld^{6/5}_{\mathrm{loc}}(\R^+;\Ld^{9/7}(\Omega))$. We end up with
$$ \mathbb{P} S= \mathbb{P} [j_f -\rho_f u -(u\cdot \nabla)u] \in \Ld^{6/5}_{\mathrm{loc}}(\R^+;\Ld^{9/7}_{\mathrm{div}}(\Omega)),$$
by the continuity of the operator $\mathbb{P}$ from $\Ld^{9/7}(\Omega)$ to $\Ld^{9/7}_{\mathrm{div}}(\Omega)$.
Thanks to this integrability of the source term and the assumption (\ref{data:hyp}), we can use the maximal $\Ld^p-\Ld^q$ regularity result for the Stokes system (see Section \ref{AnnexeMaxregStokes} of the Appendix) to deduce that
$$ u \in  \Ld^{6/5}_{\mathrm{loc}}(\R^+;\mathrm{W}^{2,9/7}(\Omega)).$$
By the Sobolev embedding $\mathrm{W}^{2,9/7}(\Omega) \hookrightarrow \Ld^9(\Omega)$, we finally end up with
$$ u \in  \Ld^{6/5}_{\mathrm{loc}}(\R^+;\Ld^9(\Omega)).$$
In particular, we obtain $u \in  \Ld^{1}_{\mathrm{loc}}(\R^+;\Ld^{9}(\Omega))$. With the estimate (\ref{interpo-moment-vlasov}) of Lemma \ref{interpo-estimate}, we first have
\begin{align*}
M_{6} f(t) \lesssim \left( {M_{6} f_0}^{1/9} + e^{\frac{t}{3}} \int_0^t \Vert u(s) \Vert_{\Ld^{9}(\Omega)} \mathrm{d}s \right)^{9} \lesssim \eta(t),
\end{align*}
which can be used together with the interpolation estimate (\ref{interpo-moment}) for moments with $k=6$ and $\ell \in \{0,1\}$ to get
\begin{align*}
\Vert \rho_f (t) \Vert_{\Ld^{3}(\Omega)} + \Vert j_f (t) \Vert_{\Ld^{9/4}(\Omega)} \lesssim \eta(t).
\end{align*}
We can now combine the last estimate with the Hölder inequality and the energy inequality (\ref{ineq-energy}) to write for all $t \geq 0$
\begin{align*}
\int_0^t \Vert j_f (s) \Vert_{\Ld^{2}(\Omega)}^2\mathrm{d}s \lesssim \int_0^t \Vert j_f (s) \Vert_{\Ld^{9/4}(\Omega)}^2\mathrm{d}s \lesssim \eta(t),
\end{align*}
and
\begin{multline*}
\int_0^t \Vert \rho_f (s) u(s)\Vert_{\Ld^{2}(\Omega)}^2\mathrm{d}s \leq \int_0^t \Vert u(s) \Vert_{\Ld^6(\Omega)}^2 \Vert \rho_f(s) \Vert_{\Ld^3(\Omega)}^2 \mathrm{d}s \\ \lesssim  \int_0^t \Vert \nabla u(s) \Vert_{\Ld^2(\Omega)}^2 \Vert \rho_f(s) \Vert_{\Ld^3(\Omega)}^2 \mathrm{d}s  \lesssim \eta(t).
\end{multline*}
 This concludes the proof of Lemma \ref{integBrinkman-1}.
\end{proof}

The next proposition tells us that any weak solution $u$ to the VNS system belongs to $\Ld^{1}_{\mathrm{loc}}(\R^+;\Ld^{\infty}(\Omega))$. This result has been already proved in \cite[Proposition 5.1]{HKMM} in the case of the torus but the proof seems to be specific to the periodic setting. 

On a bounded domain, this type of integrability of Leray solutions to the Navier-Stokes equations is actually a general property which can be deduced from \cite[p104-106]{FT} combined with the end of the proof of \cite[Lemma 8.15]{Rob}. At first sight, the result stated in these references holds under the general condition that the source term of the Navier-Stokes equations belongs to $\Ld^{\infty}_{\mathrm{loc}}(\R^+;\Ld^2(\Omega))$. Even this is not (yet) the case for the Brinkman force $j_f -\rho_f u$, we will easily adapt the proof.

\begin{proposition}\label{controlLinfiniLOC}
We have 
\begin{align}
 u \in \Ld^{1}_{\mathrm{loc}}(\R^+;\Ld^{\infty}(\Omega)), \\ \label{borneL1Linfty-u}
\rho_f, j_f \in  \Ld^{\infty}_{\mathrm{loc}} (\R^+; \Ld^{\infty}(\Omega)).
\end{align}
More precisely, there exists a continuous nondecreasing function $\eta : \R^+\rightarrow \R^+$ such that
\begin{align}
\Vert u \Vert_{\Ld^{1}(0,t;\Ld^{\infty}(\Omega))} &\lesssim \eta(t), \\
\Vert \rho_f \Vert_{\Ld^{\infty}(0,t;\Ld^{\infty}(\Omega))} + \Vert j_f \Vert_{\Ld^{\infty}(0,t;\Ld^{\infty}(\Omega))}  & \lesssim_0 \eta(t). \label{borneLinftyLinfty-rhoetj}
\end{align}
\end{proposition}
\begin{proof}
\textbullet \ \ Let $T>0$ be fixed. We want to show that $u \in \Ld^{1}(0,T;\Ld^{\infty}(\Omega))$. As already said, some parts of the proof mimic those of \cite[p104-106]{FT} but we will introduce some modifications. 

First, thanks to Lemma \ref{integBrinkman-1}, we know that the source term $j_f-\rho_f u$ of the Navier-Stokes equations satisfied by $u$ belongs to $\Ld^{2}(0,T;\Ld^{2}(\Omega))$. So, we can apply the theory of epochs of regularity (see \cite[II - Section 7 ]{FT}) for $u$ on $[0,T]$. Note that this fact is valid because, by our assumption, $u$ is a Leray solution to the Navier-Stokes equations which satisfies the strong energy inequality.

We get the existence of a subset $\sigma_T \subset [0,T]$ of full measure in $[0,T]$ with $\sigma_T=\bigcup_i ]a_i, b_i[$ and where $u$ is a strong solution on each $]a_i, b_i[$ (namely, $u \in \Ld^{\infty}_{\mathrm{loc}}(a_i,b_i;\H^1(\Omega)) \cap \Ld^2_{\mathrm{loc}}(a_i,b_i;\H^2(\Omega))$). On each interval $]a_i, b_i[$, we can take the $\Ld^2(\Omega)$ inner product (denoted by $\left\langle \cdot, \cdot \right\rangle$) of the Navier-Stokes equations with $Au$, where $A$ stands for the Stokes operator on $\Ld^2(\Omega)$, in order to obtain
\begin{align*}
\dfrac{\mathrm{d}}{\mathrm{d}t}\Vert \nabla u \Vert_{\Ld^2(\Omega)}^2 + 2\Vert  Au \Vert_{\Ld^2(\Omega)}^2+ 2\langle  \mathbb{P} (u \cdot \nabla)u , Au\rangle=2\langle \mathbb{P}( j_f-\rho_f u), Au\rangle \ \ \text{on} \ \ ]a_i, b_i[,
\end{align*}
where we have used the fact that $\frac{\mathrm{d}}{\mathrm{d}t}\Vert \nabla u \Vert_{\Ld^2(\Omega)}^2=2\langle \partial_t u , Au \rangle$ (see \cite[Lemma 6.7]{Rob}). Here and in what follows, we omit the time variable for the sake of clarity. Then, we combine the Agmon inequality (see Section \ref{Agmon-ineq} of the Appendix), the Poincaré inequality and the elliptic estimate $\Vert u \Vert_{\H^2(\Omega)} \lesssim \Vert Au \Vert_{\Ld^2(\Omega)}$ to obtain
\begin{align*}
\vert \langle \mathbb{P} (u \cdot \nabla)u , Au\rangle \vert &\leq \Vert u \Vert_{\Ld^{\infty}(\Omega)} \Vert \nabla u \Vert_{\Ld^{2}(\Omega)} \Vert \mathbb{P} Au \Vert_{\Ld^{2}(\Omega)}\\
&\lesssim \Vert \nabla u \Vert_{\Ld^{2}(\Omega)}^{3/2} \Vert  \mathbb{P} A u \Vert_{\Ld^{2}(\Omega)}^{3/2},
\end{align*}
using that the operator $\mathbb{P}$ is self-adjoint on $\Ld^2(\Omega)$. The Young inequality is now used twice to write
\begin{align*}
2\vert \langle \mathbb{P} (u \cdot \nabla)u , Au\rangle \vert &\leq C\Vert \nabla u \Vert_{\Ld^{2}(\Omega)}^{6}+\dfrac{3}{4}\Vert A u \Vert_{\Ld^{2}(\Omega)}^{2},\\
2\vert \langle \mathbb{P}(j_f-\rho_f u), Au\rangle \vert  &\leq \dfrac{1}{2}\left(8C \Vert j_f-\rho_f u \Vert^2_{\Ld^2(\Omega)}+\dfrac{1}{2}\Vert A u \Vert^2_{\Ld^2(\Omega)} \right),
\end{align*}
thanks to the continuity of the operator $\mathbb{P}$ on $\Ld^2(\Omega)$. Here, $C$ is a positive constant independent of the time variable.
All in all, we end up with
\begin{align}\label{ineg-enstrophy}
\dfrac{\mathrm{d}}{\mathrm{d}t}\Vert \nabla u \Vert_{\Ld^2(\Omega)}^2 + \Vert  Au \Vert_{\Ld^2(\Omega)}^2\lesssim \Vert j_f-\rho_f u \Vert^2_{\Ld^2(\Omega)} + \Vert \nabla u \Vert_{\Ld^{2}(\Omega)}^{6},
\end{align}
on each interval $]a_i, b_i[$, where $\lesssim$ refers to a constant independent of the time variable and independent of $i$.

\medskip

To deal with the Brinkman force $\Vert j_f-\rho_f u \Vert^2_{\Ld^2(\Omega)}$, we invoke the inequality (\ref{ineg:momentexpo}) of Lemma \ref{integBrinkman-1} according to which there exists a continuous nondecreasing function $\eta : \R^+\rightarrow \R^+$ (independent of $T$) such that for all $s \in [0,T]$
\begin{align}\label{estimate-momentL^p}
\Vert \rho_f (s) \Vert_{\Ld^{3}(\Omega)} + \Vert j_f (s) \Vert_{\Ld^{9/4}(\Omega)} \lesssim \eta(s).
\end{align}
For all $s \in [0,T]$, we can now estimate the $\Ld^2$ norm of the Brinkman force $j_f(s) -\rho_f(s) u(s)$ in the following way,
\begin{multline*} 
\Vert j_f(s) -\rho_f(s) u(s) \Vert^2_{\Ld^2(\Omega)} 
\lesssim \Vert j_f(s) \Vert^2_{\Ld^2(\Omega)} + \Vert\rho_f(s) u(s) \Vert^2_{\Ld^2(\Omega)} \\
\lesssim \Vert j_f(s) \Vert^2_{\Ld^2(\Omega)} +  \Vert \rho_f(s) \Vert_{\Ld^3(\Omega)}^2  \Vert \nabla u(s) \Vert_{\Ld^2(\Omega)}^2,
\end{multline*}
where we have used the Hölder inequality and the Sobolev embedding ${\H}^1_0(\Omega) \hookrightarrow \Ld^6(\Omega)$. So, thanks to (\ref{estimate-momentL^p}), we obtain the following inequality for all $s \in [0,T]$
\begin{align*}
\Vert j_f(s) -\rho_f(s) u(s) \Vert^2_{\Ld^2(\Omega)}  \lesssim \eta(s)^2 + \eta(s)^2  \Vert \nabla u(s) \Vert_{\Ld^2(\Omega)}^2 \leq C_T \left(1+\Vert \nabla u(s) \Vert_{\Ld^2(\Omega)}^2\right).
\end{align*}
Coming back to (\ref{ineg-enstrophy}), we get 
\begin{align*}
\dfrac{\mathrm{d}}{\mathrm{d}t}\Vert \nabla u \Vert_{\Ld^2(\Omega)}^2 + \Vert  Au \Vert_{\Ld^2(\Omega)}^2\leq C_T \left(1+\Vert \nabla u \Vert_{\Ld^2(\Omega)}^2+ \Vert \nabla u \Vert_{\Ld^{2}(\Omega)}^{6}\right),
\end{align*}
on each $]a_i,b_i[$, where $C_T \geq 0$. We divide this inequality by $(1+\Vert \nabla u \Vert_{\Ld^2(\Omega)}^2)^2$ to find that
\begin{multline*}
\dfrac{1}{(1+\Vert \nabla u \Vert_{\Ld^2(\Omega)}^2)^2}\dfrac{\mathrm{d}}{\mathrm{d}t}\Vert \nabla u \Vert_{\Ld^2(\Omega)}^2 + \dfrac{1}{(1+\Vert \nabla u \Vert_{\Ld^2(\Omega)}^2)^2}\Vert  Au \Vert_{\Ld^2(\Omega)}^2 \\
\leq C_T \left(\dfrac{1}{(1+\Vert \nabla u \Vert_{\Ld^2(\Omega)}^2)^2}+ \dfrac{\Vert \nabla u \Vert_{\Ld^{2}(\Omega)}^{6}}{(1+\Vert \nabla u \Vert_{\Ld^2(\Omega)}^2)^2}\right),
\end{multline*}
which gives
\begin{multline}\label{ineg-enstrophy2}
\dfrac{1}{(1+\Vert \nabla u \Vert_{\Ld^2(\Omega)}^2)^2}\dfrac{\mathrm{d}}{\mathrm{d}t}\Vert \nabla u \Vert_{\Ld^2(\Omega)}^2 + \dfrac{1}{(1+\Vert \nabla u \Vert_{\Ld^2(\Omega)}^2)^2}\Vert  Au \Vert_{\Ld^2(\Omega)}^2 \\
 \leq C_T \left(1+ \Vert \nabla u \Vert_{\Ld^{2}(\Omega)}^{2}\right),
\end{multline}
on each $]a_i,b_i[$. We now integrate this inequality between $a$ and $b$, with $a_i<a<b<b_i$, to get
\begin{multline*}
\dfrac{1}{1+\Vert \nabla u(a) \Vert_{\Ld^2(\Omega)}^2}-\dfrac{1}{1+\Vert \nabla u(b) \Vert_{\Ld^2(\Omega)}^2} + \int_a^b\dfrac{\Vert  Au(s) \Vert_{\Ld^2(\Omega)}^2}{(1+\Vert \nabla u(s) \Vert_{\Ld^2(\Omega)}^2)^2}\, \mathrm{d}s \\
\leq C_T \int_a^b \left(1+ \Vert \nabla u(s) \Vert_{\Ld^{2}(\Omega)}^{2}\right) \, \mathrm{d}s.
\end{multline*}
As $]a_i, b_i[$ is a maximal interval of strong regularity for $u$ (if $b_i \neq T$), we know that $\Vert \nabla u(b) \Vert_{\Ld^2(\Omega)}$ goes to infinity when $b \rightarrow b_i$. Hence, if  $\ b_i \neq T$, we have
\begin{align*}
\int_{a_i}^{b_i}\dfrac{\Vert  Au(s) \Vert_{\Ld^2(\Omega)}^2}{(1+\Vert \nabla u(s) \Vert_{\Ld^2(\Omega)}^2)^2}\, \mathrm{d}s 
\leq C_T \int_{a_i}^{b_i} \left(1+ \Vert \nabla u(s) \Vert_{\Ld^{2}(\Omega)}^{2}\right) \, \mathrm{d}s.
\end{align*}
We sum the inequality on all $i$ (for which $b_i \neq T$) and use the fact that $\sigma_T$ is of full measure in $[0,T]$ to obtain
\begin{align*}
\int_0^T\dfrac{\Vert  Au(s) \Vert_{\Ld^2(\Omega)}^2}{(1+\Vert \nabla u(s) \Vert_{\Ld^2(\Omega)}^2)^2}\, \mathrm{d}s 
&\leq C_T \int_0^T \left(1+ \Vert \nabla u(s) \Vert_{\Ld^{2}(\Omega)}^{2}\right) \, \mathrm{d}s +\dfrac{1}{1+\Vert \nabla u(T) \Vert_{\Ld^2(\Omega)}^2} \\
&\leq C_T (T+\mathrm{E}(0))+1.
\end{align*}
We conclude the proof as in \cite[p104-106]{FT} and \cite[Lemma 8.15]{Rob} by writing
\begin{align*}
\int_0^T \Vert  Au(s) \Vert_{\Ld^2(\Omega)}^{2/3}\, \mathrm{d}s 
 &=\int_0^T \dfrac{\Vert  Au(s) \Vert_{\Ld^2(\Omega)}^{2/3}}{(1+\Vert \nabla u(s) \Vert_{\Ld^2(\Omega)}^2)^{2/3}}(1+\Vert \nabla u(s) \Vert_{\Ld^2(\Omega)}^2)^{2/3}\, \mathrm{d}s  \\
 &\leq \left( \int_0^T\dfrac{\Vert  Au(s) \Vert_{\Ld^2(\Omega)}^2}{(1+\Vert \nabla u(s) \Vert_{\Ld^2(\Omega)}^2)^2}\, \mathrm{d}s \right)^{1/3}\left( \int_0^T (1+\Vert \nabla u(s) \Vert_{\Ld^2(\Omega)}^2)\, \mathrm{d}s \right)^{2/3}\\
 & \leq \left( C_T (T+\mathrm{E}(0))+1\right)^{1/3} \left( T+\mathrm{E}(0)\right)^{2/3}<\infty,
\end{align*}
and by using again the Agmon inequality (\ref{Agmon-ineq}), the Poincaré inequality,  the elliptic estimate $\Vert u \Vert_{\H^2(\Omega)} \lesssim \Vert Au \Vert_{\Ld^2(\Omega)}$ and the Young inequality, we finally deduce that
\begin{align*}
\int_0^T \Vert  u(s) \Vert_{\Ld^{\infty}(\Omega)}\, \mathrm{d}s &\lesssim \int_0^T \Vert  \nabla u(s) \Vert_{\Ld^2(\Omega)}^{2}\, \mathrm{d}s +\int_0^T \Vert  Au(s) \Vert_{\Ld^2(\Omega)}^{2/3}\, \mathrm{d}s \\
& \leq \mathrm{E}(0)+\left( C_T (T+\mathrm{E}(0))+1\right)^{1/3} \left( T+\mathrm{E}(0)\right)^{2/3}<\infty.
\end{align*}
\textbullet \ \ To prove the last estimate (\ref{borneLinftyLinfty-rhoetj}), we use Lemma \ref{inegdecal} to write that for all $s \in [0,t]$
 \begin{align*}
\Vert \rho_f (s) \Vert_{\Ld^{\infty}(\Omega)} +\Vert j_f (s) \Vert_{\Ld^{\infty}(\Omega)} \lesssim N_q(f(s)) &\lesssim e^{3s}  (1+ \Vert u \Vert_{\mathrm{L}^1 (0,s ; \mathrm{L}^{\infty}(\Omega))} ^q) N_q(f_0)  \\
&\leq  e^{3t} (1+ \Vert u \Vert_{\mathrm{L}^1 (0,t ; \mathrm{L}^{\infty}(\Omega))} ^q) N_q(f_0),
\end{align*}
so that the conclusion follows.
\end{proof}

\subsection{Higher order energy estimates for the fluid velocity}
We now state a smoothing property of the Vlasov-Navier-Stokes system for the fluid velocity $u$. We rely on the parabolic regularization for the Navier-Stokes equations: in short, there is a gain of regularity if the initial data and the source term, that is the Brinkman force $F:=j_f-\rho_f u$, are small enough.
\begin{proposition}\label{Prop:H1estimates}
There exist universal constants $C_1,C_2>0$ such that the following holds. Assume that for some $T>0$, one has
\begin{align}\label{datasmall:VNSreg}
\Vert \nabla u_0 \Vert_{\Ld^2(\Omega)}^2 + C_1 \int_0^T \Vert F(s) \Vert_{\Ld^2(\Omega)}^2 \, \mathrm{d}s \leq \dfrac{1}{\sqrt{8C_1 C_2}}.
\end{align}
Then one has for all $0 \leq t \leq T$
\begin{align}\label{ineq:VNSreg}
\Vert \nabla u(t) \Vert_{\Ld^2(\Omega)}^2 + \dfrac{1}{2}\int_0^t \Vert Au(s) \Vert_{\Ld^2(\Omega)}^2 \, \mathrm{d}s \lesssim \Vert \nabla u_0 \Vert_{\Ld^2(\Omega)}^2 + \mathrm{E}(0)\underset{s \in [0,t]}{\sup} \Vert \rho_f (s) \Vert_{\Ld^{\infty}(\Omega)},
\end{align}
where $\lesssim$ only depends on $C_1$ and $C_2$, and where $A$ stands for the Stokes operator on $\Ld^2(\Omega)$.
\end{proposition}
\begin{proof}
The estimate is a direct consequence of the parabolic regularisation for the Navier-Stokes system with source $F=j_f-\rho_f u$, that we state in Theorem \ref{RegParabNS} in Section \ref{AnnexeParabNS} of the Appendix, together with the estimate on the Brinkman force of Lemma \ref{ineg-Brinkman}.
\end{proof}
\begin{remarque}
Thanks to the smallness assumption \eqref{smallness:condition}, we can actually ensure that
\begin{align}\label{Hyp:H1small}
\Vert \nabla u_0 \Vert_{\Ld^2(\Omega)}^2 \leq \dfrac{1}{2\sqrt{8C_1 C_2}}.
\end{align}
By choosing an appropriate function $\varphi$ in (\ref{smallness:condition}), we can also reduce $\Vert \nabla u_0 \Vert_{\Ld^2(\Omega)}^2$ and $\mathrm{E}(0)$ in the sequel if necessary.
\end{remarque}
\begin{remarque}\label{control:parab}
By Proposition \ref{controlLinfiniLOC}, the r.h.s of (\ref{ineq:VNSreg}) is finite. Using the elliptic estimate $\Vert u \Vert_{\H^2(\Omega)} \lesssim \Vert Au \Vert_{\Ld^2(\Omega)}$,  we infer that 
$$ u \in \Ld^{\infty}(0,T; \H^1(\Omega)) \cap \Ld^{2}(0,T;\H^2(\Omega)),
$$
and in particular $$ \nabla u \in \Ld^{2}(0,T;\H^1(\Omega)).$$
Thus, we get the following estimate
\begin{align}
\Vert u \Vert_{\Ld^{\infty}(0,t; \Ld^6(\Omega))}^2 \lesssim \Vert \nabla u \Vert_{\Ld^{\infty}(0,t; \Ld^2(\Omega))}^2 &\lesssim \Vert \nabla u_0 \Vert_{\Ld^2(\Omega)}^2 + \mathrm{E}(0)\underset{s \in [0,t]}{\sup} \Vert \rho_f (s) \Vert_{\Ld^{\infty}(\Omega)}. \label{controlparab1}
\end{align}
\end{remarque}

\medskip

In order to ensure that the smallness condition (\ref{datasmall:VNSreg}) is satisfied for all times, we now introduce the following terminology, which has been already used in \cite{HKMM} and \cite{HK} to take advantage of the parabolic regulazation for the fluid.
\begin{definition}[Strong existence time]\label{strongtime}
A real number $T\geq 0$ is a \emph{strong existence} time whenever (\ref{datasmall:VNSreg}) holds.
\end{definition}

\begin{lemme}\label{1:strongtime}
The smallness condition (\ref{smallness:condition}) of Theorem \ref{theoreme1} ensures that $T=1$ is a strong existence time in the sense of Definition \ref{strongtime}.
\end{lemme}
\begin{proof}
Recall the meaning of the notation $\lesssim_0$ from Notation \ref{nota:lesssim_0}. We use Lemma \ref{ineg-Brinkman} and the local estimate (\ref{borneLinftyLinfty-rhoetj}) to write
\begin{align*}
\int_0^1 \Vert F(s) \Vert_{\Ld^2(\Omega)}^2 \, \mathrm{d}s \leq  \mathrm{E}(0) \underset{s \in [0,1]}{\sup} \Vert \rho_f (s) \Vert_{\Ld^{\infty}(\Omega)} \lesssim_0  \mathrm{E}(0) \leq \dfrac{1}{2C_1 \sqrt{8C_1 C_2}},
\end{align*}
where we have used the assumption (\ref{smallness:condition}). Combining this inequality with (\ref{Hyp:H1small}) leads to the result.
\end{proof}

\section{Estimate for the second derivatives of the fluid velocity}\label{Section:estimateforces}

In this section, we provide a crucial estimate on the second order derivatives of $u$: to do so, we look at the convection term $\left( u \cdot \nabla \right) u$ and at the Brinkman force $F=j_f-\rho_f u$ as source terms for the Navier-Stokes equations. We rely on some maximal $\Ld^p\Ld^q$ regularity for the Stokes equation on a bounded domain in order to get estimates involving moments of $f$.

\medskip

We first introduce the following useful notations involving the moments of the kinetic distribution.
\begin{definition}
We set for all $t \geq 1$
\begin{align*}
\M_{\rho_f}(t) &:= \underset{s \in [1,t]}{\sup} \Vert \rho_f(s) \Vert_{\Ld^{\infty}(\Omega)}, \\ \M_{j_f}(t) &:= \underset{s \in [1,t]}{\sup} \Vert j_f(s) \Vert_{\Ld^{\infty}(\Omega)} ,\\
\M_{\rho_f,j_f}(t) &:=\M_{\rho_f}(t)+\M_{j_f}(t).
\end{align*}

\begin{remarque}\label{controlLinfiniLOC:constant}
Note that thanks to Proposition \ref{controlLinfiniLOC}, we can control $\rho_f(s)$ and $j_f(s)$ on $[0,1]$ in $\Ld^{\infty}(\Omega)$. Therefore, if $t>1$, 
we will make a constant use of
\begin{align*}
\underset{[0,t]}{\sup} \,  \Vert \rho_f(s) \Vert_{\Ld^{\infty}(\Omega)} &\lesssim_0 1 +\M_{\rho_f}(t), \\
\underset{[0,t]}{\sup} \,  \Vert j_f(s) \Vert_{\Ld^{\infty}(\Omega)} &\lesssim_0 1 +\M_{j_f}(t). 
\end{align*}
\end{remarque}

\end{definition}
\begin{proposition}\label{inegHessienn:expo}
Let $a,b,r \in ]1,\infty[$ and $\lambda>0$ fixed. For all $t \geq 1/2$ and all $0< q \leq a,b$, we have 
\begin{align*}
\begin{split}
 \int_{1/2}^t e^{- \lambda s} \Vert \mathrm{D}^2u(s) &\Vert_{\Ld^r(\Omega)}^q   \, \mathrm{d}s \lesssim_0 \Phi(\lambda) \left(1 +\Vert \left( u \cdot \nabla \right) u\Vert_{\Ld^a(0,t;\Ld^{r}(\Omega))}^q + \Vert F\Vert_{\Ld^b(0,t;\Ld^{r}(\Omega))}^q   \right), 
 \end{split}
\end{align*}
where $\Phi : \R^+ \rightarrow \R^+$ is nonincreasing.
\end{proposition}

\begin{proof}
We introduce $w_1$ et $w_2$ which are the unique divergence-free solutions on $[0,+\infty[$ to the following Cauchy problems
\begin{align}\label{Cauchyprob:w12}
 \left\{
    \begin{array}{ll}
         \partial_t w_1 +A w_1= -\mathbb{P} \left( u \cdot \nabla \right) u,\\
        w_1(0)=0, &
    \end{array}
\right.
\text{and}  \ \ \ \ \ \
\left\{
    \begin{array}{ll}
         \partial_t w_2 +A w_2 = \mathbb{P}F, \\
       w_2(0)=0,   &
    \end{array}
\right.
\end{align}
where $A$ stands for the Stokes operator on $\Ld^2(\Omega)$ and where $\mathbb{P}$ is the Leray projection on vector fields.

Thanks to the maximal regularity for the Stokes system and the continuity of the operator $\mathbb{P}$ on $\Ld^r(\Omega)$ (see Section \ref{AnnexeMaxregStokes} of the Appendix), we infer the following estimates for all $t \geq 0$
\begin{align}
 \left( \int_{0}^t  \Vert \mathrm{D}^2 w_1(s) \Vert_{\Ld^r(\Omega)}^a   \, \mathrm{d}s \right)^{1/a} &\lesssim  \left( \int_{0}^t  \Vert \left( u \cdot \nabla \right) u \Vert_{\Ld^r(\Omega)}^a   \, \mathrm{d}s \right)^{1/a} , \label{Hess-w1}\\ 
 \left( \int_{0}^t  \Vert \mathrm{D}^2 w_2(s) \Vert_{\Ld^r(\Omega)}^b   \, \mathrm{d}s \right)^{1/b} &\lesssim \left( \int_{0}^t  \Vert F(s) \Vert_{\Ld^r(\Omega)}^b   \, \mathrm{d}s \right)^{1/b}. \label{Hess-w2}
\end{align}

On the other hand, if we set $\widetilde{u}:=u-(w_1+w_2)$, we have $\widetilde{u}(t)=e^{-t A}u_0 \in D(A^s)$ for all $s \geq 0$. Let $(a_k)_{k \geq 0}$ be an Hilbertian basis of $\Ld^2_{\mathrm{div}}(\Omega)$ made up of eigenfunctions of $A$ with associated positive eigenvalues $(\lambda_k)_{k \geq 0}$ such that the sequence $(\lambda_k)_{k \geq 0}$ is nondecreasing and  $\lambda_k \overset{k \rightarrow + \infty}{\longrightarrow} + \infty$. We can write in $\Ld^2(\Omega)$
 $$ u_0(x)=\sum_{k =0}^{\infty} c_k a_k, \ \ c_k:=\langle u_0, a_k \rangle_{\Ld^2(\Omega)},$$
and the continuity of $e^{-t A}$ on $\Ld^2(\Omega)$ yields
$$\widetilde{u}(t,x)=\sum_{k =0}^{\infty} c_k e^{-\lambda_k t} a_k, \ \ t \geq 0.$$
Then, we use the fact for all  $t \geq 0$ and $\ell \geq 0$, we have 
$$ \Vert \widetilde{u}(t) \Vert_{{\H}^{\ell}(\Omega)} \lesssim \Vert A^{\ell /2}\widetilde{u}(t)  \Vert_{\Ld^2(\Omega)},$$
(see e.g: \cite[Chapter 2 - p69]{Rob}) so that we get for all $t \geq 0$ and $\ell \geq 1$
\begin{align*}
\Vert \widetilde{u}(t) \Vert_{{\H}^{\ell}(\Omega)}^2& \lesssim \sum_{k =0}^{\infty}\vert c_k \vert ^2 \lambda_k^{\ell} e^{-2 \lambda_k t}.
\end{align*}
Note that there exists a constant $C>0$ independent of $k$ and $t$ such that
$$ \lambda_k^{\ell} e^{-2 \lambda_k t} \leq C e^{- \lambda_k t}, \ \ t \geq 1/2, \ \ \ell \geq 1.$$
Indeed, for all $t \geq 1/2$ 
$$\lambda_k^{\ell} e^{- \lambda_k t}  \leq \lambda_k^{\ell} e^{- \lambda_k/2},$$
and the r.h.s of this inequality tends to $0$ as $ \rightarrow + \infty$ because $\lambda_k \overset{k \rightarrow + \infty}{\longrightarrow} + \infty$. Thus, for all $t \geq 1/2$ and $\ell \geq 1$, we have the following estimate
\begin{align*}
\Vert \widetilde{u}(t) \Vert_{{\H}^{\ell}(\Omega)}^2
\lesssim \left( \sum_{k =0}^{\infty} \vert c_k\vert^2 \right)  e^{-\lambda_1 t} 
=\Vert u_0 \Vert_{\Ld^2(\Omega)}^2 e^{-\lambda_1 t},
\end{align*}
where we have used the Plancherel-Parseval theorem.
We therefore have, for all $t \geq 1/2$ and $\ell \geq 1$
\begin{align}\label{borneHl}
\int_{1/2}^t  \Vert \widetilde{u}(s) \Vert_{{\H}^{\ell}(\Omega)} ^q   \, \mathrm{d}s \lesssim \Vert u_0 \Vert_{\Ld^2(\Omega)}^q \int_{1/2}^{+\infty} e^{-q\lambda_1 s/2}   \, \mathrm{d}s \lesssim \Vert u_0 \Vert_{\Ld^2(\Omega)}^q.
\end{align}
By using the estimate (\ref{borneHl}) with $\ell$ large enough, together with Sobolev embedding, we thus get
\begin{align}\label{Hess:ushift}
 \left(\int_{1/2}^t \Vert \mathrm{D}^2 \widetilde{u}(s)\Vert_{\Ld^r(\Omega)}^q   \, \mathrm{d}s\right) ^{1/q}\lesssim_0 1.
\end{align}
We then write $u=w_1+w_2+\widetilde{u}$ and it follows that
$$ \Vert \mathrm{D}^2u(s) \Vert_{\Ld^r(\Omega)}^q 
\lesssim \Vert\mathrm{D}^2 w_1(s) \Vert_{\Ld^r(\Omega)}^q 
+ \Vert \mathrm{D}^2 w_2(s) \Vert_{\Ld^r(\Omega)}^q
+\Vert \mathrm{D}^2 \widetilde{u}(s) \Vert_{\Ld^r(\Omega)}^q.$$
We have just dealt with the last term. For the other ones, and for $a\neq q, b\neq q$, we use the Hölder inequality, which is justified since $\frac{a}{q}, \frac{b}{q} \geq 1$, to write
\begin{align}
\label{ineg:integexpo1}\int_{1/2} ^t e^{-\lambda s}\Vert \mathrm{D}^2 w_1(s) \Vert_{\Ld^r(\Omega)}^q   \, \mathrm{d}s \leq \left( \int_0 ^t e^{-\lambda s\frac{a}{a-q}}   \, \mathrm{d}s \right)^{1-q/a} \left( \int_0 ^t \Vert \mathrm{D}^2 w_1(s) \Vert_{\Ld^r(\Omega)}^a   \, \mathrm{d}s \right)^{q/a}, \\
\label{ineg:integexpo2}\int_{1/2} ^t e^{-\lambda s}\Vert \mathrm{D}^2 w_2(s) \Vert_{\Ld^r(\Omega)}^q   \, \mathrm{d}s \leq \left( \int_0 ^t e^{-\lambda s\frac{b}{b-q}}   \, \mathrm{d}s \right)^{1-q/b} \left( \int_0 ^t \Vert \mathrm{D}^2 w_2(s) \Vert_{\Ld^r(\Omega)}^b   \, \mathrm{d}s \right)^{q/b}.
\end{align}
The first integral in the r.h.s of (\ref{ineg:integexpo1}) (resp of (\ref{ineg:integexpo2})) is equal to $\frac{a-q}{a}\frac{1}{\lambda}$ (resp equal to$\frac{b-q}{b}\frac{1}{\lambda}$) and these expression are indeed nonnegative and nonincreasing in $\lambda$.
If $a=q$, we have to to replace the previous inequality by
\begin{align*}
\int_{1/2} ^t e^{-\lambda s}\Vert \mathrm{D}^2 w_1(s) \Vert_{\Ld^r(\Omega)}^q   \, \mathrm{d}s \leq e^{-\lambda/2} \left( \int_0 ^t \Vert \mathrm{D}^2 w_1(s) \Vert_{\Ld^r(\Omega)}^a   \, \mathrm{d}s \right)^{q/a},
\end{align*}
and in a similar way for the term with $w_2$ if $b=q$.

Combining the inequalities (\ref{Hess-w1}), (\ref{Hess-w2}) et (\ref{Hess:ushift}), we end up with 
\begin{multline*}
\int_{1/2}^t e^{- \lambda s} \Vert \mathrm{D}^2u(s) \Vert_{\Ld^r(\Omega)}^q   \, \mathrm{d}s 
\\ \lesssim (\lambda^{-1}+e^{-\lambda/2} ) \left(\Vert \mathrm{D}^2 w_1\Vert_{\Ld^a(0,t;\Ld^{r}(\Omega))}^q +\Vert \mathrm{D}^2 w_2\Vert_{\Ld^b(0,t;\Ld^{r}(\Omega))}^q \right) + \Vert  \mathrm{D}^2 \widetilde{u}\Vert_{\Ld^q(1/2,t;\Ld^{r}(\Omega))}^q \\ 
\lesssim_0  \Phi(\lambda) \left(1 +\Vert \left( u \cdot \nabla \right) u\Vert_{\Ld^a(0,t;\Ld^{r}(\Omega))}^q + \Vert F\Vert_{\Ld^b(0,t;\Ld^{r}(\Omega))}^q   \right),
\end{multline*}
for some nonincreasing function $\Phi : \R^+ \rightarrow \R^+$.

\end{proof}

\begin{remarque}\label{remark:interval}
Note that we can also state a variant of the result of Proposition \ref{inegHessienn:expo} under the form 
\begin{align*}
\begin{split}
 \int_{1}^t e^{- \lambda s} \Vert \mathrm{D}^2u(s) &\Vert_{\Ld^r(\Omega)}^q   \, \mathrm{d}s \lesssim_0 \Psi(\lambda) \left(1 +\Vert \left( u \cdot \nabla \right) u\Vert_{\Ld^a(1/2,t;\Ld^{r}(\Omega))}^q + \Vert F\Vert_{\Ld^b(1/2,t;\Ld^{r}(\Omega))}^q   \right), 
 \end{split}
\end{align*}
for all $t \geq 1$ and for some nonincreasing function $\Psi : \R^+ \rightarrow \R^+$. To do so, we have to change the previous proof in the following way: we define $w_1$ and $w_2$ as solutions to the Cauchy problems (\ref{Cauchyprob:w12}) with $w_1(1/2)=w_2(1/2)=0$ and we still set $\widetilde{u}:=u-(w_1+w_2)$. We thus have $\widetilde{u}(t)=e^{-t A}u(1/2)$ and we can write 
$$ \widetilde{u}(t)=\sum_{k =1}^{\infty} \gamma_ke^{-\lambda_k (t-1/2)} a_k, \ \ \gamma_k:=\langle u(1/2), a_k \rangle_{\Ld^2(\Omega)}.$$
The rest of a the proof is then similar.
\end{remarque}

\medskip

In the two following lemmas, we give $\Ld^p \Ld^q$ estimates on the convection term and on the Brinkman force $F=j_f-\rho_f u$. These results are obtained in the very same way as in the case of the torus \cite[Lemmas 6.2 - 6.3]{HKMM} by using interpolation arguments and the parabolic regularization of (\ref{ineq:VNSreg}) which is valid for strong times. Indeed, Remark \ref{control:parab} entails
\begin{align*}
\Vert u \Vert_{\Ld^{\infty}(0,t; \Ld^6(\Omega))}^2 &\lesssim \Vert \nabla u_0 \Vert_{\Ld^2(\Omega)}^2 + \mathrm{E}(0)\underset{s \in [0,t]}{\sup} \Vert \rho_f (s) \Vert_{\Ld^{\infty}(\Omega)} \\
&\lesssim (\Vert \nabla u_0 \Vert_{\Ld^2(\Omega)}^2 + \mathrm{E}(0) )(1+ \underset{s \in [0,t]}{\sup} \Vert \rho_f (s) \Vert_{\Ld^{\infty}(\Omega)}) \\
&\lesssim_0 1+\underset{s \in [0,t]}{\sup} \Vert \rho_f (s) \Vert_{\Ld^{\infty}(\Omega)}.
\end{align*}
This yields the following results.
\begin{lemme}\label{interpo:convection}
There exist $a \in (2,4)$ and $r_a \in (2,3)$ such that for all strong existence times $t \geq 0$
\begin{align*}
\Vert \left( u \cdot \nabla \right) u\Vert_{\Ld^a(0,t;\Ld^{r_a}(\Omega))} \lesssim_0 1 +\M_{\rho_f,j_f}(t).
\end{align*}
\end{lemme}

\begin{lemme}\label{interpo:Brinkman}
For all $b >4$, there exists $r_b >4$ such that for all strong existence times $t \geq 0$
\begin{align*}
\Vert F\Vert_{\Ld^b(0,t;\Ld^{r_b}(\Omega))} \lesssim_0 1 +\M_{\rho_f,j_f}(t)^{\frac{3}{2}-\frac{2}{b}}.
\end{align*}
\end{lemme}

\section{End of the proof of Theorem \ref{theoreme1}}\label{section:bootstrap}
We are now able to set up the bootstrap procedure we have mentioned in the end of the Introduction. In order to get a control on $\nabla u$, we interpolate the higher regularity estimates with the poinstwise $\Ld^2(\Omega)$ bound on $u$ provided by the exponential decay of the total kinetic energy.

We first state the following result on $\nabla u$, which is non-uniform in time for the moment.
\begin{lemme}\label{grad:nonunif}
For any strong existence time $t \geq 1$, one has
\begin{align*}
\nabla u \in \Ld^1(1,t;\Ld^{\infty}(\Omega)).
\end{align*}
\end{lemme}
\begin{proof}
Let $t \geq 1$ be a strong existence time. Since the trace of $u$ is $0$ on $\partial \Omega$, we use the Gagliardo-Nirenberg-Sobolev inequality (see Theorem \ref{gagliardo-nirenberg} in Appendix) with $(j,m,q)=(1,2,2)$ to write
\begin{align*}
\Vert \nabla u(s) \Vert_{\Ld^p(\Omega)} \lesssim \Vert \mathrm{D}^2 u(s) \Vert_{\Ld^r(\Omega)} ^{\alpha} \Vert u(s) \Vert_{\Ld^2(\Omega)}^{1-\alpha}, \ \ \ \ s\geq 1/2,
\end{align*}
for all $p \in [1,\infty]$ and $r \in [1,\infty]$ satisfying the relation
\begin{align}\label{gagliardo-grad}
\dfrac{1}{p}=\dfrac{1}{3}+\alpha\left( \dfrac{1}{r}-\dfrac{2}{3} \right)+\dfrac{1-\alpha}{2},
\end{align}
and where $\alpha \in [1/2,1[$. With the energy inequality (\ref{ineq-energy}) and a Hölder inequality in time, we get
\begin{align*}
\int_{1/2}^t \Vert \nabla u(s) \Vert_{\Ld^p(\Omega)}^c \, \mathrm{d}s \leq \E(0)^{c\frac{1-\alpha}{2}}c^{-1}(e^{ct}-e^{\frac{c}{2}})\int_{1/2}^t e^{-c  s} \Vert \mathrm{D}^2 u(s) \Vert_{\Ld^r(\Omega)} ^{c \alpha}\, \mathrm{d}s,
\end{align*}
which turns into
\begin{multline*}
\int_{1/2}^t \Vert \nabla u(s) \Vert_{\Ld^p(\Omega)}^c \, \mathrm{d}s \\
\lesssim_0 \E(0)^{c\frac{1-\alpha}{2}}c^{-1}(e^{ct}-e^{\frac{c}{2}}) \Phi(c)  \left(1 +\Vert \left( u \cdot \nabla \right) u\Vert_{\Ld^a(0,t;\Ld^{r}(\Omega))}^{c\alpha} + \Vert F\Vert_{\Ld^b(0,t;\Ld^{r}(\Omega))}^{c\alpha}  \right),
\end{multline*}
thanks to Proposition \ref{inegHessienn:expo}, for all $c \in [1,+\infty[$ such that $c \alpha \leq a,b$ and exponents $1<a,b<\infty$. Now using Lemma \ref{interpo:convection} and Lemma \ref{interpo:Brinkman}, we obtain $a,b,r_a,r_b$ such that $b>4>a>2$ and $r:=\min(r_a,r_b)>2$ and for which we can write
\begin{align}\label{boundgrad:nonunif}
\begin{split}
& \int_{1/2}^t \Vert \nabla u(s) \Vert_{\Ld^p(\Omega)}^c \, \mathrm{d}s
 \lesssim_0 \gamma(t), \\[2mm]
\gamma(t):=\E(0)^{c\frac{1-\alpha}{2}}c^{-1}(e^{ct}-e^{\frac{c}{2}})& \Phi(c)  \left[1 +(1 +\M_{\rho_f,j_f}(t))^{c\alpha} + \left(1 +\M_{\rho_f,j_f}(t)^{\frac{3}{2}-\frac{2}{b}}\right)^{c\alpha}  \right],
\end{split}
\end{align}
provided that $\alpha \in [1/2,1[$ and $p \in [1,\infty]$ satisfy $\alpha c \leq \min(a,b)$ and the relation (\ref{gagliardo-grad}).

However, reaching $p=\infty$ is not possible yet because the relation (\ref{gagliardo-grad}) would imply
$$ \alpha=5 \left( 7-\dfrac{6}{r} \right)^{-1},$$
so that the condition  $\alpha \in [1/2,1[$ is actually equivalent to the condition $r>3$. This last inequality is not \textit{a priori} satisfied by our choice $r=\min(r_a,r_b)$. Nevertheless, we will replace $r_a$ by a larger exponent. To do so, we first use (\ref{boundgrad:nonunif}) with $c=a<b$ and by carefully looking at the relation (\ref{gagliardo-grad}) when $\alpha$ is close enough to $1$, we see there exists $\alpha \in [1/2,1[$ related to $p>6$ by (\ref{gagliardo-grad}) and such that 
\begin{align*}
\Vert \nabla u \Vert_{\Ld^a(1/2,t;\Ld^p(\Omega))} \lesssim_0 \gamma(t).
\end{align*}
Since $p>6$, we now use the Hölder inequality to find $\tilde{r}_a:= 6p/(p-6)>3$ and write 
\begin{align*}
\left( \int_{1/2}^t \Vert \left( u \cdot \nabla \right) u(s)\Vert_{\Ld^{\tilde{r}_a}(\Omega)}^a \, \mathrm{d}s  \right)^{1/a} &\leq  \Vert u\Vert_{\Ld^{\infty}(1/2,t;\Ld^{6}(\Omega))}  \Vert \nabla u\Vert_{\Ld^a(1/2,t;\Ld^{p}(\Omega))}  \\
&   \leq \gamma(t) \Vert u\Vert_{\Ld^{\infty}(0,t;\Ld^{6}(\Omega))}  \\
& \leq  \gamma(t) \left(\Vert \nabla u_0 \Vert^2_{\Ld^2(\Omega)}+\mathrm{E}(0) \eta(t)\right),
\end{align*}
thanks to the estimates (\ref{controlparab1}) and (\ref{borneLinftyLinfty-rhoetj}).

It allows us to replace Lemma \ref{interpo:convection} by the previous inequality so that we can get $\tilde{r}_a>3$ instead of $r_a$. Now, this yields $\tilde{r}:=\min(r_b,\tilde{r}_a)>3$ and by taking
$$ \tilde{\alpha}=5 \left( 7-\dfrac{6}{\tilde{r}} \right)^{-1},$$
we have $$\tilde{\alpha}  \in[1/2,1[, \ \ 0=\dfrac{1}{3}+\tilde{\alpha}\left( \dfrac{1}{\tilde{r}}-\dfrac{2}{3} \right)+\dfrac{1-\tilde{\alpha}}{2}.$$
Arguing as in the beginning of the proof, we now use Remark \ref{remark:interval} with $\tilde{r}>3$, $c=1$ and $\tilde{\alpha}$ and we get
\begin{align*}
\int_1^t \Vert \nabla u(s)\Vert_{\Ld^{\infty}(\Omega)} \, \mathrm{d}s \lesssim_0 \tilde{\gamma}(t),
\end{align*}
where $$ \widetilde{\gamma}(t):=\E(0)^{\frac{1-\tilde{\alpha}}{2}}(e^{t}-e) \Phi(1)  \left[1 +(1 +\M_{\rho_f,j_f}(t))^{\tilde{\alpha}} + \left(1 +\M_{\rho_f,j_f}(t)^{\frac{3}{2}-\frac{2}{b}}\right)^{\tilde{\alpha}}  \right].$$
Because of $(\ref{borneLinftyLinfty-rhoetj})$, the quantity $\widetilde{\gamma}(t)$ is finite so that this concludes the proof.

\end{proof}

\bigskip

In order to set up a bootstrap argument, we naturally introduce the following quantity:
\begin{align}\label{def:tstar}
t^{\star} :=\sup \left\lbrace \text{strong existence times } t \in \R^+  \text{ such that } \int_1 ^t \Vert \nabla u(s) \Vert_{\Ld^{\infty}(\Omega)}  \, \mathrm{d}s < \delta \right\rbrace,
\end{align}
where $\delta >0$ is taken small enough in order to satisfy $\delta e^{\delta}\leq 1/9$. Our main goal is now to show that $t^{\star} = + \infty$.

\begin{lemme}\label{tstar>1}
We have $t^{\star}>1$ and for any $t<t^{\star}$, the estimate $\mathrm{M}_{\rho_f,j_f}(t) \lesssim_0 1$ holds.
\end{lemme}
\begin{proof}
By reducing $\E(0)$ in Lemma \ref{ineg-Brinkman} and by the same proof as in Lemma \ref{1:strongtime} , we observe that $t=1+\varepsilon$ is still a strong existence time for $\varepsilon$ small enough. Thus, thanks to Lemma \ref{grad:nonunif}, we can find $\varepsilon$ small enough such that $1+\varepsilon$ is a strong existence time and such that $\Vert \nabla u \Vert_{\Ld^1(1,1+\varepsilon;\Ld^{\infty}(\Omega))} <\delta$. It therefore implies that $t^{\star}>1$.

Now, we take $t\in [1,t^{\star}[$ and we write it as $t=t^{\star}-r$ with $r>0$. By definition of $t^{\star}$, there exists a time $\tilde{t}$ such that $t=t^{\star}-r<\tilde{t}<t^{\star}$ and such that $\tilde{t}$ is a strong existence satisfying $\Vert \nabla u \Vert_{\Ld^1(1,\tilde{t};\Ld^{\infty}(\Omega))}<\delta $. Now, we can use Lemma \ref{lastLMgradient} with $t_0=1$ and the estimate (\ref{borneL1Linfty-u}) for $u$ to get
$$ \Vert \rho_f(t) \Vert_{\Ld^{\infty}(\Omega)}  \lesssim N_q(f_0) (1+ \Vert u \Vert_{\mathrm{L}^1 (0,1 ; \mathrm{L}^{\infty}(\Omega))} ^q) \lesssim_0 N_q(f_0) (1+ \eta(1) ^q), $$
therefore $\M_{\rho_f}(t) \lesssim_0 1$, uniformly in $t$.
Similarly, for $j_f$, we have
\begin{align*}
\Vert j_f(t) \Vert_{\Ld^{\infty}(\Omega)} &\lesssim  e^{-t} \left(\int_{1} ^t e^s \Vert u(s) \Vert_{\Ld^{\infty}(\Omega)}   \, \mathrm{d}s +1 \right) N_q(f_0) (1+ \Vert u \Vert_{\mathrm{L}^1 (0,1 ; \mathrm{L}^{\infty}(\Omega))} ^q) \\
&\lesssim_0 e^{-t} \left(\int_{1} ^t e^s \Vert u(s) \Vert_{\Ld^{\infty}(\Omega)}   \, \mathrm{d}s +1 \right).
\end{align*}
The Sobolev embedding $\H^2(\Omega) \hookrightarrow \Ld^{\infty}(\Omega)$ and the elliptic estimate $\Vert u \Vert_{\H^2(\Omega)} \lesssim \Vert Au \Vert_{\Ld^2(\Omega)}$ ensure that for all $t \in [1,t^{\star}[$
\begin{align*}
 e^{-t}\int_{1} ^t e^s \Vert u(s) \Vert_{\Ld^{\infty}(\Omega)}  \, \mathrm{d}s  
\lesssim  e^{-t} \int_{1} ^t e^s \Vert Au(s) \Vert_{\Ld^{2}(\Omega)}   \, \mathrm{d}s.
\end{align*}
Thanks to Cauchy-Schwarz inequality, we get for all $t \in [1,t^{\star}[$ 
\begin{align*}
e^{-t} \int_{1} ^t e^s \Vert u(s) \Vert_{\Ld^{\infty}(\Omega)}   \, \mathrm{d}s  &\leq
e^{-t}\left( \int_{1}^{t} e^{2s} \, \mathrm{d}s \right)^{1/2} \left(\int_{1} ^t  \Vert A u(s) \Vert_{\Ld^{2}(\Omega)}^2   \, \mathrm{d}s \right)^{1/2}  \\
&  \lesssim \left(\int_{1} ^t  \Vert A u(s) \Vert_{\Ld^{2}(\Omega)}^2   \, \mathrm{d}s \right)^{1/2} \\
& \lesssim_0 \left(1+\underset{s \in [0,t]}{\sup} \Vert \rho_f (s) \Vert_{\Ld^{\infty}(\Omega)} \right)^{1/2} \\
& \lesssim_0 \left(1 + \M_{\rho_f}(t) \right)^{1/2},
\end{align*}
where we have used the parabolic estimate (\ref{ineq:VNSreg}) (which holds since $t<\tilde{t}$). This concludes the proof as we have already proved that $\M_{\rho_f}(t) \lesssim_0 1$.
\end{proof}
\begin{remarque}\label{remark:M}
Note that the estimate we have just proved is uniform in time. Therefore, by considering for example a sequence $(t_n)_n$ of strong existence times with $\Vert \nabla u \Vert_{\Ld^1(1,t_n;\Ld^{\infty}(\Omega))} <\delta $ such that $t_n\rightarrow t^{\star}$ and $t_n \leq t^{\star}$, we get $\mathrm{M}_{\rho_f,j_f}(t^{\star}) \lesssim_0 1$.
\end{remarque}

\begin{proposition}\label{grad:unif}
If $t^{\star}<\infty$, there exists $\gamma >0$ such that the following estimate holds
\begin{align}\label{boundgrad-bootstrap}
\int_1^{t^{\star}} \Vert \nabla u(s) \Vert_{\Ld^{\infty}(\Omega)} \, \mathrm{d}s \lesssim_0 \E(0)^{\gamma}.
\end{align}
\end{proposition}
\begin{proof}
We modify the proof of Lemma \ref{grad:nonunif} and take advantage of the exponential decay of the kinetic energy on $[1,t^{\star}]$ in order to get uniform in time estimates. We start again using the Gagliardo-Nirenberg-Sobolev inequality (see Theorem \ref{gagliardo-nirenberg} in Appendix) with $(j,m,q)=(1,2,2)$ to write
\begin{align*}
\Vert \nabla u(s) \Vert_{\Ld^p(\Omega)} \lesssim \Vert \mathrm{D}^2 u(s) \Vert_{\Ld^r(\Omega)} ^{\alpha} \Vert u(s) \Vert_{\Ld^2(\Omega)}^{1-\alpha}, \ \ \ \ s\geq 1/2,
\end{align*}
for all $p \in [1,\infty]$ and $r \in [1,\infty]$ satisfying the relation (\ref{gagliardo-grad})
and where $\alpha \in [1/2,1[$. Now, since $t^{\star}<\infty$, we can combine Proposition \ref{decrexpo:energy} on $[0,t^{\star}]$ and Lemma \ref{tstar>1} to get a rate $\lambda^{\star}$ such that $\E(t) \leq e^{- \lambda^{\star} t} \E(0)$ on $[0,t^{\star}]$. By looking at the definition of the kinetic energy and by setting $\lambda:=\lambda^{\star}(1-\alpha)/2$, we have
\begin{align}\label{boundgrad:unif}
\Vert \nabla u(s) \Vert_{\Ld^p(\Omega)} \lesssim \E(0)^{\frac{1-\alpha}{2}} e^{-\lambda s} \Vert \mathrm{D}^2 u(s) \Vert_{\Ld^r(\Omega)} ^{\alpha} , \ \ \ \ s \in [1/2,t^{\star}].
\end{align}
This inequality is the key to get uniform in time estimates.

Moreover, by taking a sequence $(t_n)_n$ of strong existence times with $\Vert \nabla u \Vert_{\Ld^1(1,t_n;\Ld^{\infty}(\Omega))} <\delta $ such that $t_n\rightarrow t^{\star}$ and $t_n \leq t^{\star}$, we see that estimate (\ref{ineq:VNSreg}) holds for $t^{\star}$ thanks to Lemma \ref{tstar>1} so that the statements of Proposition \ref{inegHessienn:expo} and Lemmas \ref{interpo:convection}-\ref{interpo:Brinkman} still hold at time $t=t^{\star
}$.

%
%

We now perfom the same arguments as in the proof of Lemma \ref{grad:nonunif}: namely, we replace the inequality (\ref{boundgrad:nonunif}) by the inequality (\ref{boundgrad:unif}).  We eventually end up with a control on $\nabla u$ of the form (\ref{boundgrad-bootstrap}) and this concludes the proof.
\end{proof}

\medskip

In view of Lemma \ref{tstar>1}, it remains to show the following statement.
\begin{proposition}
We have $t^{\star}=+\infty$.
\end{proposition}
\begin{proof}
By contradiction, let assume that $t^{\star}< + \infty$. We will get a contradiction by proving (provided that $\E(0)$ is small enough) the existence of a strong existence  time $t>t^{\star}$ such that the estimate $\Vert \nabla u \Vert_{\Ld^1(1,t;\Ld^{\infty}(\Omega)}<\delta$ holds.

\textbullet \ \ We first prove that there exists a strong existence time larger than $t^{\star}$. Recall that we work under the assumption $\Vert \nabla u_0 \Vert_{\Ld^2(\Omega)}^2 \leq (2 \sqrt{8C_1 C_2})^{-1}$ (see (\ref{Hyp:H1small})). Recall also that the estimate (\ref{ineg-Brinkman}) yields
\begin{align*}
\int_0 ^t \Vert F(s) \Vert_{\Ld^2(\Omega)}^2   \, \mathrm{d}s \lesssim_0 \E(0) \left( 1 + \mathrm{M}_{\rho_f,j_f}(t)\right),
\end{align*}
for all $t \in [1,t^{\star}]$, so that Lemma \ref{tstar>1}, Remark \ref{remark:M} and the Notation (\ref{Lesssim0}) for the symbol $\lesssim_0$ provide the existence of a nondecreasing function $\varphi : \R^+ \rightarrow \R^+$ which is onto and such that
\begin{align*}
\int_0 ^{t^{\star}} \Vert F(s) \Vert_{\Ld^2(\Omega)}^2   \, \mathrm{d}s \leq \E(0) \varphi \Big( 1+\Vert \nabla u_0 \Vert_{\Ld^2(\Omega)}+ N_q f_0  + \mathrm{E}(0) \Big).
\end{align*}
We thus infer that
\begin{align*}
\Vert \nabla u_0 \Vert_{\Ld^2(\Omega)}^2 & + C_1 \int_0^{t^{\star}} \Vert F(s) \Vert_{\Ld^2(\Omega)}^2 \, \mathrm{d}s  \\
 & \leq \Vert \nabla u_0 \Vert_{\Ld^2(\Omega)}^2 +  C_1 \E(0) \varphi \Big( 1+\Vert \nabla u_0 \Vert_{\Ld^2(\Omega)}+ N_q f_0  + \mathrm{E}(0) \Big) \\
 & \leq \dfrac{1}{2\sqrt{8C_1 C_2}} +C_1 \E(0) \varphi \Big( 1+\Vert \nabla u_0 \Vert_{\Ld^2(\Omega)}+ N_q f_0  + \mathrm{E}(0) \Big).
\end{align*}
Thanks to the smallness condition (\ref{smallness:condition}), we can choose $\E(0)$ and $\Vert \nabla u_0 \Vert_{\Ld^2(\Omega)}$ small enough so that
\begin{align*}
\Vert \nabla u_0 \Vert_{\Ld^2(\Omega)}^2 & + C_1 \int_0^{t^{\star}} \Vert F(s) \Vert_{\Ld^2(\Omega)}^2 \, \mathrm{d}s< \dfrac{1}{\sqrt{8C_1 C_2}}.
\end{align*}
Using again the integrability (\ref{ineg-Brinkman}) of $s \mapsto \Vert F(s) \Vert_{\Ld^2(\Omega)}^2$ on $[0,T]$ for all $T>0$, we obtain by continuity a strong existence time (strictly) larger than $t^{\star}$.

\textbullet \ \ Now, we turn to the existence of a strong existence time larger than $t^{\star}$ which satisfies (\ref{def:tstar}). We use the uniform control of Proposition \ref{grad:unif} to get the existence of a nondecreasing continuous and onto function $\varphi : \R^+ \rightarrow \R^+$ such that 
\begin{align*}
\int_1^{t^{\star}} \Vert \nabla u(s) \Vert_{\Ld^{\infty}(\Omega)} \leq\varphi \Big( 1+\Vert \nabla u_0 \Vert_{\Ld^2(\Omega)}+ N_q f_0   + \mathrm{E}(0) \Big) \E(0)^{\gamma},
\end{align*}
with $\gamma >0$. A smallness condition such as (\ref{smallness:condition})
again ensures that 
\begin{align*}
\int_1^{t^{\star}} \Vert \nabla u(s) \Vert_{\Ld^{\infty}(\Omega)} \leq \dfrac{\delta}{2}.
\end{align*}
Thanks to Lemma \ref{grad:nonunif},
a continuity argument shows that there exists a strong existence time $t>t^{\star}$ such that
\begin{align*}
\int_1^{t} \Vert \nabla u(s) \Vert_{\Ld^{\infty}(\Omega)} <\delta.
\end{align*}
This is a contradiction with the definition of $t^{\star}$. Therefore we must have $t^{\star}=+\infty$ and the proof of Theorem \ref{theoreme1} is finally complete.
\end{proof}

\section{Further description of the asymptotic local density}\label{Section:asympprofile}
This section aims at providing a further description of the asymptotic behavior of $f$ in the space variable, that is to say at proving Theorem \ref{thm-infini}. Indeed, we have obtained the existence of a spatial profile $\rho^{\infty}$ in an abstract framework in Corollary \ref{coro:result}. However, a careful study of the particle trajectory will bring more information about this asymptotic state. We will see that, in full generality, the description of this profile will depend on the whole evolution of the system.

\bigskip

We first refer to Section \ref{Section:particle:traj} where we have defined 
the characteristic curves $s \mapsto (\mathrm{X}(s), \mathrm{V}(s))$ in $\R^3 \times \R^3$ for the Vlasov equation (associated to the natural extension $Pu$ for the fluid velocity $u$). We will make a constant use of the notations used in this section.

\bigskip

Suppose that $u \in \Ld^1_{\mathrm{loc}}(\R^+;\H^1_0 \cap \W^{1,\infty}(\Omega))$. The characteristic curves for the Vlasov equation are given for all $t,s \geq 0$ by
\begin{equation}\label{expr:Zt}
\left\{
      \begin{aligned}
        \mathrm{X}(s;t,x,v)&=x+(1-e^{-s+t})v+\int_t^s (1-e^{\tau-s}) (Pu)(\tau,\mathrm{X}(\tau;t,x,v))  \, \mathrm{d}\tau,  \\
        \mathrm{V}(s;t,x,v)&= e^{-s+t}v+\int_t ^s e^{\tau-s}(Pu)(\tau,\mathrm{X}(\tau;t,x,v))   \, \mathrm{d}\tau.   \\
      \end{aligned}
    \right.
\end{equation}
We will also use the notation $\mathrm{Z}_{s,t}$ defined in (\ref{notation:Z}).

We then recall the following representation formula for a given weak solution $(u,f)$ to the Vlasov-Navier-Stokes system (see Proposition \ref{Prop:formulerep}). We have
\begin{align}\label{formule-rep2}
f(t,x,v)=e^{3t} \mathbf{1}_{\mathcal{O}^t}(x,v) \, f_0(\mathrm{X}(0;t,x,v),\mathrm{V}(0;t,x,v)),
\end{align}
where $\mathcal{O}^t$ is defined in (\ref{equivOt}).

\medskip

In order to describe the asymptotic profile $\rho^{\infty}$ of Corollary \ref{coro:result}, we take a test function $\psi\in \mathscr{C}^{\infty}_c(\Omega)$ and we look at the following quantity for $t \geq 0$
\begin{align*}
\int_{\Omega} \rho_f(t,x)\psi(x)   \, \mathrm{d}x.
\end{align*}
Thanks to the previous representation formula (\ref{formule-rep2}) for the Vlasov equation, we can write for all $t>0$
\begin{align*}
 \int_{\Omega} \rho_f(t,x)\psi(x)   \, \mathrm{d}x &= \int_{\Omega \times \R^3} f(t,x,v)\psi(x)  \, \mathrm{d}x \,\mathrm{d}v \\
 & = \int_{\Omega \times \R^3} e^{3t} \mathbf{1}_{\mathcal{O}^t } (x,v) \, f_0(\mathrm{X}(0;t,x,v),\mathrm{V}(0;t,x,v)) \psi(x)  \,\mathrm{d}x \,\mathrm{d}v.
\end{align*}
We now use the natural change of variable $z=\mathrm{Z}_{0,t}(x,v)$, by remembering that  $\mathrm{Z}_{0,t}(\mathcal{O}^t)=\left\lbrace  (x,v) \in \Omega \times \R^3 \ \mid \ \tau^{+}(0,x,v) >t  \right\rbrace$, where the forward exit time $\tau^{+}(0,x,v)$ is defined by
\begin{align*}
\tau^{+}(0,x,v)=\sup \left\lbrace s \geq 0\ \mid \forall \sigma \in [0,s], \  \mathrm{X}(\sigma;0,x,v) \in \Omega \right \rbrace,
\end{align*}
(see Section \ref{Annexe:proofrep} of the Appendix). We infer that for all $t>0$
\begin{align*}
 \int_{\Omega} \rho_f(t,x)\psi(x)   \, \mathrm{d}x =\int_{\Omega \times \R^3} \mathbf{1}_{\tau^{+}(0,x,v)>t} \, f_0(x,v) \psi(\mathrm{X}_{t,0}(x,v)) \, \mathrm{d}x \, \mathrm{d}v,
\end{align*}
as in the proof of Proposition \ref{Prop:formulerep} in Section \ref{Annexe:proofrep}. We now remark that for any $t>0$, $\mathbf{1}_{\tau^{+}(0,x,v)>t}=0$ if $x \notin \Omega$. We thus get for all $t>0$
\begin{align}\label{rhof:test} 
 \int_{\Omega} \rho_f(t,x)\psi(x)   \, \mathrm{d}x =\int_{\R^3 \times \R^3} \mathbf{1}_{\tau^{+}(0,x,v)>t} f_0(x,v) \psi(\mathrm{X}_{t,0}(x,v)) \, \mathrm{d}x \, \mathrm{d}v.
\end{align}

To go further, we need to understand the behavior of the family of curves $\mathrm{X}_{t,0}(x,v)$ when $t \rightarrow \infty$. This is the aim of the following Lemma, which is in the same spirit as \cite[Lemma 8.1]{HKMM}. 
%
%

\begin{lemme}\label{curve:asymp}
There exists $\delta >0$ such that if $u \in \Ld^1(\R^+; \H^1_0 \cap\mathrm{W}^{1,\infty}(\Omega))$ with
\begin{align}\label{borne-gradient-asymp}
\int_0 ^{\infty} \Vert \nabla u(s) \Vert_{\Ld^{\infty}(\Omega)}  \, \mathrm{d}s \leq \delta,
\end{align}
then the family of functions $(x,v) \mapsto \mathrm{X}_{t,0}(x,v)$ converges in $\mathscr{C}^1(\R^3 \times \R^3)$, when $t \rightarrow + \infty$, towards $\mathrm{X}_{\infty} : \R^3 \times \R^3 \mapsto \mathrm{X}_{\infty}(x,v) \in \R^3$ with
\begin{align*}
\mathrm{X}_{\infty}(x,v) =x+v+\int_0^{\infty}(Pu)(\tau,\mathrm{X}_{\tau,0}(x,v))  \, \mathrm{d}\tau.
\end{align*}
\end{lemme}
\begin{proof}

In order to show that the family of curves $\left( \mathrm{X}_{t,0}(x,v) \right)_t$ has a limit in $ \mathscr{C}^1(\R^3 \times \R^3)$ when $t \rightarrow + \infty$, we will show that is satisfies a local Cauchy's criterion when $t \rightarrow + \infty$.

By coming back to the expression (\ref{expr:Zt}), we first write
\begin{align}\label{expr:Zt2}
\mathrm{X}_{t,0}(x,v)=x+v+\int_0^{\infty} \mathbf{1}_{\tau \leq t} (Pu)(\tau,\mathrm{X}_{\tau,0}(x,v))  \, \mathrm{d}\tau+ \varepsilon(t,x,v),
\end{align}
where
\begin{align*}
\varepsilon(t,x,v):=e^{-t}v -\int_{0}^{\infty} e^{\tau-t} \mathbf{1}_{\tau \leq t} (Pu)(\tau,\mathrm{X}_{\tau,0}(x,v))   \, \mathrm{d}\tau.
\end{align*}
Henceforth, we will use the notation $\mathrm{X}_{\tau,0}(z):=\mathrm{X}_{\tau,0}(x,v)$, where $z=(x,v) \times \R^3 \times \R^3$ as usual.

To show that a Cauchy's criterion is satisfied, we compute the difference between two expressions (\ref{expr:Zt2}) at $t_1$ and $t_2$ with $0<t_1 <t_2$ and we get 
\begin{align*}
\vert \mathrm{X}_{t_2,0}(z)-\mathrm{X}_{t_1,0}(z)\vert &\leq  \int_{t_1}^{t_2}  \vert (P u)(\tau,\mathrm{X}_{\tau,0}(z)) \vert   \, \mathrm{d}\tau+ \vert \varepsilon(t_2,x,v) -\varepsilon(t_1,x,v) \vert \\
&\leq  \int_{t_1}^{t_2}  \Vert (Pu)(\tau)\Vert_{\Ld^{\infty}(\R^3)}    \, \mathrm{d}\tau + \vert \varepsilon(t_2,x,v) -\varepsilon(t_1,x,v) \vert\\
&\leq \int_{t_1}^{t_2}  \Vert u(\tau)\Vert_{\Ld^{\infty}(\Omega)}    \, \mathrm{d}\tau + \vert \varepsilon(t_2,x,v) -\varepsilon(t_1,x,v) \vert.
\end{align*}

We now fix a compact set $ K \subset \R^3  \times \R^3$. The previous computation thus yields
\begin{align*}
 \underset{z \in K}{\sup} \,\vert \mathrm{X}_{t_2,0}(z)-\mathrm{X}_{t_1,0}(z)\vert &\leq  \int_{t_1}^{t_2}  \Vert u(\tau)\Vert_{\Ld^{\infty}(\Omega)}    \, \mathrm{d}\tau+\underset{z \in K}{\sup} \, \vert \varepsilon(t_2,z) -\varepsilon(t_1,z) \vert \\
 &:=(\mathrm{I}) + (\mathrm{II}),
\end{align*}
and it remains to see how the terms $(\mathrm{I})$ and $(\mathrm{II})$ behave when $\min(t_1,t_2) \rightarrow + \infty$.

Thanks to Poincaré inequality on $\H^1_0(\Omega) \cap \W^{1,\infty}(\Omega)$ and to the assumption (\ref{borne-gradient-asymp}), we observe that $t \mapsto \Vert u(t) \Vert_{\Ld^{\infty}(\Omega)} $ is integrable on $\R^+$. We thus infer that the term $(\mathrm{I})$ converges to $0$ when $\min(t_1,t_2) \rightarrow + \infty$.

For the term $(\mathrm{II})$, we remark that for all $t>0$
\begin{align*}
\underset{z \in K}{\sup} \,\vert \varepsilon(t,z) \vert  \lesssim_K e^{-t} +\int_{0}^{\infty} e^{\tau-t} \mathbf{1}_{\tau \leq t}  \Vert (P u)(\tau) \Vert_{\Ld^{\infty}(\R^3)}    \, \mathrm{d}\tau.
\end{align*}
Since $e^{\tau-t} \mathbf{1}_{\tau \leq t}  \Vert (P u)(\tau) \Vert_{\Ld^{\infty}(\R^3)}   \overset{t \rightarrow + \infty}{\longrightarrow} 0$ for almost any $\tau \in \R^+$ and since it is controlled by the integrable function $\tau \mapsto \Vert u(\tau) \Vert_{\Ld^{\infty}(\Omega)}$ on $\R^+$, we can apply the dominated convergence theorem to show that the previous r.h.s converges to $0$ when $t \rightarrow + \infty$.
Thanks to the triangular inequality, the term $(\mathrm{II})$ finally tends to $0$ when $\min(t_1,t_2) \rightarrow + \infty$. 

Therefore we have shown that
\begin{align*}
\underset{z \in K}{\sup} \, \vert \mathrm{X}_{t_2,0}(z)-\mathrm{X}_{t_1,0}(z) \vert\longrightarrow 0 \ \ \text{when} \ \ \min(t_1,t_2) \longrightarrow + \infty.
\end{align*}
This yields the existence of a mapping $(s,z) \mapsto \mathrm{X}_{\infty}(x,v) \in \R^3$ which is the uniform limit in $\mathscr{C}(\R^3\times \R^3)$ of the maps $z \mapsto \mathrm{X}_{t,0}(z)$ quand $t \rightarrow + \infty$. Again thanks to the dominated convergence theorem, we can pass to the limit when $t \rightarrow + \infty$ in the expression (\ref{expr:Zt}) to get, as announced
\begin{align*}
\mathrm{X}_{\infty}(x,v) = x+v+\int_0^{\infty}(Pu)(\tau,\mathrm{X}_{\tau,0}(x,v))\, \mathrm{d}\tau.
\end{align*}

For the moment, the limit $\mathrm{X}_{\infty}(x,v)$ is only continuous in its variables $(x,v)$ as a uniform limit of such functions. A (local) Cauchy's criterion for the family of derivatives in space-velocity $z \mapsto \mathrm{D}_z   \mathrm{X}_{t,0}(z) $ will be actually satisfied when $t \rightarrow + \infty$. Indeed, we differentiate the expression (\ref{expr:Zt}) to get
\begin{align*}
\mathrm{D}_z  \mathrm{X}_{t,0}(z) = 2\mathrm{A}+\int_{0}^{\infty} \mathbf{1}_{ \tau \leq t}\nabla (Pu)(\tau,\mathrm{X}_{\tau,0}(z)) \mathrm{D}_z   \mathrm{X}_{\tau,0}(z)  \, \mathrm{d}\tau + \widetilde{\varepsilon}(t,z),
\end{align*}
where
\begin{align*}
\widetilde{\varepsilon}(t,z)&:=e^{-t}\mathrm{B}-\int_{0}^{\infty} e^{\tau-t} \mathbf{1}_{\tau \leq t} \nabla (Pu)(\tau,\mathrm{X}_{\tau,0}(z)) \mathrm{D}_z   \mathrm{X}_{\tau,0}(z)   \, \mathrm{d}\tau, \\[2mm]
\mathrm{A}&:=  \begin{pmatrix} \mathrm{I}_3  \mid \mathrm{I}_3  \end{pmatrix}
\in \mathrm{M}_{3,6}(\R), \\[2mm]
\mathrm{B}&:=  \begin{pmatrix} \mathrm{0}_{\mathrm{M}_{3,3}(\R)}  \mid \mathrm{I}_3  \end{pmatrix}\in \mathrm{M}_{3,6}(\R),
\end{align*}
this identity being valid thanks to the integrability of $\tau \mapsto \Vert \nabla(P u)(\tau)\Vert_{\Ld^{\infty}(\R^3)}$ on $\R^+$. As before, if $t_1 < t_2$, we have the following inequality
\begin{align*}
\vert \mathrm{D}_z  \mathrm{X}_{t_2,0}(z)- \mathrm{D}_z  \mathrm{X}_{t_1,0}(z)\vert &\leq \int_{t_1}^{t_2}  \vert \nabla (Pu)(\tau,\mathrm{X}_{\tau,0}(z)) \mathrm{D}_z   \mathrm{X}_{\tau,0}(z) \vert  \, \mathrm{d}\tau  + \vert \widetilde{\varepsilon}(t_2,z)- \widetilde{\varepsilon}(t_1,z) \vert.
\end{align*}
We thus need to derive some bounds on $\mathrm{D}_z   \mathrm{X}_{\tau,0}(z)$. We have for all $0 \leq s \leq t$
\begin{align*}
\vert \mathrm{D}_x   \mathrm{X}_{s,0}(z)\vert & \leq 1+\int_0^{\infty} \mathbf{1}_{\tau \leq t}(1-e^{\tau-t}) \vert \nabla (Pu)(\tau,\mathrm{X}_{\tau,0}(z) \vert \vert \mathrm{D}_x   \mathrm{X}_{\tau,0}(z)\vert  \, \mathrm{d}\tau \\
&  \leq 1+\underset{0 \leq \tau \leq t}{\sup}\vert \mathrm{D}_x   \mathrm{X}_{\tau,0}(z) \vert \int_0^{\infty}  \Vert \nabla (Pu)(\tau) \Vert_{\Ld^{\infty}(\R^3)}\, \mathrm{d}\tau \\
& \leq 1+\underset{0 \leq \tau \leq t}{\sup}\vert \mathrm{D}_x   \mathrm{X}_{\tau,0}(z) \vert \int_0^{\infty}  \Vert \nabla u(\tau) \Vert_{\Ld^{\infty}(\Omega)}\, \mathrm{d}\tau.
\end{align*}
Thanks to the assumption (\ref{borne-gradient-asymp}) on $\nabla u$, we obtain 
\begin{align}\label{bound:Dx}
\underset{0 \leq \tau \leq t}{\sup}\vert \mathrm{D}_x   \mathrm{X}_{\tau,0}(z) \vert \leq \dfrac{1}{1-\delta}.
\end{align}
We observe that the very same procedure leads to
\begin{align}\label{bound:Dv}
\underset{0 \leq \tau \leq t}{\sup}\vert \mathrm{D}_v   \mathrm{X}_{\tau,0}(z) \vert \leq \dfrac{1}{1-\delta}.
\end{align}
Imposing $\delta \leq 1/2$, the two previous r.h.s become smaller that $2$ and we will use this uniform bound in what follows.
Again, we fix a compact set $ K \subset \R^3 \times \R^3$. We have
\begin{align*}
 \underset{z \in K}{\sup} \,\vert \mathrm{D}_z \mathrm{X}_{t_2,0}(z)-\mathrm{D}_z\mathrm{X}_{t_1,0}(z)\vert 
 &\leq  \int_{t_1}^{t_2}  \Vert \nabla (Pu)(\tau)\Vert_{\Ld^{\infty}(\R^3)}    \, \mathrm{d}\tau+\underset{z \in K}{\sup} \, \vert \widetilde{\varepsilon}(t_2,z) -\widetilde{\varepsilon}(t_1,z) \vert \\
& \leq \int_{t_1}^{t_2}  \Vert \nabla u (\tau)\Vert_{\Ld^{\infty}(\Omega)}    \, \mathrm{d}\tau+\underset{z \in K}{\sup} \, \vert \widetilde{\varepsilon}(t_2,z) -\widetilde{\varepsilon}(t_1,z) \vert \\
 &:=(\widetilde{\mathrm{I}}) + (\widetilde{\mathrm{II}}).
\end{align*}

As before, the first term $(\widetilde{\mathrm{I}}) $ tends to $0$ when $\min(t_1,t_2)$ goes to $+ \infty$ because of the integrability assumption (\ref{borne-gradient-asymp}).

For the second term $(\widetilde{\mathrm{II}})$, we proceed as before and we write for all $t>0$
\begin{align*}
\underset{z \in K}{\sup} \,\vert \widetilde{\varepsilon}(t,z) \vert  \leq e^{-t} +\int_{0}^{\infty} e^{\tau-t} \mathbf{1}_{\tau \leq t}  \Vert \nabla (P u)(\tau) \Vert_{\Ld^{\infty}(\R^3)}    \, \mathrm{d}\tau.
\end{align*}
Again, the dominated convergence theorem shows that this expression tends to $0$ when $t$ tends to $+ \infty$. and this implies that the term $(\widetilde{\mathrm{II}})$ converges to $0$ when $\min(t_1,t_2)$ goes to $0$.

All in all, we have shown that
\begin{align*}
\underset{z \in K}{\sup} \, \vert \mathrm{D}_z \mathrm{X}_{t_2,0}(z)-\mathrm{D}_z \mathrm{X}_{t_1,0}(z) \vert\longrightarrow 0 \ \  \text{when} \ \  \min(t_1,t_2) \longrightarrow + \infty.
\end{align*}
Cauchy's criterion then applies for the derivatives and this concludes the proof.
\end{proof}

In what follows, we use the notation $\X_{t,0,v}(x):=\X(t;0,x,v)$ and $\X_{\infty,v}(x):=\X_{\infty}(x,v)$.
\begin{lemme}\label{LM:chvar-Xinfini}
Under the assumption (\ref{borne-gradient-asymp}), for all $v \in \R^3$ and $t \in \R^+$, the maps $\X_{t,0,v} : x \mapsto \X_{t,0,v}(x)$ and the map $\X_{\infty,v} : x \mapsto \X_{\infty}(x,v)$ are $\mathscr{C}^1$-diffeomorphisms from $\R^3$ to $\R^3$.
\end{lemme}
\begin{proof}
The proof follows the same steps as that of \cite[Lemma 8.2]{HKMM}. Indeed, from the expression (\ref{expr:Zt}), we get 
\begin{align*}
\mathrm{D}_x   \mathrm{X}_{t,0,v}(x)- \mathrm{I}_3=\int_0^{\infty} \mathbf{1}_{\tau \leq t}(1-e^{\tau-t})  \nabla (Pu)(\tau,\mathrm{X}_{\tau,0}(x,v))  \mathrm{D}_x   \mathrm{X}_{\tau,0}(x,v)   \, \mathrm{d}\tau,
\end{align*}
and then, thanks to the bound (\ref{bound:Dx}), we can write
\begin{align*}
\Vert \mathrm{D}_x   \mathrm{X}_{t,0,v}- \mathrm{I}_3 \Vert_{\infty} \leq \dfrac{\delta}{1-\delta},
\end{align*}
so that up to taking $\delta$ small enough in the assumption (\ref{borne-gradient-asymp}), we can obtain  
\begin{align}\label{bound:C1diff-t}
\Vert \mathrm{D}_x   \mathrm{X}_{t,0,v}- \mathrm{I}_3 \Vert_{\infty} \leq \dfrac{1}{9}.
\end{align}
A variant of the global inversion theorem about perturbation of the identity mapping (see \cite[Lemma 9.4]{HKMM}) leads to the conclusion.
\end{proof}

In order to pass to the limit in (\ref{rhof:test}) when $t\rightarrow + \infty$, we need to determine the limit of $\mathbf{1}_{\tau^{+}(0,x,v)>t}$ when $t \rightarrow +\infty$.
\begin{lemme}\label{LM:indicator}
For almost every $(x,v) \in \R^3 \times \R^3$, we have
\begin{align*}
\mathbf{1}_{\tau^{+}(0,x,v)>t} \underset{t \rightarrow +\infty}{\longrightarrow} \mathbf{1}_{\mathcal{O}_{\infty}}, 
\end{align*}
where 
$\displaystyle 
\mathcal{O}_{\infty}:=\left\lbrace (x,v) \in \Omega \times \R^3\mid \forall t \geq 0, \  \X(t;0,x,v) \in \Omega \right\rbrace.
$
\end{lemme}
\begin{proof}
First, we note that for every $(x,v) \in \Omega \times \R^3$, we have $\tau^{+}(0,x,v)>0$ while for every $(x,v) \in \Omega^c \times \R^3$, we have $\tau^{+}(0,x,v) = 0$ and $\X(0;0,x,v) \in \Omega^c$. We then focus on the following alternative when $(x,v) \in \Omega \times \R^3$.
If $\tau^{+}(0,x,v)=+\infty$, then $t \mapsto \mathbf{1}_{\tau^{+}(0,x,v)>t} $ is constant equal to $1$ for all $t>0$ therefore for all $\sigma \geq 0$, $\X(\sigma;0,x,v) \in \Omega$. 
If $\tau^{+}(0,x,v)\neq+\infty$, then $\mathbf{1}_{\tau^{+}(0,x,v)>t} \underset{t \rightarrow +\infty}{\longrightarrow} 0$ and we have $\X( \tau^{+}(0,x,v);0,x,v) \in \partial \Omega$ so that the conclusion holds.
\end{proof}

\medskip

\medskip

We can now conclude the proof of Theorem \ref{thm-infini} by passing to the limit in the expression (\ref{rhof:test}) when $t \rightarrow + \infty$. Indeed, thanks to Lemma \ref{curve:asymp} and \ref{LM:indicator} and since $f_0 \in \Ld^{1}(\Omega \times \R^3)$, we can use the dominated convergence theorem when $t \rightarrow + \infty$ to get
\begin{align*}
 \int_{\Omega} \rho_f(t,x)\psi(x)   \, \mathrm{d}x  \underset{t \rightarrow + \infty}{\longrightarrow} \int_{\R^3 \times \R^3} \mathbf{1}_{\mathcal{O}_{\infty}}(x,v) f_0(x,v) \psi(\X_{\infty,v}(x)) \, \mathrm{d}x \, \mathrm{d}v.
\end{align*}
We then use the reverse change of variable $x=\X^{-1}_{\infty,v}(y)$ for all velocity $v \in \R^3 $ in the previous integral thanks to Lemma \ref{LM:chvar-Xinfini} in order to obtain
\begin{align*}
&\int_{\R^3 \times \R^3} \mathbf{1}_{\mathcal{O}_{\infty}}(x,v) f_0(x,v) \psi(\X_{\infty}(x,v)) \, \mathrm{d}x \, \mathrm{d}v \\
&=\int_{\R^3 \times \R^3} \mathbf{1}_{\mathcal{O}_{\infty}}(\X^{-1}_{\infty,v}(x),v) f_0(\X^{-1}_{\infty,v}(x),v)  \vert \det \, \mathrm{D}_x \X^{-1}_{\infty,v}(x) \vert \psi(x) \, \mathrm{d}x \, \mathrm{d}v \\
&=\int_{\Omega} \int_{\R^3} \mathbf{1}_{\mathcal{O}_{\infty}}(\X^{-1}_{\infty,v}(x),v) f_0(\X^{-1}_{\infty,v}(x),v) \vert \det \, \mathrm{D}_x \X^{-1}_{\infty,v}(x) \vert \psi(x) \, \mathrm{d}x \, \mathrm{d}v,
\end{align*}
because $\psi$ is compactly supported in $\Omega$. Thanks to the convergence $\rho_f(t) \overset{t \rightarrow +\infty}{\longrightarrow} \rho^{\infty}$ in $\H^{-1}(\Omega)$ of Corollary \ref{coro:result}, we thus get by uniqueness of the limit that for a.e $x \in \Omega$,
\begin{align*}
\rho^{\infty}(x)=\int_{\R^3} \mathbf{1}_{\mathcal{O}_{\infty}}(\X^{-1}_{\infty,v}(x),v) f_0(\X^{-1}_{\infty,v}(x),v) \vert \det \, \mathrm{D}_x \X^{-1}_{\infty,v}(x) \vert \, \mathrm{d}v.
\end{align*}
Note that the meaning of the previous indicator function is the following. If $(x,v) \in \Omega \times \R^3$, 
\begin{align*}
(\X^{-1}_{\infty,v}(x),v) \in \mathcal{O}_{\infty} 
\Leftrightarrow \, & \X^{-1}_{\infty,v}(x) \in \Omega, \ \ \text{and} \ \  \forall t \geq 0, \ \  \X(t;0,\X^{-1}_{\infty,v}(x),v) \in \Omega\\
\Leftrightarrow \, & \exists ! y \in \Omega, \ \   \X_{\infty,v}(y)=x, \ \ \text{and} \ \  \forall t \geq 0, \ \  \X(t;0,y,v) \in \Omega.
\end{align*}
This concludes the proof of Theorem \ref{thm-infini}. 

\section{Asymptotic profiles with a prescribed mass}\label{section:prescribe:mass}
As already explained in the Introduction, the asymptotic profile $\rho^{\infty}$ has a total mass which is not known \textit{a priori}, because $t \mapsto \Vert \rho_f(t) \Vert_{\Ld^1(\Omega)}$ can be decreasing. We can actually find some class of initial data for which any prescribed mass (which is less than of equal to the initial mass $1$) will be achieved, that is Proposition \ref{prop:massalpha}.

 In order to do so, we first study two scenarios for the support of the initial density, namely  we impose $\mathrm{supp} \, f_0 \subset \mathrm{K}_1 \times \mathrm{K}_2$ where
\begin{itemize}
\item $\mathrm{K}_1$ is a compact of $\Omega$, at strictly positive distance from $\partial \Omega$, say included in $\B_x(a,\varepsilon) \subset \Omega$,
\item $\mathrm{K}_2$ is a compact of $\R^3$, say included in $\B_v(0,R)$ for small $R>0$, \textit{or} is an exterior domain of $\R^3$, say included in $(\R^3 \setminus \B_v(0,R) )$ for large $R>0$.
\end{itemize}  

Thanks to the representation formula (\ref{formule-rep}), we observe that for all $t \geq 0$, $f(t)$ has a support in space and velocity which is transported by the flow from the initial support of $f_0$.
More precisely, we have for all $t \geq 0$
\begin{align*}
\mathrm{supp} \, f(t) \subset\ \Big( \mathrm{X}(t;0,\mathrm{supp} \, f_0) \times  \mathrm{V}(t;0,\mathrm{supp} \, f_0) \Big)\cap \mathcal{O}^t \subset \R^3 \times \R^3.
\end{align*}
Our main idea is that small initial velocities will lead to a limit profile which is compactly supported in $\Omega$ while high initial velocities will give a profile vanishing in $\Omega$.

\medskip

In this section, we shall use several times the DiPerna-Lions theory for transport equation in the same fashion as that of the proof of Lemma \ref{inegdecal}: this allows us to define the characteric curves in a classical sense. Note that this procedure is actually only required on the interval of time $[0,1]$ because the smallness condition \eqref{smallness:condition} will ensure that $u \in \Ld^1(1,+\infty;\W^{1,\infty}_0(\Omega))$ (see the proof of Propositions \ref{prop:inside}-\ref{prop:outside}).

\subsection{The case of small initial velocities}
We first investigate the case in which the initial support in velocity is included in a small ball, leading to a situation where particles stay confined in the domain $\Omega$,  and far from the boundary.

\begin{proposition}\label{prop:inside}
Let $(u_0,f_0)$ be an admissible initial condition in the sense of Definition \ref{CIadmissible} satisfying 
\begin{align*}
\begin{split}
\mathrm{supp} \, f_0 \subset \B_x(a,\varepsilon) \times \B_v(0,R) \subset \Omega \times \R^3,
\end{split}
\end{align*}
where $\varepsilon, R>0$ satisfy the geometric condition
\begin{align}\label{geom:cond:thm}
& d(\overline{\B}_x(a,\varepsilon),\partial \Omega)>0, \ \ 2R<d(\overline{\B}_x(a,\varepsilon),\partial \Omega).
\end{align}
If the smallness condition (\ref{smallness:condition}) is satisfied, then, for any weak solution $(u,f)$ to the Vlasov-Navier-Stokes system with initial data $(u_0,f_0)$, there exists $d=d(\varepsilon,R,\Omega)>0$ such that we have for all $t \geq 0$
\begin{align*}
&\mathrm{supp} \,  f(t) \subset \B_x(a,d) \times \R^3  \subsetneq  \Omega \times \R^3.
\end{align*}
Furthermore, there exists $\rho^{\infty} \in \Ld^{\infty}_c(\Omega)$ with $\int_{\Omega} \rho^{\infty} (x) \, \mathrm{d}x=1$ and
such that 
\begin{align*}
\mathrm{W}_1\left( f(t), \rho^{\infty}\otimes \delta_{v=0} \right) \underset{t \rightarrow +\infty}{\longrightarrow}0,
\end{align*}
exponentially fast.
\end{proposition}

We state two lemmas that give estimates for the size of the image of the initial support $\B_x(a,\varepsilon) \times \B_v(0,R)$ under the flow $t \mapsto \mathrm{X}(t;0,x,v)$. This is achieved thanks to an appropriate control on the $\Ld^1\Ld^{\infty}$ norm of $u$.

\begin{lemme}\label{LM:geom1}
Let $u \in \Ld^2_{\mathrm{loc}}(\R^+;\H^1_0(\Omega)) \cap \Ld^1_{\mathrm{loc}}(\R^+;\W^{1,\infty}(\Omega))$. For all $t \geq 0$, for all $x \in \Omega$ and for all $v \in \B(0,R)$, we have 
\begin{align*}
\vert \mathrm{V}(t;0,x,v) -v \vert &\leq    \Vert  u \Vert_{\Ld^1(0,t;\Ld^{\infty}(\Omega))} +R,\\
\vert \mathrm{X}(t;0,x,v) -x \vert &\leq  2  \Vert   u \Vert_{\Ld^1(0,t;\Ld^{\infty}(\Omega))} +R.
\end{align*}
Moreover, we have
\begin{align*}
\mathrm{V}(t;0,x,v) \in \B\left(0,    \Vert   u \Vert_{\Ld^1(0,t;\Ld^{\infty}(\Omega))} +2R \right),
\end{align*}
and if $x\in \B(a,\varepsilon)$, we have 
\begin{align*}
\mathrm{X}(t;0,x,v) \in \B\left(a,  2  \Vert   u \Vert_{\Ld^1(0,t;\Ld^{\infty}(\Omega))} +R+\varepsilon\right).
\end{align*}
\end{lemme}
\begin{proof}
Recall that we have the formula for $t \geq 0$
\begin{align*}
        \mathrm{V}(t;0,x,v)&= e^{-t}v+\int_0^t e^{\tau-t}(P u)(\tau,\mathrm{X}(\tau;0,x,v)) \,  \mathrm{d}\tau,\\
        \mathrm{X}(t;0,x,v)&=x+v -\mathrm{V}(t;0,x,v)+\int_0^t  (P u)(\tau,\mathrm{X}(\tau;0,x,v)) \, \mathrm{d}\tau.
\end{align*}
Therefore, for $t \geq 0$, the triangular inequality first leads to
\begin{align*}
\vert \mathrm{V}(t;0,x,v) -v \vert &\leq \vert \mathrm{V}(t;0,x,v) -e^{-t}v \vert +\vert e^{-t}v -v \vert\\
&\leq \Vert P u \Vert_{\Ld^1(0,t;\Ld^{\infty}(\R^3))} + (1-e^{-t})\vert v \vert \\
&\leq \Vert  u \Vert_{\Ld^1(0,t;\Ld^{\infty}(\Omega))} +R,
\end{align*}
if $v \in \B_v(0,R)$. Then, in a similar way we get
\begin{align*}
\vert \mathrm{X}(t;0,x,v) -x \vert &\leq \Vert P u \Vert_{\Ld^1(0,t;\Ld^{\infty}(\R^3))}  +  \vert \mathrm{V}(t;0,x,v) -v \vert \\
&\leq \Vert P u \Vert_{\Ld^1(0,t;\Ld^{\infty}(\R^3))}  + \Vert  u \Vert_{\Ld^1(0,t;\Ld^{\infty}(\Omega))} +R\\
&\leq 2 \Vert u \Vert_{\Ld^1(0,t;\Ld^{\infty}(\Omega))} +R.
\end{align*}
Finally, the triangular inequality gives the last statements.
\end{proof}

\begin{corollaire}\label{coro:geom1}
Let $u \in \Ld^2_{\mathrm{loc}}(\R^+;\H^1_0(\Omega)) \cap \Ld^1_{\mathrm{loc}}(\R^+;\W^{1,\infty}(\Omega))$. Let $\varepsilon>0$ and $a \in \Omega$ such that $\B_x(a,\varepsilon) \subset \Omega$ with $d(\overline{\B}_x(a,\varepsilon),\partial \Omega)>0$ and $R>0$ such that
\begin{align}\label{cond:geom}
2R<d(\overline{\B}_x(a,\varepsilon),\partial \Omega).
\end{align}
Suppose that the velocity field $u$ satisfies
\begin{align}
\Vert u \Vert_{\Ld^1(\R^+;\Ld^{\infty}(\Omega))} \leq \delta, \label{cond:veloc-u}
\end{align}
where 
\begin{align*}
\delta:=\dfrac{1}{2}\left(\dfrac{d(\overline{\B}_x(a,\varepsilon),\partial \Omega)}{2}-R\right).
\end{align*}
Then for all $T \geq 0$ and for all $(x,v) \in \B_x(a,\varepsilon) \times \B_v(0,R)$, we have
\begin{align}\label{charac:stayint}
\mathrm{X}(T;0,x,v) \in \Omega,
\end{align}
and more precisely
\begin{align}\label{locball:charac}
\mathrm{X}(T;0,x,v) \in \B \Big( a,\varepsilon+R+2\delta \Big) \subsetneq \Omega.
\end{align}
\end{corollaire}
\begin{proof}
Let $T \geq 0$. We first apply Lemma \ref{LM:geom1} until time $T$ to see that for all $\sigma \in [0,T]$
$$\mathrm{X}(\sigma;0,x,v) \in \B\left(a, \varepsilon + R + 2\Vert   u \Vert_{\Ld^1(0,T;\Ld^{\infty}(\Omega))} \right),$$ 
if $(x,v) \in \B_x(a,\varepsilon) \times \B_v(0,R)$.

Then, thanks to the global assumption (\ref{cond:veloc-u}), we get for all $\sigma \in [0,T]$ and for all $(x,v) \in \B_x(a,\varepsilon) \times \B_v(0,R)$
$$\mathrm{X}(\sigma;0,x,v)  \in \B(a,L),
$$
where 
$$L:=\varepsilon+R+2\delta.$$
Thus, the point is just to ensure that 
\begin{align*}
\B(a,  L) \subsetneq \Omega.
\end{align*}
Since $d(a,\partial \Omega)=\varepsilon+d(\overline{\B}_x(a,\varepsilon),\partial \Omega)$,
the inequality $L<d(a,\partial \Omega)$ will be satisfied if we prove that
$$R+2\delta <d(\overline{\B}_x(a,\varepsilon),\partial \Omega).$$
By the definition of $\delta$, we have indeed 
$$R+2 \delta=\dfrac{d(\overline{\B}_x(a,\varepsilon),\partial \Omega)}{2}<d(\overline{\B}_x(a,\varepsilon),\partial \Omega),$$
and thus infer that
\begin{align*}
\B(a,L) \subsetneq  \Omega, \ \ \mathrm{with } \ \ d ( \B(a,  L ), \partial \Omega ) >0,
\end{align*}
leading to the conclusion (\ref{charac:stayint}).
\end{proof}
It means that, for a well chosen small support of $f_0$, the support of $f(t)$ stays inside $\Omega \times \R^3$ and far from the boundary as long as $ \Vert  u \Vert_{\Ld^1_T\Ld^{\infty}_x(\Omega)}$ is small enough for all positive times. This assumption will be actually essentially satisfied by using the previous results of the bootstrap argument of Section \ref{section:bootstrap}.

\bigskip

We eventually turn to the proof of Proposition \ref{prop:inside}. 
\begin{proof}\textbullet \ \ Recall that we work with an initial distribution $f_0$ such that $\mathrm{supp} \, f_0 \subset \B_x(a,\varepsilon) \times \B_v(0,R) \subset \Omega \times \R^3$ and under the assumption (\ref{geom:cond:thm}), namely
$$ 2R<d(\overline{\B}_x(a,\varepsilon),\partial \Omega)
,$$
which is exactly the previous condition (\ref{cond:geom}).

First, we define $\delta=\delta(\varepsilon,R,\Omega)>0$ as in (\ref{cond:veloc-u}). We then use the fact that $u \in \Ld^1(0,1,\Ld^{\infty}(\Omega))$ (see Proposition \ref{controlLinfiniLOC}) together with the dominated convergence theorem to obtain $\underline{t} \in ]0,1[$ such that the condition $$\Vert u \Vert_{\Ld^1(0,\underline{t};\Ld^{\infty}(\Omega))}<\dfrac{\delta}{2},$$ holds.

Furthermore, from the assumption (\ref{smallness:condition}) on the initial data, we can perform the same analysis as in the bootstrap argument of Section \ref{section:bootstrap} to get that the quantity $\Vert \nabla u \Vert_{\Ld^1(\underline{t} ,+\infty;\Ld^{\infty}(\Omega))}$ is as small as we want if we reduce $\E(0)$ and $\Vert \nabla u_0 \Vert_{\Ld^2(\Omega)}$. Thanks to the Poincaré inequality on $\H^1_0\cap \W^{1,\infty}(\Omega)$, we see that we are able to get the control
\begin{align*}
\Vert u \Vert_{\Ld^1(\underline{t},+\infty;\Ld^{\infty}(\Omega))} \leq C_{\Omega} \Vert \nabla u \Vert_{\Ld^1(\underline{t},+\infty;\Ld^{\infty}(\Omega))}<\dfrac{\delta}{2},
\end{align*}
(this is possible because we work on an interval of time far from $0$).
Combining these two pieces, we obtain the fact that the global condition (\ref{cond:veloc-u}) is satisfied.
Thanks to Corollary \ref{coro:geom1} together with the inclusion
\begin{align*}
\mathrm{supp} \, f(t) \subset \mathrm{X}(t;0,\mathrm{supp} \, f_0) \times  \mathrm{V}(t;0,\mathrm{supp} \, f_0),
\end{align*}
we get for all $t \in \R^{+}$
\begin{align*}
\mathrm{supp} \,  f(t) \subsetneq  \Omega \times \R^3.
\end{align*}
In particular, with (\ref{locball:charac}) of Corollary \ref{coro:geom1}, we have for all $t \in \R^{+}$
\begin{align*}
\mathrm{supp} \,  f(t) \subset  \B_x \Big( a,\varepsilon+R+2\delta \Big) \times \R^3. 
\end{align*}
\end{proof}

\begin{remarque}\label{conserv:mass}
In such a situation where the particles never leave the domain, we observe that the trace of $f$ vanishes for all times therefore there is conservation of the mass. Thus, it is now possible to prove 
a convergence of the type
\begin{align*}
\mathrm{W}_{1,\overline{\Omega} \times \R^3}\left( f(t), \rho^{\infty}\otimes \delta_{v=0} \right)
\leq \E(0)^{1/2} C_{\lambda} \exp(- \lambda t),
\end{align*} 
in Corollary \ref{coro:result}, where $\mathrm{W}_{1,\overline{\Omega} \times \R^3}$ stands for the $1$-Wasserstein distance on $\overline{\Omega} \times \R^3$, and with
\begin{align}\label{cvgence:densite:cpct}
\mathrm{W}_{1,\overline{\Omega} }\left( \rho_f(t), \rho^{\infty} \right) \overset{ \rightarrow + \infty}{\longrightarrow} 0,
\end{align} 
where $\mathrm{W}_{1,\overline{\Omega}}$ stands for the $1$-Wasserstein distance on $\overline{\Omega}$, as it was already observed in the case of the torus (see \cite{HKMM}). Indeed, thanks to the vanishing trace of $f$, we can now take smooth test functions on $\overline{\Omega}$, whose Lipschitz constant is less than $1$, in the proof of Proposition \ref{asymp-profil} because there are no boundary terms.

In particular, the previous convergence (\ref{cvgence:densite:cpct}) and the fact that the support of each $\rho_f(t)$ is uniformly included in a ball imply that the limit $\rho^{\infty}$ is also compactly supported in the same ball.
\end{remarque}

%

\subsection{The case of high initial velocities}
Here, we deal with the situation where all particles escape from the domain $\Omega$ after a finite time so that the kinetic distribution vanishes uniformly after this time. More precisely, we have the following Proposition.

%
\begin{proposition}\label{prop:outside}
 Let $(u_0,f_0)$ be an admissible initial condition in the sense of Definition \ref{CIadmissible} and take any weak solution $(u,f)$ to the Vlasov-Navier-Stokes system with initial data $(u_0,f_0)$. There exists $L=L(\Omega)$ for which if $\varepsilon,T,R>0$ are choosen such that
\begin{align}
& \mathrm{supp} \, f_0 \subset \B_x(a,\varepsilon) \times \Big( \R^3 \setminus \B_v(0,R) \Big) \subset \Omega \times \R^3, \label{cond:supportTHM4} \\
& L>\varepsilon, \ \ R > \dfrac{2L+\varepsilon}{1-e^{-T}} , \label{cond:velocTHM4}
\end{align}
then the following holds. If the smallness condition (\ref{smallness:condition}) is satisfied, then for any weak solution $(u,f)$ to the Vlasov-Navier-Stokes system with initial data $(u_0,f_0)$, we have
\begin{align}\label{result:densityvanish}
\forall t \geq T, \ \ f(t)=0 \ a.e.
\end{align}
\end{proposition}
We first show that the particle trajectory always enters in the complementary of any ball after a finite time, provided that the initial support in velocity contains only "high" velocities and that the $\Ld^1\Ld^{\infty}$ norm of the fluid is small enough.

\begin{lemme}\label{LM:exit1point}
Let $L,\varepsilon,R,T>0$ such that
\begin{align}
 L &>\varepsilon, \label{LMhighvitesse1} \\
 R &>\dfrac{\varepsilon+2L}{1-e^{-T}}.\label{LMhighvitesse2}
\end{align} 
If
\begin{align}\label{nabla:u-OUT}
\Vert   u \Vert_{\Ld^1(\R^+;\Ld^{\infty}(\Omega))}< \dfrac{L}{8},
\end{align}
then for all $(x,v) \in \B_x(a,\varepsilon) \times (\R^3 \setminus \B_v(0,R))$, there exists $T_{x,v} \in (0,T]$ such that
\begin{align*}
\mathrm{X}(T_{x,v};0,x,v) \in \R^3 \setminus \B(a,L).
\end{align*}
\end{lemme}
\begin{proof}
First, we write the formula (\ref{expr:Zt}) under the form
\begin{equation}\label{eq:X+V}
\left\{
      \begin{aligned}
        \mathrm{X}(s;t,x,v)+\mathrm{V}(s;t,x,v)&=x+v+\int_t^s(Pu)(\tau,\mathrm{X}(\tau;t,x,v))  \, \mathrm{d}\tau,  \\
        \mathrm{V}(s;t,x,v)&= e^{-s+t}v+\int_t ^s e^{\tau-s}(Pu)(\tau,\mathrm{X}(\tau;t,x,v))   \, \mathrm{d}\tau.   \\
      \end{aligned}
    \right.
\end{equation}
For all $(x,v) \in \B_x(a,\varepsilon) \times (\R^3 \setminus \B_v(0,R))$ and $t \in [0,T]$, we use (\ref{eq:X+V}) and the triangular inequality to write

\begin{align*}
\vert \mathrm{X}(t;0,x,v) -a \vert
&\geq \left\vert (e^{-t}-1)v+\int_0^t  e^{\tau-t}(P u)(\tau,\mathrm{X}(\tau;0,x,v))\mathrm{d}\tau \right \vert  
 \\  
 & \quad -\left\vert x-a \right\vert -  \left\vert \int_0^t  (P u)(\tau,\mathrm{X}(\tau;0,x,v))\mathrm{d}\tau \right\vert ,
\end{align*}
so that 
\begin{align*}
\vert \mathrm{X}(t;0,x,v) -a \vert &\geq (1-e^{-t})R-\varepsilon -2 \int_0^t  \vert (P u)(\tau,\mathrm{X}(\tau;0,x,v)) \vert \mathrm{d}\tau \\
& \geq (1-e^{-t})R-\varepsilon -2\Vert   u \Vert_{\Ld^1(\R^+;\Ld^{\infty}(\Omega))}.
\end{align*}
We thus end up with
\begin{align*}
\underset{\tau \in [0,T]}{\sup} \ \vert \mathrm{X}(\tau;0,x,v)-a \vert \geq (1-e^{-T})R-\varepsilon -2\Vert   u \Vert_{\Ld^1(\R^+;\Ld^{\infty}(\Omega))}.
\end{align*}
Since we have $(1-e^{-T})R -\varepsilon>2L$ thanks to assumption (\ref{LMhighvitesse2}), we observe that the condition (\ref{nabla:u-OUT}) implies 
\begin{align*}
\forall (x,v) \in \B_x(a,\varepsilon) \times (\R^3 \setminus \B_v(0,R)), \ \ \underset{\tau \in [0,T]}{\sup} \ \vert \mathrm{X}(\tau;0,x,v)-a \vert >L.
\end{align*}
It means that that for all $(x,v) \in \B_x(a,\varepsilon) \times (\R^3 \setminus \B_v(0,R))$, there exists at least one time $T_{x,v} \geq 0$ such that 
\begin{align*}
\vert \mathrm{X}(T_{x,v};0,x,v)-a \vert >L.
\end{align*}
We observe that $T_{x,v}>0$ since $L>\varepsilon$ by assumption (\ref{LMhighvitesse1}) , and this concludes the proof.
\end{proof}

We now give a proof of Proposition \ref{prop:outside}.
\begin{proof}
Let us recall that we work in a small data regime (\ref{smallness:condition}) and under the assumption of high initial velocity (\ref{cond:supportTHM4})-(\ref{cond:velocTHM4}), that is 
\begin{align*}
& R> \dfrac{\varepsilon+2L}{1-e^{-T}},\\[1mm]
\mathrm{supp} \, f_0 & \subset \B_x(a,\varepsilon) \times (\R^3 \setminus \B_v(0,R)) \subset \Omega \times \R^3,
\end{align*}
for some $L$ to be determined later. In short, this condition will allow us to ensure that all the particles have left the domain juste after time $T$ (which is fixed). If we want this time $T$ to be small, it imposes high initial speeds ($L$ being fixed).

Thanks to the representation formula for the distribution function and by using the same change of variable as in the proof of Proposition \ref{Prop:formulerep}, we have for any $t \geq T$
\begin{align*}
\int_{\Omega \times \R^3} f(t,x,v) \, \mathrm{d}x \mathrm{d}v &= \int_{\mathcal{O}^t} e^{3t} f_0(\mathrm{Z}_{0,t}(x,v)) \, \mathrm{d}x  \, \mathrm{d}v \\
&= \int_{\Omega \times \R^3} \mathbf{1}_{\tau^+(0,z)>t}  \, f_0(z) \, \mathrm{d}z \\
& \leq \int_{\mathrm{supp} \, f_0} \mathbf{1}_{\tau^+(0,z)>T} \, f_0(z) \, \mathrm{d}z,
\end{align*}
because $f_0 \geq 0 \ a.e.$ Since $f \geq 0 \ a.e$, it is enough to show that the previous time $T>0$ (depending on $\mathrm{supp} \,f_0$ and $\Omega$) satisfies
$$  \ \forall z \in \mathrm{supp} \, f_0, \ \ \tau^+(0,z) \leq T,$$
for the integral $\int_{\Omega \times \R^3} f(t,x,v) \, \mathrm{d}x \mathrm{d}v$ to vanish for all $t \geq T$.

%
Recall that because of the definition of the extension operator $P$ and because $\Omega$ is bounded, we have
\begin{align*}
\mathrm{supp} \, Pu \subset\Omega  \subsetneq \B(a,L),
\end{align*}
for some $L=L(\Omega)>\varepsilon$.

First, we use the fact that $u \in \Ld^1_{\mathrm{loc}}(\R^+,\Ld^{\infty}(\Omega))$ (see Proposition \ref{controlLinfiniLOC}) together with the dominated convergence theorem to obtain a time $\overline{t} \in ]0,T[$ such that $\Vert u \Vert_{\Ld^1(0,\overline{t};\Ld^{\infty}(\Omega))}  < L/16$. Then, if we reduce $\E(0)$ and $\Vert \nabla u_0 \Vert_{\Ld^2(\Omega)}$, we can perform the same analysis as in the bootstrap argument of Section \ref{section:bootstrap} to get that the quantity $\Vert \nabla u\Vert_{\Ld^1(\underline{t} ,+\infty;\Ld^{\infty}(\Omega))}$ is as small as we want. Combining this argument with the Poincaré inequality on $\H^1_0\cap \W^{1,\infty}(\Omega)$, we get 
\begin{align*}
\Vert u \Vert_{\Ld^1(\underline{t},+\infty;\Ld^{\infty}(\Omega))} \leq C_{\Omega} \Vert \nabla u \Vert_{\Ld^1(\underline{t},+\infty;\Ld^{\infty}(\Omega))}<\dfrac{L}{16}.
\end{align*}
Therefore, the condition (\ref{nabla:u-OUT}) is indeed satisfied.

Take now $(x,v) \in \mathrm{supp} \,f_0$. 
We use Lemma \ref{LM:exit1point} with our choice of $L$ to obtain a finite time $t=t(x,v)>0$ such that
\begin{align*}
t \leq T, \ \ \mathrm{X}(t;0,x,v) \in \R^3 \setminus \B(a,L).
\end{align*} 
We emphasize the fact that $T$ does not depend on $(x,v)$. By definition of $\tau^+$ and by continuity in time of the trajectory, we have $\tau^+(0,x,v) \leq t(x,v)$ so that the control 
$$  \ \forall (x,v) \in \mathrm{supp} \, f_0, \  \tau^+(0,x,v) \leq T,$$
is indeed satisfied. It therefore implies that $f \equiv 0$ almost everywhere after time $T$, which is what we wanted to prove.
\end{proof}
\begin{remarque}
In such a scenario, we have $\rho^{\infty}=0$ in the statement of Corollary \ref{coro:result}.
\end{remarque}

\subsection{Proof of Proposition \ref{prop:massalpha}}
The aim of this subsection is to provide a proof of Proposition \ref{prop:massalpha}. We want to show that, for a fixed $\alpha \in [0,1]$, there exists an initial data $(u_0,f_0)$ which gives rise to weak solutions whose kinetic part concentrates in velocity, with a spatial asymptotic profile of total mass $\alpha$.

To do so, in view of the two previous sections, we consider $L=L(\Omega)$ given in Proposition \ref{prop:outside} and we choose $\varepsilon, T, R_1, R_2>0$ and initial admissible data in the sense of (\ref{data:hyp}) such that 
\begin{align*}
& R_1<R_2,\\
& \mathrm{supp} \, f_0 \subset \B_x(a,\varepsilon) \times\Big( \B_v(0,R_1) \sqcup   \left( \R^3 \setminus \B_v(0,R_2)  \right) \Big),  \\
&\int_{\Omega \times \R^3} f_0 \mathbf{1}_{\vert v \vert <R_1} \, \mathrm{d}x \, \mathrm{d}v =\alpha, \\[2mm]
& \B_x(a,\varepsilon) \subset \Omega, \ \ d(\overline{\B}_x(a,\varepsilon),\partial \Omega)>0, \ \ 2R_1<d(\overline{\B}_x(a,\varepsilon),\partial \Omega), \\[2mm]
& L>\varepsilon,  \ \ R_2 > \dfrac{2L+\varepsilon}{1-e^{-T}}.
\end{align*}
Furthermore, we consider that the initial kinetic energy $\E(0)$ and $\Vert \nabla u_0 \Vert_{\Ld^2(\Omega))}$ are small enough in the sense of (\ref{smallness:condition}).
%
%

We then take $(u,f)$ a weak solution to the Vlasov-Navier-Stokes system starting at $(u_0,f_0)$. We write
\begin{align*}
f_0=f_0 \mathbf{1}_{\vert v \vert < R_1}+f_0 \mathbf{1}_{\vert v \vert \geq R_2}:=f^{int}_0+f^{ext}_0,
\end{align*}
and we denote by $f^{int}$ (resp. $f^{ext}$) the renormalized solution to the Vlasov equation associated to the field $u$ and starting at $f^{int}_0$ (resp. starting at $f^{ext}_0$). Thanks to the well-posedness and the linearity of the Vlasov equation (for a fixed field $u$), we get
$$
f=f^{int}+f^{ext}.
$$
From Propositions \ref{prop:inside} and \ref{prop:outside}, we get
\begin{align*}
& \exists c >0, \ \ \forall t \geq 0, \ \ d(\mathrm{supp} \, f^{int}(t),\partial \Omega) \geq c>0,\\[2mm]
& \exists T:=T_{\varepsilon,R_2}>0, \ \ \forall t \geq T, \ \ f^{ext}(t)=0.
\end{align*}
Furthermore, since $f(t)=f^{int}(t)$ after time $t \geq T$, we can apply Remark \ref{conserv:mass} so that there exists $\rho^{\infty} \in \Ld^{\infty}(\Omega)$ compactly supported in $\Omega$ such that
$$ \W_{1,\overline{\Omega}}(\rho_f(t), \rho^{\infty}) \underset{t\rightarrow + \infty}{\longrightarrow} 0.$$
In particular, this weak convergence yields
\begin{align*}
\int_{\Omega} \rho_f(t) \, \mathrm{d}x  \underset{t\rightarrow + \infty}{\longrightarrow}\int_{\Omega} \rho^{\infty} \, \mathrm{d}x.
\end{align*}
But since we have conservation of the mass for the interior part $f^{int}$, we get for all $t \geq T$, omitting the variables,
\begin{align*}
\int_{\Omega} \rho_f(t) \, \mathrm{d}x =\int_{\Omega \times \R^3} f(t)\, \mathrm{d}x \, \mathrm{d}v=\int_{\Omega \times \R^3} f^{int}(t)\, \mathrm{d}x \, \mathrm{d}v=\int_{\Omega \times \R^3} f^{int}_0(t) \, \mathrm{d}x \, \mathrm{d}v=\alpha,
\end{align*}
from which we infer that the asymptotic profile $\rho^{\infty}$ satisfies
$$\int_{\Omega} \rho^{\infty} \, \mathrm{d}x=\alpha.$$
The proof of Proposition \ref{prop:massalpha} is therefore complete.
\begin{appendix}
\section{Appendix}

\subsection{Boundary value problem in $\Omega \times \R^3$ for the kinetic equation }\label{DiPernaLions}
\begin{theoreme}
Take $f_0 \in \Ld^1 \cap \Ld^{\infty}(\Omega \times \R^3)$ and a vector field $u \in \Ld^1_{\mathrm{loc}}(\R^+; \W^{1,1}(\Omega))$. Consider the following kinetic boundary value problem on $\Omega \times \R^3$.
\begin{align*}
\partial_t f +v\cdot \nabla_x f + {\rm div}_v((u-v)f)&=0,\\
f_{\mid t=0}&=f_0,\\
f &=0, \ \mathrm{on} \ \Sigma^{-}.
\end{align*}
Then we have, for all fixed $T>0$ 

\medskip

\textbullet \ \ \underline{Well-posedness}: There exists a unique 
$f \in \Ld^{\infty}_{\mathrm{loc}}(\R^+;\Ld^1 \cap\Ld^{\infty}(\Omega \times \R^3)) $
which is a weak solution of the previous Cauchy problem. Furthermore, $$ f \in \mathscr{C}(\R^+;\Ld^p_{\mathrm{loc}}(\overline{\Omega} \times \R^3)),$$ for all $p \in [1,\infty)$ and the function $f$ has a trace on $\partial \Omega \times \R^3$ defined in the following sense: there exists a unique element $\gamma f \in \Ld^{\infty}([0,T] \times \partial \Omega \times \R^3 )$ such that for any test function $\psi \in \mathscr{C}^{\infty}([0,T] \times \overline{\Omega } \times \R^3)$ with compact support in velocity, and for all $0 \leq t_1 \leq t_2 \leq T$
\begin{multline*}
\int_{t_1}^{t_2} \int_{\Omega \times \R^3}  f(t,x,v) \left[  \partial_t \psi + v \cdot \nabla_x \psi +(u-v) \cdot \nabla_v \psi \right](t,x,v) \, \mathrm{d}v \, \mathrm{d} x \, \mathrm{d}t \\ 
=\int_{\Omega \times \R^3} f(t_2,x,v) \psi(t_2,x,v) \, \mathrm{d}v \,\mathrm{d}x-\int_{\Omega \times \R^3} f(t_1,x,v) \psi(t_1,x,v) \, \mathrm{d}v \,\mathrm{d}x \\
+ \int_{t_1}^{t_2}\int_{\partial \Omega \times \R^3} \left[ (\gamma f) \psi(t,x,v) \right] v \cdot n(x) \, \mathrm{d}v \, \mathrm{d\sigma}(x) \, \mathrm{d}t.
\end{multline*}

\medskip

\textbullet \ \ \underline{Renormalization}: For every $\beta \in \mathscr{C}^1(\R)$, for all test function $\psi \in \mathscr{C}^{\infty}([0,T] \times \overline{\Omega } \times \R^3)$ with compact support in velocity, and for all $0 \leq t_1 \leq t_2 \leq T$, we have
\begin{multline*}
\int_{t_1}^{t_2} \int_{\Omega \times \R^3}  \beta(f(t,x,v)) \left[  \partial_t \psi + v \cdot \nabla_x \psi +(u-v) \cdot \nabla_v \psi \right](t,x,v) \, \mathrm{d}v \, \mathrm{d} x \, \mathrm{d}t \\ 
=\int_{\Omega \times \R^3} \beta(f(t_2,x,v)) \psi(t_2,x,v) \, \mathrm{d}v \,\mathrm{d}x-\int_{\Omega \times \R^3} \beta(f(t_1,x,v)) \psi(t_1,x,v) \, \mathrm{d}v \,\mathrm{d}x \\
+ \int_{t_1}^{t_2}\int_{\partial \Omega \times \R^3} \left[ \beta(\gamma f) \psi(t,x,v) \right] v \cdot n(x) \, \mathrm{d}v \, \mathrm{d\sigma}(x) \, \mathrm{d}t \\
-3 \int_{t_1}^{t_2} \int_{\Omega \times \R^3} \psi \left[f\beta'(f)-\beta(f) \right](t,x,v)\, \mathrm{d}v \, \mathrm{d} x \, \mathrm{d}t .
\end{multline*}
\medskip

\textbullet \ \ \underline{Stability}: If
\begin{align*}
u_n \rightarrow u  \ \ \text{in} \ \ \Ld^1_{\mathrm{loc}}(\R^+;\Ld^1(\Omega)) \ \text{ and } \
f_{0,n} \rightarrow f_0  \ \ \text{in} \ \ \Ld^1_{\mathrm{loc}}(\Omega \times \R^3),
\end{align*}
the corresponding sequence of solutions $(f_n)_n$ satisfies for all $p<\infty$,
\begin{align*}
f_n \longrightarrow f  \ \ \text{in} \ \Ld^{\infty}_{\mathrm{loc}}(\R^+;\Ld^p(\Omega \times \R^3)).
\end{align*}
\end{theoreme}
Such a result can be found in \cite[Theorem 3.2 - Proposition 3.2]{BGM} for the well-posedness and renormalization properties and in \cite[Theorem VI.1.9]{BF} for the stability property.

\subsection{Proof of Proposition \ref{Prop:formulerep}}\label{Annexe:proofrep}
This Section aims at giving a proof for the representation formula (\ref{formule-rep}), which holds for the weak solution to the Vlasov equation. We use the notations and definitions of Section \ref{Section:particle:traj}.

For $(x,v) \in \Omega \times \R^3$ and for any $t \geq 0$, we first define
\begin{align*}
\tau^{+}(t,x,v)&:=\sup \left\lbrace s \geq t \ \mid \forall \sigma \in [t,s], \  \mathrm{X}(\sigma;t,x,v) \in \Omega \right \rbrace.
\end{align*}
If $t \geq 0$ is fixed, we recall that 
\begin{align*}
\mathcal{O}^t &:= \left\lbrace  (x,v) \in \Omega \times \R^3 \ \mid \ \tau^{-}(t,x,v)<0  \right\rbrace,
\end{align*}
where $\tau^{-}(t,x,v)$ has been defined in (\ref{def:tau-}),
If $t \geq 0$ is fixed, we observe that
\begin{align*}
\mathcal{O}^t &=\underset{\sigma \in [0,t]}{\bigcap} \mathrm{Z}_{t,\sigma}(\Omega \times \R^3), \\ 
\mathrm{Z}_{0,t}(\mathcal{O}^t)&=\left\lbrace  (x,v) \in \Omega \times \R^3 \ \mid \ \tau^{+}(0,x,v) >t  \right\rbrace = \underset{\sigma \in [0,t]}{\bigcap} \mathrm{Z}_{0,\sigma}(\Omega \times \R^3).
\end{align*}
By continuity, we have $\mathrm{X}_{\tau^{+}(0,z),0}(z) \in \partial \Omega$ if $z \in \Omega \times \R^3$ and $\tau^{+}(0,z)< +\infty$. More precisely, we have the following result.
\begin{lemme}
For $z=(x,v) \in \Omega \times \R^3$, if $\tau^{+}(0,z)< +\infty$ then we have
\begin{align}\label{traj:exit}
\mathrm{Z}_{\tau^{+}(0,z),0}(z)=(\mathrm{X}_{\tau^{+}(0,z),0}(z),\mathrm{V}_{\tau^{+}(0,z),0}(z)) \in \Sigma^+ \cup \Sigma^0.
\end{align}
\end{lemme}
\begin{proof}
Let us suppose by contradiction that $\mathrm{V}_{\tau^{+}(0,z),0}(z) \cdot n( \mathrm{X}_{\tau^{+}(0,z),0}(z)) <0$. Note that since $z \in \Omega \times \R^3$, we have $\tau^{+}(0,z)>0$.
Since $\Omega$ is smooth, there exists $r>0$ and $\Psi \in  \mathscr{C}^1(\B(\mathrm{X}_{\tau^{+}(0,z),0}(z),r))$ such that for all $y \in \R^3$,
\begin{align*}
y \in \B\left(\mathrm{X}_{\tau^{+}(0,z),0}(z),r \right) \cap \Omega \Leftrightarrow \Psi(y)>0,\\
y \in \B\left(\mathrm{X}_{\tau^{+}(0,z),0}(z),r \right) \cap \partial \Omega \Leftrightarrow \Psi(y)=0.
\end{align*}
For all $\tau>0$ close to $\tau^+:=\tau^{+}(0,z)$, we have
\begin{align*}
\Psi \left(\mathrm{X}_{\tau,0}(z)\right)&=\Psi\left( \mathrm{X}_{\tau^{+},0}(z) \right) +(\tau-\tau^{+}) \dot{\mathrm{X}}_{\tau^{+},0}(z) \cdot \nabla \Psi \left( \mathrm{X}_{\tau^{+},0}(z) \right)  + o(\tau -\tau^{+})\\[2mm]
&=(\tau^+-\tau) \mathrm{V}_{\tau^{+},0}(z)  \cdot n \left( \mathrm{X}_{\tau^{+},0}(z) \right) \left\vert \nabla \Psi \left( \mathrm{X}_{\tau^{+},0}(z) \right) \right\vert + o(\tau -\tau^{+}).
\end{align*}
Thus, for $\tau>0$ close to $\tau^+$ with $\tau <\tau^+$, we get $\Psi \left(\mathrm{X}_{\tau,0}(z)\right)<0$.
This means there exists $\varepsilon>0$ such that for all $\tau \in (\tau^{+}(0,z)-\varepsilon, \tau^{+}(0,z)) \subset \R^+$, $\mathrm{X}_{\tau,0}(z) \notin \overline{\Omega}$, which is in contradiction with the definition of $\tau^{+}(0,z)$.
\end{proof}

\medskip

We now turn to the proof of Proposition \ref{Prop:formulerep}.
Since the Vlasov equation has a unique solution when $u$ is fixed (see Appendix \ref{DiPernaLions}), we have to check that
\begin{equation*}
   \begin{aligned}
      & \R^+ \times \Omega \times \R^3 &&\longrightarrow \R \\
      & \ (t,\ \ x,\ \ v) && \longmapsto e^{3t} \mathbf{1}_{\mathcal{O}^t}(x,v) f_0(\mathrm{Z}_{0,t}(x,v)),
      \end{aligned}
\end{equation*}
is a weak solution to this equation associated to the velocity field $u$ and starting at $f_0$. Let us fix $T>0$. We thus take $\varphi \in \mathscr{C}_c^{\infty}([0,T] \times \overline{\Omega} \times \R^3)$ such that $\varphi(T)=0$ and vanishing on $[0,T] \times (\Sigma^+ \cup \Sigma_0)$ and we want to show that
\begin{multline*}
\int_0^T \int_{\Omega \times \R^3}  e^{3s} \mathbf{1}_{\mathcal{O}^s}(x,v) f_0(\mathrm{Z}_{0,s}(x,v))\Big[ \partial_t \varphi + v \cdot \nabla_x \varphi + (u-v) \cdot \nabla_v  \varphi \Big](s,x,v)\, \mathrm{d}v  \, \mathrm{d}x  \, \mathrm{d}s \\
=-\int_{\Omega \times \R^3} f_0(x,v)\varphi(0,x,v) \, \mathrm{d}v  \, \mathrm{d}x.
\end{multline*}
First, we use the fact that $$\mathrm{Z}_{0,s}(\mathcal{O}^s)=\left\lbrace  (x,v) \in \Omega \times \R^3 \ \mid \ \tau^{+}(0,x,v) >s  \right\rbrace,$$ to write
\begin{align*}
\int_0^T &\int_{\Omega \times \R^3}  e^{3s} \mathbf{1}_{\mathcal{O}^s}(x,v) f_0(\mathrm{Z}_{0,s}(x,v))\Big[ \partial_t \varphi + v \cdot \nabla_x \varphi + (u-v) \cdot \nabla_v  \varphi \Big](s,x,v)\, \mathrm{d}v  \, \mathrm{d}x  \, \mathrm{d}s \\
&= \int_0^T \int_{\mathcal{O}^s}  e^{3s} f_0(\mathrm{Z}_{0,s}(x,v))\Big[ \partial_t \varphi + v \cdot \nabla_x \varphi + (u-v) \cdot \nabla_v  \varphi \Big](s,x,v)\, \mathrm{d}v  \, \mathrm{d}x  \, \mathrm{d}s \\
&=\int_0^T \int_{\Omega \times \R^3} \mathbf{1}_{\tau^+(0,z)>s}  e^{3s} f_0(z)\Big[ \partial_t \varphi + v \cdot \nabla_x \varphi + (u-v) \cdot \nabla_v  \varphi \Big](s,\mathrm{Z}_{s,0}(z) )\mathrm{J}_s(z) \, \mathrm{d}z \, \mathrm{d}s,
\end{align*}
where we have used the change of variable $z=\mathrm{Z}_{0,s}(x,v)$ and where $\mathrm{J}_s(z)$ stands for the Jacobian of the map $z\mapsto \mathrm{Z}_{s,0}(z)$, whose value is $\mathrm{J}_s(z)=e^{-3s}$ for all $s \geq 0$.
Furthermore, by using
\begin{align*}
\dfrac{\mathrm{d}}{\mathrm{d}s}  \Big[ \varphi(s,\mathrm{Z}_{s,0}(z) ) \Big]=\Big[ \partial_t \varphi + v \cdot \nabla_x \varphi + (u-v) \cdot \nabla_v  \varphi \Big]  (s,\mathrm{Z}_{s,0}(z)),
\end{align*}
we get
\begin{align*}
\int_0^T & \int_{\Omega \times \R^3} \mathbf{1}_{\tau^+(0,z)>s}  e^{3s} f_0(z)\Big[ \partial_t \varphi + v \cdot \nabla_x \varphi + (u-v) \cdot \nabla_v  \varphi \Big](s,\mathrm{Z}_{s,0}(z)) \mathrm{J}_s(z) \, \mathrm{d}z \, \mathrm{d}s \\
&=\int_0^T \int_{\Omega \times \R^3} \mathbf{1}_{\tau^+(0,z)>s}  f_0(z) \dfrac{\mathrm{d}}{\mathrm{d}s} \Big[ \varphi(s,\mathrm{Z}_{s,0}(z) ) \Big] \, \mathrm{d}z \, \mathrm{d}s \\
&=\int_0^T \int_{\Omega \times \R^3} \mathbf{1}_{\tau^+(0,z)>s}   \dfrac{\mathrm{d}}{\mathrm{d}s} \Big[ f_0(z) \varphi(s,\mathrm{Z}_{s,0}(z) ) \Big] \, \mathrm{d}z \, \mathrm{d}s \\
&=\int_{\Omega \times \R^3}  \left\lbrace\int_0^{\min(T,\tau^+(0,z)) } \dfrac{\mathrm{d}}{\mathrm{d}s} \Big[ f_0(z) \varphi(s,\mathrm{Z}_{s,0}(z) ) \Big] \, \mathrm{d}s \right\rbrace \, \mathrm{d}z,
\end{align*}
where we have used Fubini's Theorem. We then write
\begin{align*}
\int_{\Omega \times \R^3}  \left\lbrace\int_0^{\min(T,\tau^+(0,z)) } \dfrac{\mathrm{d}}{\mathrm{d}s} \Big[ f_0(z) \varphi(s,\mathrm{Z}_{s,0}(z) ) \Big] \, \mathrm{d}s \right\rbrace \, \mathrm{d}z = (\mathrm{I})+(\mathrm{II}),
\end{align*}
where
\begin{align*}
(\mathrm{I})&:= \int_{\Omega \times \R^3} \mathbf{1}_{\tau^+(0,z)>T}  \left\lbrace\int_0^{\min(t,\tau^+(0,z)) } \dfrac{\mathrm{d}}{\mathrm{d}s} \Big[ f_0(z) \varphi(s,\mathrm{Z}_{s,0}(z) ) \Big] \, \mathrm{d}s \right\rbrace \, \mathrm{d}z,\\
(\mathrm{II})&:=\int_{\Omega \times \R^3} \mathbf{1}_{\tau^+(0,z) \leq T}  \left\lbrace\int_0^{\min(t,\tau^+(0,z)) } \dfrac{\mathrm{d}}{\mathrm{d}s} \Big[ f_0(z) \varphi(s,\mathrm{Z}_{s,0}(z) ) \Big] \, \mathrm{d}s \right\rbrace \, \mathrm{d}z.
\end{align*}

\medskip
First, we have
\begin{align*} (\mathrm{I})&=\int_{\Omega \times \R^3} \mathbf{1}_{\tau^+(0,z)>T} f_0(z) \varphi(T,\mathrm{Z}_{T,0}(z)) \, \mathrm{d}z-\int_{\Omega \times \R^3} \mathbf{1}_{\tau^+(0,z)>T} f_0(z) \varphi(0,z) \, \mathrm{d}z .
\end{align*}
Since $\varphi(T)=0$, the first integral vanishes and we obtain
\begin{align*} (\mathrm{I})=  -\int_{\Omega \times \R^3} \mathbf{1}_{\tau^+(0,z) > T} f_0(z) \varphi(0,z) \, \mathrm{d}z.
\end{align*}
For the second term, we have
\begin{align*} (\mathrm{II})=\int_{\Omega \times \R^3} \mathbf{1}_{\tau^+(0,z) \leq T} f_0(z) \varphi(\tau^+(0,z),\mathrm{Z}_{\tau^+(0,z),0}(z)) \, \mathrm{d}z -\int_{\Omega \times \R^3} \mathbf{1}_{\tau^+(0,z) \leq T} f_0(z) \varphi(0,z) \, \mathrm{d}z.
\end{align*}
Thanks to (\ref{traj:exit}) and to the fact that $\varphi$ vanishes on $\Sigma^+ \cup \Sigma^0$, we see that the first integral is actually $0$.

\medskip

Eventually, gathering all the previous pieces, we end up with
\begin{align*}
\mathrm{(I})+(\mathrm{II}) 
= -\int_{\Omega \times \R^3} f_0(x,v)\varphi(0,x,v) \, \mathrm{d}v  \, \mathrm{d}x.
\end{align*}
The proof of Proposition (\ref{Prop:formulerep}) is finally complete.

\subsection{The Wasserstein distance}\label{Appendix:Wasserstein}
In this section, $X$ stands for a separable and complete subset of $\R^d$ or $\R^d \times \R^d$ (in the previous sections, we used $X=\overline{\Omega}$ or $X=\overline{\Omega} \times \R^d$ where $\Omega$ is an open subset of $\R^d$).

We recall the definition of the 1-Wasserstein distance on $X$ and the useful and classical Monge-Kantorovich formula
(see \cite[Section 11.8]{Dudley} for instance).

\begin{definition}
For all $m>0$, we define $\mathcal{M}_{1,m}(X)$ the set of  positive measures $\mu$ on $X$ such that
\begin{align*}
\int_X \vert x \vert \,\mathrm{d}\mu(x) < \infty, \ \ \ \  \mu(X)=m.
\end{align*}
\end{definition}

\begin{definition}
For all $m>0$, if $\mu$ et $\nu$ are two measures belonging to $\mathcal{M}_{1,m}(X)$, we define the Wasserstein distance $\mathrm{W}_1(\mu,\nu)$ as the quantity
\begin{align*}
\mathrm{W}_1(\mu,\nu) := \underset{\gamma \in \Pi(\mu,\nu)}{\inf}\int_{X \times X} \vert x-x'\vert \, \mathrm{d}\gamma(x,x'),
\end{align*}
where $\Pi(\mu,\nu)$ stands for the set of positive measures on $X \times X$ whose first marginal is $\mu$ and second marginal $\nu$.
\end{definition}

\begin{proposition}
Fix $m>0$. Given $(\mu_n)_n\in\mathcal{M}_{1,m}(X)^\N$ and $\mu\in\mathcal{M}_{1,m}(X)$, the two following statements are equivalent
\begin{itemize}
\item[(i)] For all $f\in\mathscr{C}_b(X)$, 
\begin{align*}
\int_X (f(z)+|z|)\,\mathrm{d}\mu_n(z) \overset{n\rightarrow + \infty}{\longrightarrow} \int_X (f(z)+|z|)\,\mathrm{d}\mu(z).
\end{align*}
\item[(ii)] $(\W_1(\mu_n,\mu))_n \overset{n\rightarrow + \infty}{\longrightarrow} 0$.
\end{itemize}
\end{proposition}

\begin{theoreme}[Duality formula of Monge-Kantorovich]\label{MongeKanto}
If $\mu$ et $\nu$ are two measures belonging to $\mathcal{M}_{1,m}(X)$, we have the following formula
\begin{align*}
\mathrm{W}_1(\mu,\nu)=\sup \enstq{ \left\vert \int_X \psi(x)\mathrm{d}\mu(x) - \int_X \psi(x) \mathrm{d}\nu(x) \right\vert} {\psi \in \mathrm{Lip}(X), \ \ \Vert \nabla \psi \Vert_{\infty} \leq 1}.
\end{align*}
\end{theoreme}

\subsection{Extension of Lipschitz functions vanishing at the boundary}\label{Appendix:ExtensionLip}

\begin{theoreme}
Let $\Omega$ be a smooth bounded and open subset of $\R^d$. Take $u \in \W^{1,\infty}(\overline{\Omega})$ such that $u_{\mid \partial \Omega}=0$. If we set \begin{equation*}
\forall x \in \R^d, \ \ (Pu)(x):=\left\{
      \begin{aligned}
        & u(x), \ &&\text{if} \ x \in \Omega, \\
	& 0\ &&\text{if} \ x \in \R^3 \setminus \Omega,
      \end{aligned}
    \right.
\end{equation*} 
then the function $Pu\in \W^{1,\infty}(\R^d)$ and we have
\begin{align*}
\Vert Pu\Vert_{\Ld^{\infty}(\R^d)}=\Vert u \Vert_{\Ld^{\infty}(\Omega)}, 
\ \ \Vert \nabla (Pu) \Vert_{\Ld^{\infty}(\R^3)} \leq \Vert \nabla u \Vert_{\Ld^{\infty}(\Omega)}.
\end{align*}
\end{theoreme}
\begin{proof}

Let us show that the extension $Pu$ is Lipschitz on $\R^d$. To do so, take $x,y \in \R^d$.
\begin{itemize}
\item If $x \in \overline{\Omega}$ and $y \in \overline{\Omega}$, then we have $$\vert (Pu)(x)-(Pu)(y) \vert =\vert u(x)-u(y) \vert \leq \Vert \nabla u \Vert_{\Ld^{\infty}(\Omega)}\vert x -y \vert.$$
\item If $x \notin \overline{\Omega}$ and $y \notin \overline{\Omega}$, then $\vert (Pu)(x)-(Pu)(y) \vert =0 \leq\Vert \nabla u \Vert_{\Ld^{\infty}(\Omega)} \vert x -y \vert$.
\item If $x \in \overline{\Omega}$ and $y \notin \overline{\Omega}$, then we use the connectedness of the segment $[x,y]$ to find a point $z \in [x,y] \cap \partial \Omega$. We thus write 
$$\vert (Pu)(x)-(Pu)(y) \vert =\vert u(x)\vert=\vert u(x)-u(z) \vert \leq \Vert \nabla u \Vert_{\Ld^{\infty}(\Omega)} \vert x -z \vert \leq \Vert \nabla u \Vert_{\Ld^{\infty}(\Omega)} \vert x -y \vert.$$
\end{itemize}
Here, we have used the fact that the Lipschitz semi-norm of $u$ on $\Omega$ is smaller than $\Vert \nabla u \Vert_{\Ld^{\infty}(\Omega)}$ (see e.g: \cite[Proposition III.2.9]{BF}). All in all, we end up with $$\vert (Pu)(x)-(Pu)(y) \vert \leq \Vert \nabla u \Vert_{\Ld^{\infty}(\Omega)} \vert x -y \vert,$$ so that $Pu\in \W^{1,\infty}(\R^3)$ with the desired inequalities, because of the identification of $\W^{1,\infty}(\R^3)$ with the space of bounded Lipschitz functions on $\R^3$.
\end{proof}

\subsection{A variant of Gronwall lemma}\label{Annexe:gronwall}
We state a Gronwall lemma under integral form in which exponential decay is obtained.
\begin{lemme}\label{variante:gronwall}
Let $g : \R^{+} \rightarrow \R^{+}$ be an integrable and decaying function such that $g(0)$ is well defined and which satisfies, for some $\lambda >0$, for almost every $t \geq 0$
\begin{align*}
\lambda \int_t ^{ + \infty} g(\tau)  \, \mathrm{d}\tau \leq g(t).
\end{align*}
Then, for almost every $t \geq 0$, we have
\begin{align*}
g(t) \lesssim_{\lambda} g(0) e^{-\lambda t},
\end{align*}
where $\lesssim_{\lambda}$ refers to a constant only depending on $\lambda$.
\end{lemme}
A proof can be found in \cite[Appendix]{HKMM}.

\subsection{Agmon inequality on a bounded domain of $\R^3$}\label{Agmon-ineq}

\begin{proposition}
Let $\Omega$ be a smooth bounded domain of $\R^3$. For all $u \in  \mathrm{H}^2(\Omega)$, we have the inequality
\begin{align*}
\|u\|_{\Ld^\infty(\Omega)}\lesssim \|u\|_{\H^1(\Omega)}^{1/2} \|u\|_{\H^2(\Omega)}^{1/2},
\end{align*}
where $\lesssim$ only depends on $\Omega$.
\end{proposition}
We refer to \cite[Lemma 4.10]{ConstFoias} for a proof.
\medskip

\subsection{Gagliardo-Nirenberg-Sobolev inequality on a bounded domain}
\begin{theoreme}\label{gagliardo-nirenberg}
Let $\Omega$ be a smooth bounded domain of $\R^d$. Let $1 \leq p,q,r \leq \infty$ and $m \in \N$. Suppose $j \in \N$ and $\alpha \in [0,1]$ satisfy the relations
\begin{align*}
&\dfrac{1}{p}=\dfrac{j}{d}+\left( \dfrac{1}{r}-\dfrac{m}{d} \right)\alpha+\dfrac{1-\alpha}{q},\\
&\dfrac{j}{m} \leq \alpha \leq 1,
\end{align*}
with the exception $\alpha<1$ if $m-j-d/r \in \N$. 

Then for all $g \in \Ld^q(\Omega)$, if $\mathrm{D}^m g \in \Ld^r(\Omega)$, we have $\mathrm{D}^j g \in \Ld^p(\Omega)$ with the estimate 
$$ \Vert \mathrm{D}^j g  \Vert_{\Ld^p(\Omega)} \lesssim \Vert \mathrm{D}^m g \Vert_{\Ld^r(\Omega)} ^{\alpha} \Vert g \Vert_{\Ld^q(\Omega)}^{1-\alpha} +\Vert g \Vert_{\Ld^s(\Omega)},$$
where $1 \leq s \leq \max \{ q,r \}$ and where $\lesssim$ only depends on $\Omega$. 

Moreover, if $g$ has a vanishing trace at $\partial \Omega$, we can drop the last term $\Vert g \Vert_{\Ld^s(\Omega)}$ in the r.h.s of the inequality.
\end{theoreme}
This result can be found in \cite[Thm 1.5.2]{CherMil}.

\subsection{Maximal $\Ld^p \Ld^q$ regularity for the Stokes system on a bounded domain}\label{AnnexeMaxregStokes}

Let $\Omega$ be a smooth bounded domain of $\R^d$ and $1 < q <\infty$. Each vector field $u \in \Ld^q(\Omega)$ is uniquely decomposed as 
\begin{align*}
& u=\widetilde{u} + \nabla p, \\
& \widetilde{u}  \in \Ld^q_{\mathrm{div}}(\Omega), \ p \in \Ld^q(\Omega), \ \nabla p \in \Ld^q(\Omega),
\end{align*}
where $\Ld^q_{\mathrm{div}}(\Omega)$ stands for the closure in $\Ld^q(\Omega)$ of $\mathscr{D}_{\mathrm{div}}(\Omega)$.
In this so called Helmoltz decomposition, we recall that the projection $\mathbb{P}_q : u \mapsto \widetilde{u} $ is continuous from $\Ld^q(\Omega)$ to $\Ld^q_{\mathrm{div}}(\Omega)$.

For $1<q<\infty$, we consider the following Stokes operator:
\begin{align*}
A_q := -\mathbb{P}_q \Delta u, \ \  D(A_q):=\Ld^q_{\mathrm{div}}(\Omega) \cap \W^{1,q}_0(\Omega) \cap \W^{2,q}(\Omega).
\end{align*}
We also set
\begin{align*}
\mathrm{D}_q^{1-\frac{1}{s},s}(\Omega):=\left( D(A_q),\Ld^q_{\mathrm{div}}(\Omega) \right)_{1/s,s},
\end{align*}
where $( \ , \ )_{1/s,s}$ refers to the real interpolation space of exponents $(1/s,s)$.

\begin{theoreme}
Consider $0 < T \leq \infty$ and $1 <q,s <\infty$. Then, for every $u_0 \in D_q^{1-\frac{1}{s},s}(\Omega)$ which is divergence free and $f \in \Ld^s(0,T;\Ld^q_{\mathrm{div}}(\Omega))$, there exists a unique solution $u$ of the Stokes system
\begin{align*}
\partial_t u +A_q u &=f, \\ 
u(0,x)&=u_0(x),
\end{align*}
satisfying 
\begin{align*}
& u \in \Ld^s(0,T';D(A_q)) \ \text{for all finite } T' \leq T, \\
& \partial_t u \in \Ld^s(0,T;\Ld^q(\Omega)),
\end{align*}
and
\begin{align*}
\Vert \partial_t u \Vert_{\Ld^s(0,T;\Ld^q(\Omega))}  + \Vert \mathrm{D}^2 u\Vert_{\Ld^s(0,T;\Ld^q(\Omega))}\leq C\left(\Vert u_0 \Vert_{D_q^{1-\frac{1}{s},s}(\Omega)} + \Vert f \Vert_{\Ld^s(0,T;\Ld^q(\Omega))} \right),
\end{align*}
where $C=C(q,s,\Omega)$.

Furthermore, if $u_0  \in \H^1_{\mathrm{div}}(\Omega)$ and if $s \in (1,2)$, the statement holds and we can replace $\Vert u_0 \Vert_{D_q^{1-1/s,s}(\Omega)}$ by $\Vert u_0 \Vert_{\H^1_{0}(\Omega)}$ in the right hand side of the previous inequality.
\end{theoreme}

%
%

\medskip

A proof of this result and further references on the theory can be found in \cite{Giga}. The last statement about the fact that $\H^1_{\mathrm{div}}(\Omega)=D(A^{\frac{1}{2}}) \hookrightarrow D_q^{1-1/s,s}(\Omega)$ for $s \in (1,2)$ comes from \cite[Remark 2.5]{Giga}.
\end{appendix}

\subsection{Parabolic regularization for the Navier-Stokes system with a source term on a bounded domain}\label{AnnexeParabNS}
\begin{theoreme}\label{RegParabNS}
Let $\Omega$ be a regular bounded domain of $\R^3$. There exist two universal constants $C_1, C_2 >0$ such that the following holds. 
Consider $u_0 \in \H^1_{\mathrm{div}}(\Omega)$ and $F \in \Ld^2_{\mathrm{loc}}(\R^+;\Ld^2(\Omega))$ and $T>0$ such that
\begin{align}\label{Hyp:AnnexNSreg}
\Vert \nabla u_0 \Vert_{\Ld^2(\Omega)}^2 + C_1 \int_0^T \Vert F(s) \Vert_{\Ld^2(\Omega)}^2 \, \mathrm{d}s \leq \dfrac{1}{\sqrt{8C_1 C_2}}.
\end{align}
Then, there exists on $[0,T]$ a unique Leray solution to the Navier-Stokes equations with initial data $u_0$ and source $F$. This solution $u$ belongs to $\Ld^{\infty}(0,T;\H^1_{\mathrm{div}}(\Omega)) \cap \Ld^{2}(0,T;\H^2(\Omega))$ and satisfies for a.e. $0 \leq t \leq T$
\begin{align}\label{Ineg:AnnexNSreg}
\Vert \nabla u(t) \Vert_{\Ld^2(\Omega)}^2 + \dfrac{1}{2}\int_0^t \Vert A_2 u(s) \Vert_{\Ld^2(\Omega)}^2 \, \mathrm{d}s \leq \Vert \nabla u_0 \Vert_{\Ld^2(\Omega)}^2 + C_1 \int_0^t \Vert F(s) \Vert_{\Ld^2(\Omega)}^2 \, \mathrm{d}s,
\end{align}
where $A_2$ stand for the Stokes operator on $\Ld^2(\Omega)$.
\end{theoreme}
\begin{proof}
First, if such a Leray solution to the Navier-Stokes equations exists, we have directly $u \in \Ld^{4}(0,T;\H^1_{\mathrm{div}}(\Omega))$, which corresponds to a classical tridimensional case of weak-strong uniqueness (see for instance \cite[Theorem 3.3]{NS-chemin}). Hence, it remains to prove that such a solution does exist on $[0,T]$. We proceed in the following way.

We rely on an approximation procedure by regularising the data $F$ and $u_0$ and by considering a standard Galerkine approximation $(u_N)_N$ of the corresponding Navier-Stokes system. The classical idea is to obtain the desired parabolic estimations on the sequence $(u_N)_N$ on $[0,T]$. Combining these new estimates with the classical energy estimates for the Leray solution of the Navier-Stokes system on $[0,T]$, we can use a compactness argument to produce a solution with the $\Ld^{\infty}(0,T;\H^1_{\mathrm{div}}(\Omega)) \cap \Ld^{2}(0,T;\H^2(\Omega))$ regularity and which satisfies the estimate (\ref{Ineg:AnnexNSreg}) on $[0,T]$.

Thus, in the following, it is sufficient to work with a smooth solution $u$ and with smooth data.
We first apply the Leray projection $\mathbb{P}$ to the Navier-Stokes equations and multiply this by $Au$ to obtain, as in the proof of Proposition \ref{controlLinfiniLOC}
\begin{align}\label{enstrophyNS}
\dfrac{\mathrm{d}}{\mathrm{d}t}\Vert \nabla u \Vert_{\Ld^2(\Omega)}^2 + \Vert  Au \Vert_{\Ld^2(\Omega)}^2 \leq  C_1\Vert F \Vert^2_{\Ld^2(\Omega)} + C_1\Vert \nabla u \Vert_{\Ld^{2}(\Omega)}^{6},
\end{align}
where $C_1>0$. Moreover, by Poincaré inequality, there exists an another universal constant $C_2>0$ such that $\Vert \nabla u \Vert_{\Ld^2(\Omega)}^2 \leq C_2 \Vert Au \Vert_{\Ld^2(\Omega)}^2$. So we rewrite the inequality (\ref{enstrophyNS}) as
\begin{align*}
\dfrac{\mathrm{d}}{\mathrm{d}t}\Vert \nabla u \Vert_{\Ld^2(\Omega)}^2 + \dfrac{1}{2}\Vert  Au \Vert_{\Ld^2(\Omega)}^2 \leq  C_1\Vert F \Vert^2_{\Ld^2(\Omega)} + C_1\Vert \nabla u \Vert_{\Ld^{2}(\Omega)}^{6}-\dfrac{1}{2C_2} \Vert \nabla u \Vert_{\Ld^2(\Omega)}^2,
\end{align*}
on $[0,T]$. Now, if we set
\begin{align}\label{def:xyz}
x(t):= \Vert \nabla u(t) \Vert_{\Ld^2(\Omega)}^2,  \ \ \ \ y(t):=\Vert  Au (t)\Vert_{\Ld^2(\Omega)}^2, \ \ \ \ z(t):=\Vert F(t) \Vert^2_{\Ld^2(\Omega)},
\end{align}
and if we integrate the previous inequality between $0$ and $t$, we infer that for almost every $0 \leq t \leq T$,
\begin{align}\label{Inegdiff}
x(t)+ \dfrac{1}{2}\int_0^t y(s) \, \mathrm{d}s \leq x(0)+\int_0^t C_1 z(s) \, \mathrm{d}s +\int_0^t C_1 x(s) \left(x(s)^2-\dfrac{1}{2C_1 C_2} \right) \, \mathrm{d}s.
\end{align}
To obtain (\ref{Ineg:AnnexNSreg}), it is sufficient to prove that $x(t) \leq 1/\sqrt{2C_1 C_2}$ for $t \in [0,T)$.
To do so, we remark that the assumption (\ref{Hyp:AnnexNSreg}) precisely corresponds to
\begin{align*}
 x(0)+\int_0^t C_1 z(s) \, \mathrm{d}s \leq \dfrac{1}{2\sqrt{2C_1 C_2}},
\end{align*}
and in particular,  $x(0) \leq 1/2\sqrt{2C_1 C_2}$. So, by continuity of $t \mapsto x(t)$, there exists a maximal time $T_1 \in ]0,T]$ such that $x(t) \leq 1/\sqrt{2C_1 C_2}$ on $[0,T_1[$. If $T=T_1$, there is nothing to do. So, we argue by contradiction by assuming that $T_1<T$. We thus have for all $0 \leq s \leq T_1$
$$ x(s)^2-\dfrac{1}{2C_1 C_2} \leq 0,$$ and the inequality (\ref{Inegdiff}) with $t=T_1$ turns into
\begin{align*}
x(T_1)+ \dfrac{1}{2}\int_0^{T_1} y(s) \, \mathrm{d}s &\leq x(0)+\int_0^{T_1} C_1 z(s) \, \mathrm{d}s, 
 \leq  x(0)+\int_0^{T} C_1 z(s) \, \mathrm{d}s, 
\end{align*}
which implies in particular
\begin{align*}
x(T_1) \leq \dfrac{1}{2 \sqrt{2C_1 C_2}} <\dfrac{1}{\sqrt{2C_1 C_2}}.
\end{align*}
Again by continuity of $t \mapsto x(t)$, there exists $\varepsilon >0$ such that for all $t \in [T_1, T_1+\varepsilon]$, $x(t) \leq  1/2\sqrt{C_1 C_2}$, which is a contradiction with the definition of $T_1$. It thus implies that for all $t \in [0,T)$, 
$x(t) \leq 1/\sqrt{2C_1 C_2}$ and we finally end up with the following inequality on $[0,T]$
\begin{align*}
x(t)+ \dfrac{1}{2}\int_0^t y(s) \, \mathrm{d}s \leq x(0)+\int_0^t C_1 z(s) \, \mathrm{d}s \leq \dfrac{1}{2\sqrt{2C_1 C_2}}.
\end{align*}
By a view of the definition (\ref{def:xyz}), this finally proves the inequality (\ref{Ineg:AnnexNSreg}).
\end{proof}

\noindent {\bf Acknowledgements.} Partial support by the grant ANR-19-CE40-0004 is acknowledged.

\bibliographystyle{abbrv}


\bibliography{biblio}

\end{document}